\documentclass[a4paper, 10.5pt]{article}
\usepackage{amsmath}
\usepackage{amssymb,esint}
\usepackage{amscd}
\usepackage{xspace}
\usepackage{fancyhdr}
\usepackage{color}
\newtheorem{theorem}{Theorem}[section]

\newtheorem{corollary}[theorem]{Corollary}

\newtheorem{definition}[theorem]{Definition}

\newtheorem{lemma}[theorem]{Lemma}

\newtheorem{proposition}[theorem]{Proposition}

{\catcode`\@=11\global\let\AddToReset=\@addtoreset
\AddToReset{equation}{section}

\AddToReset{theorem}{section}

\newcommand{\be}{\begin{equation}}
\newcommand{\bel}[1]{\begin{equation}\label{#1}}
\newcommand{\ee}{\end{equation}}
\newcommand{\barr}{\begin{eqnarray}}
\newcommand{\earr}{\end{eqnarray}}
\newcommand{\bars}{\begin{eqnarray*}}
\newcommand{\ears}{\end{eqnarray*}}

\newtheorem{subn}{\name}

\newcommand{\bsn}[1]{\def\name{#1}\begin{subn}}
\newcommand{\esn}{\end{subn}}
\newtheorem{sub}{\name}[section]
\newcommand{\bs}{\begin{sub}}
\newcommand{\es}{\end{sub}}
\newcommand{\bsl}[1]{\begin{sub}\label{#1}}
\newcommand{\bth}[1]{\def\name{Theorem}
\begin{sub}\label{t:#1}}
\newcommand{\blemma}[1]{\def\name{Lemma}
\begin{sub}\label{l:#1}}
\newcommand{\bcor}[1]{\def\name{Corollary}
\begin{sub}\label{c:#1}}
\newcommand{\bdef}[1]{\def\name{Definition}
\begin{sub}\label{d:#1}}
\newcommand{\bprop}[1]{\def\name{Proposition}
\begin{sub}\label{p:#1}}

\newcommand{\BA}{\begin{array}}
\newcommand{\EA}{\end{array}}
\newcommand{\BAN}{\renewcommand{\arraystretch}{1.2}
\setlength{\arraycolsep}{2pt}\begin{array}}
\newcommand{\BAV}[2]{\renewcommand{\arraystretch}{#1}
\setlength{\arraycolsep}{#2}\begin{array}}
\newcommand{\BSA}{\begin{subarray}}
\newcommand{\ESA}{\end{subarray}}
\newcommand{\BAL}{\begin{aligned}}
\newcommand{\EAL}{\end{aligned}}
\newcommand{\BALG}{\begin{alignat}}
\newcommand{\EALG}{\end{alignat}}
\newcommand{\BALGN}{\begin{alignat*}}
\newcommand{\EALGN}{\end{alignat*}}
\newcommand{\note}[1]{\textit{#1.}\hspace{2mm}}
\newcommand{\Proof}{\note{Proof}}
\newcommand{\qeda}{\hspace{10mm}\hfill $\square$}

\newcommand{\abs}[1]{\left |#1\right |}
\newcommand{\norm}[1]{\left \|#1\right \|}

\def\angb<#1>{\langle #1 \rangle}

\newcommand{\opname}[1]{\mbox{\rm #1}\,}
\newcommand{\supp}{\opname{supp}}
\newcommand{\dist}{\opname{dist}}
\newcommand{\myfrac}[2]{{\displaystyle \frac{#1}{#2} }}
\newcommand{\myint}[2]{{\displaystyle \int_{#1}^{#2}}}

\newcommand{\prt}{\partial}

\newcommand{\ti}{\times}

\newcommand{\nind}{\noindent}

\def\ga{\alpha}     \def\gb{\beta}       
             \def\ge{\epsilon}

            \def\gl{\lambda}
\def\gm{\mu}        \def\gn{\nu}

\def\Gl{\Lambda}          
\def\Gw{\Omega}              

   \def\CM^+{{\mathcal M}}

      \def\CF{{\mathcal F}}
      
   \def\CK{{\mathcal K}}

   \def\BBN {\mathbb N}    
   \def\BBR {\mathbb R}

\date{}
\begin{document}
\title{Quasilinear elliptic  equations  with a source reaction term involving the function and its gradient and  measure data
}
\author{
	{\bf Marie-Fran\c{c}oise Bidaut-V\'eron\thanks{E-mail address: veronmf@univ-tours.fr, Laboratoire de Math\'ematiques et Physique Th\'eorique,Universit\'e Fran\c{c}ois Rabelais,  Tours,  France}}\\[0.5mm]
	{\bf Quoc-Hung Nguyen\thanks{E-mail address: qhnguyen@shanghaitech.edu.cn, ShanghaiTech University, 393 Middle Huaxia Road, Pudong,
			Shanghai, 201210, China.}}\\[0.5mm]
	{\bf Laurent V\'eron\thanks{ E-mail address: Laurent.Veron@lmpt.univ-tours.fr, Laboratoire de Math\'ematiques et Physique Th\'eorique,Universit\'e Fran\c{c}ois Rabelais,  Tours,  France}}\\[2mm]}
\maketitle

{\abstract We study the equation $-\text{div}(A(x,\nabla u))=|u|^{q_1-1}u|\nabla u|^{q_2}+\mu$
where $A(x,\nabla u)\sim |\nabla u|^{p-2}\nabla u$ in some suitable sense, $\gm$ is a measure and  $q_1$, $q_2$ are nonnegative real numbers and satisfy $q_1+q_2>p-1$. We give sufficient conditions for existence of solutions expressed in terms of the Wolff potential or the Riesz potentials of the measure. Finally we connect the potential estimates on the measure with Lipchitz estimates with respect to some Bessel or Riesz capacity.}\smallskip

\nind{\small key-words: Quasilinear equations; Wolff and Riesz potentials; Hardy-Littlewood maximal function; renormalized solutions; Bessel and Riesz capacities.}\smallskip

\nind {\small 2010 Mathematics Subject Classication: 31C15, 35J62, 35J92, 35R06, 45G15.}
\tableofcontents
\section{Introduction and main results}
This article is devoted to the study of existence of solutions of some second order quasilinear equations with measure data with a source-reaction term involving the function and its gradient. First we consider the problem with a  Radon measure $\mu$ in $\mathbb{R}^N$ in the whole space
\begin{equation}\label{In-1}
-\text{div}(A(x,\nabla u))=\abs u^{q_1-1}u\abs{\nabla u}^{q_2}+\mu~~\text{ in }~\mathbb{R}^N.
\end{equation}
In this setting, $(x,\xi)\mapsto A(x,\xi)$ from $\mathbb{R}^N\times\mathbb{R}^N$ to $\mathbb{R}^N$ is a Carath\'eodory vector field  satisfying for almost all $x\in\BBR^N$ the growth and ellipticity conditions
\bel{In-3}\BA {llll}
(i)&\qquad\qquad|\mathcal{A}(x,\xi)|+|\xi||\nabla_\xi \mathcal{A}(x,\xi)||\leq \Lambda_1|\xi|^{p-1}\quad\text{for all }\;\xi\in\BBR^N,\\[1mm]

(ii)&\qquad\qquad \langle\mathcal{A}(x,\xi)-\mathcal{A}(x,\eta),\xi-\eta \rangle\geq   \Lambda_2 (|\xi|^2+|\eta|^2)^{\frac{p-2}{2}}|\xi-\eta|^2 ~\text{for all}\;\xi,\eta\in\BBR^N,\\[1mm]
(iii)&\qquad\qquad |A(x,\xi)-A(y,\xi)|\leq \Lambda_1 |x-y|^{\alpha_0}|\xi|^{p-1}~\text{for all}\;\xi\in\BBR^N,\\[1mm]

(iv) &\qquad\qquad A(x,\gl \xi)=|\gl|^{p-2}\gl A(x, \xi)\quad\text{for all }(\gl,\xi)\;\in\BBR\ti\BBR^N,
\EA\ee
where $\Gl_1\geq\Gl_2>0$ are constants and  $\frac{3N-2}{2N-1}<p<N$, and where $q_1,q_2>0$ satisfy  $q_1+q_2>p-1$, and $\alpha_0\in (0,1)$.   The special case $\mathcal{A}(x, \xi)=|\xi|^{p-2}\xi$ gives 
rise to the standard $p$-Laplacian $\Delta_p u= {\rm div}\,( |\nabla u|^{p-2} \nabla u)$. Note that these conditions imply that  $A(x,0)=0$ for a.e. $x\in \mathbb{R}^N$, and 
\begin{equation*}
\langle \nabla_{\xi} A(x,\xi)\lambda, \eta\rangle \geq 2^{\frac{p-2}{2}}   \Lambda_2 |\xi|^{p-2} |\eta|^{2}
\end{equation*}
for every $(\eta,\xi)\in\mathbb{R}^N \times \mathbb{R}^N\setminus\{(0,0)\}$ and  a.e. $x \in \mathbb{R}^N$.\\

 When $p=2,q_1=q_2=1$, we obtain a toy model of the forced stationary Navier—Stokes equations describing the motion of incompressible fluid in the whole space $\mathbb{R}^N$:
 \begin{equation}\label{In-2*}
 \begin{array}{lll}%
 &-\Delta U+\nabla p=-U.\nabla U+F,\\&
~~~~~~\text{div}(U)=0,
 \end{array}
 \end{equation}
 in $\mathbb{R}^N$,where  $U=(U_1,...,U_N):\mathbb{R}^N\to \mathbb{R}^N$
 is an unknown velocity of the fluid, $P:\mathbb{R}^N\to \mathbb{R}$ is an unknown pressure, and $ F=(F_1,...,F_N):\mathbb{R}^N\to \mathbb{R}^N$ is a given external force. \\
 
We also consider the homogeneous Dirichlet problem with measure data in a bounded domain $\Omega\subset\mathbb{R}^N$
\begin{equation}\label{In-2}
\begin{array}{lll}%
-\text{div}(A(x,\nabla u))=\abs u^{q_1-1}u\abs{\nabla u}^{q_2}+\mu\qquad&\text{in }\Omega,\\
\phantom{-\text{div}(A(x,\nabla ))}u=0&\text{on } \partial\Omega,
\end{array}
 \end{equation}
where, in this setting, $A:(x,\xi)\mapsto A(x,\xi)$ is a Carath\'eodory vector field defined in $\Gw\ti\BBR^N$ satisfying \eqref{In-3} (i)-(iv) in $\Gw\ti\BBR^N$, and $q_1$, $q_2$ are as in the first case; and   $\Omega\subset \mathbb{R}^N$ is a bounded domain with a $C^{1,\beta_0}$ boundary for $\beta_0\in (0,1)$ and such that $\Omega\subset B_{R}(x_0)$ for some $R>0$ and $x_0\in \Omega$. 
The two specific cases,
\bel{In-6}\BA {lll}-\text{div}(A(x,\nabla u))=\abs u^{q_1-1}u+\mu\qquad&\text{in }\Gw\\
\phantom{-\text{div}(A(x,\nabla ))}u=0&\text{in }\prt\Gw,
\EA\ee 
and 
\bel{In-7}\BA {lll}-\text{div}(A(x,\nabla u))=\abs {\nabla u}^{q_2}+\mu\qquad&\text{in }\Gw,\\
\phantom{-\text{div}(A(x,\nabla ))}u=0&\text{in }\prt\Gw,
\EA\ee  
have been studied thoroughly in the last decade. Each of these equations carries a critical exponent $q_j^c$: $q_1^c=\frac{N(p-1)}{N-p}$ for equation (\ref{In-6}) and $q_2^c=\frac{N(p-1)}{N-1}$ for (\ref{In-7}). These critical thresholds mean that if $0<p-1<q_1<q_1^c$ for (\ref{In-6})  and $1-\frac 1N<p-1<q_1<q_1^c$ for (\ref{In-7}) any nonnegative bounded measure is eligible for the respective equation, provided it is small enough. Concerning equation (\ref{In-1}),   
the criticality is expressed by a linear relation $0<q_1(N-p)+q_2(N-1)<N(p-1)$. Then, if  $p>2-\frac 1N$ and $q_1+q_2>p-1$, problem (\ref{In-2}) any Radon measure small enough, see \cite[Chap 6-2]{Ve1} and references therein. The treatment of the supercritical case for equations  (\ref{In-6})
and (\ref{In-7}) have been treated more recently. In these cases not only the measure $\gm$ has to be small enough, but also it cannot be too concentrated with respect to some Bessel capacity, specific to each problem. It is proved in \cite{22PhVe} that if $\gm$ is a nonnegative Radon measure with compact support in $\Gw$, a necessary and sufficient condition for the existence of a renormalized solution to  (\ref{In-6}) is that there exists some $c_1>0$ depending on the structural constants and $\norm\gm_{\mathfrak M}$ such that 
\bel{In-9}
\gm(K)\leq c_1Cap_{\mathbf{G}_p,\frac{q_{_1}}{q_{_1}+1-p}}(K)\quad\text{for all compact set }K\subset\overline\Gw,
\ee 
where  $Cap_{\mathbf{G}_p,\frac{q_{_1}}{q_{_1}+1-p}}$ denotes some Bessel capacity. Concerning  (\ref{In-7}), assuming $\max\{p-1,1\}<q_{_2}$, it is proved in \cite{Ph2,Ph1,QH4,QH6,QH5} that there exists a structural constant $c_2>0$ as above such that if 
\bel{In-10}
|\gm|(K)\leq c_2Cap_{\mathbf{G}_1,\frac{q_{_2}}{q_{_2}+1-p}}(K)\quad\text{for all compact set }K\subset\Gw,
\ee 
there exists a renormalized solution to  (\ref{In-7}) with the property that 
\bel{In-11}
\myint{K}{}\abs{\nabla u}^{q_2}dx\leq c_3Cap_{\mathbf{G}_1,\frac{q_{_2}}{q_{_2}+1-p}}(K)\quad\text{for all compact set }K\subset\overline\Gw,
\ee
for some $c_3>0$. \medskip

The complete expression of these results as well as the ones we will state below necessitates the introduction of several definitions and notations from harmonic analysis such as Wolff potential, Riesz potentials, Bessel spaces and maximal functions. The role of these operators has appeared to be a key-stone for conducting a fine analysis of quasilinear equations with measure data; this is very clearly presented in the introduction of the seminal paper \cite{22PhVe}. 
If $D$ is either a bounded domain or whole $\mathbb{R}^{N}$, we denote by $\mathfrak{M}(D)$ (resp. $\mathfrak{M}_b(D)$) the set of Radon measures (resp. bounded Radon measures) in $D$. Their positive cones are $\mathfrak{M}^+(D)$ and  $\mathfrak{M}_b^+(D)$ respectively. For $R\in (0,\infty]$, we define the $R$- truncated Wolff potential $\mathbf{W}^R_{\alpha,p}$
                $(\alpha\in (0,N/p), p>1)$ and the  $R$-truncated Riesz potential $\mathbf{I}^R_{\beta}$ $(\beta\in (0,N))$ of a measure $\mu\in\mathfrak{M}^+(\mathbb{R}^{N})$ by 
 \begin{equation}\label{In-12}
\mathbf{W}_{\alpha,p}^R[\mu](x)=\int_{0}^{R}\left(\frac{\mu(B_\rho(x))}{\rho^{N-\alpha p}}\right)^{\frac{1}{p-1}}\frac{d\rho}{\rho}~\text{ and }~ \mathbf{I}_{\beta}^R[\mu](x)=\int_{0}^{R}\frac{\mu(B_\rho(x))}{\rho^{N-\beta}}\frac{d\rho}{\rho},
 \end{equation}
for all $x$ in $\mathbb{R}^{N}$. If $R=\infty$, we drop it in the expressions of \eqref{In-12}. We write $  \mathbf{W}_{\alpha,p}^R[f]$,  $\mathbf{I}_{\beta}^R[f]$ in place of  $  \mathbf{W}_{\alpha,p}^R[\mu]$,  $\mathbf{I}_{\beta}^R[\mu]$ whenever $d\mu=fdx$, where $f\in L^1_{loc}(\mathbb{R}^N)$.                                                    
                    
For $\alpha>0$, $p>1$, the $(\mathbf{I}_{\alpha},p)$-capacity, $(\mathbf{G}_{\alpha},p)$-capacity of a Borel set $O\subset \mathbb{R}^{N}$ are defined by
\begin{align*}
& \text{Cap}_{\mathbf{I}_{\alpha},p}(O)=\inf\left\{\int_{\mathbb{R}^{N}}|g|^pdx: g\in L^p_+(\mathbb{R}^{N}), \mathbf{I}_{\alpha}[g]\geq \chi_{O}\right\},\\&
\text{Cap}_{\mathbf{G}_{\alpha},p}(O)=\inf\left\{\int_{\mathbb{R}^{N}}|g|^pdx: g\in L^p_+(\mathbb{R}^{N}), \mathbf{G}_{\alpha}*g\geq \chi_{O}\right\},
\end{align*}
where $\mathbf{G}_\alpha=\CF^{-1}\left((1+\abs\xi^2)^{-\frac\ga 2}\right)$ is the Bessel kernel of order $\alpha$, see \cite{55AH} (and $\CF$ and $\CF^{-1}$ are respectively the Fourier transform and its inverse).\\

The results we prove  consist in obtaining sufficient conditions for the solvability of (\ref{In-1}) or (\ref{In-2}) where $A$ is of the form  (\ref{In-3}) in $\BBR^N$ (or $\Gw$)  expressed in terms of inequalities between Wolff or Riesz potentials of $\gm$. In order to obtain these inequalities we will develop a series of sharp relations between these potentials and will connect them with some specific capacities. We recall that a Radon measure 
$\gm$ in $\BBR^N$ (or $\Gw$) is {\it absolutely continuous} with respect to some capacity Cap in $\BBR^N$ (or $\Gw$) if for a Borel set $E$
 \bel{In-13}
\text{Cap} (E)=0\Longrightarrow \abs\gm(E)=0,
\ee
and it is {\it Lipschitz continuous} (with constant $c>0$) if 
 \bel{In-14}
 \abs\gm(E)\leq c\text{Cap} (E)\quad\text{for all Borel set } E.
\ee
The capacity associated to the Sobolev space $W^{1,p}(\BBR^N)$ is denoted by Cap$_{1,p}$. It coincides with Cap$_{\mathbf{G}_1,p}$, \cite[Th 1.2.3]{55AH}. A measure is called {\it diffuse} if it is absolutely continuous with respect to Cap$_{1,p}$\medskip

Our first result deals with the equation in the whole space,
  \begin{theorem}\label{mainthem1}Let $q_1,q_2>0, q_1+q_2>p-1,0<q_2<\frac{N(p-1)}{N-1}$ and  $\mu\in \mathfrak{M}(\mathbb{R}^N)$. Assume that $A(x,\xi)=A(\xi)$ for any $(x,\xi)\in \mathbb{R}^N\times\mathbb{R}^N$.  If for some 
  	$C>0$  depending on $p$, $N$, $q_j$ and $\Gl_j$ (j=1,2), there holds
  \begin{align}
    |\mu|(K)\leq C\text{Cap}_{\mathbf{I}_{\frac{q_1p+q_2}{q_1+q_2}},\frac{q_1+q_2}{q_1+q_2-p+1}}(K)\quad\text{for all compact }K\subset\mathbb{R}^N,
  \end{align}
	then problem   \eqref{In-1} admits a  distributional solution $u$ which satisfies
  \bel{In-18}
  |u|\leq C_0\mathbf{W}_{1,p}[|\mu|],\quad |\nabla u|\leq C_0\mathbf{W}_{\frac{1}{p},p}[|\mu|],
  \ee
  if $p>2$, and 
  \bel{In-19}
  |u|\leq C_0\left(\mathbf{I}_{p}[|\mu|]\right)^{\frac{1}{p-1}},~~~|\nabla u|\leq C_0\left(\mathbf{I}_{1}[|\mu|]\right)^{\frac{1}{p-1}},
  \ee
  if $\frac{3N-2}{2N-1}<p\leq 2$.  	Moreover, if $\mu\geq 0$, then $u\geq 0$.     
  \end{theorem}
Notice also that if $\mu\geq 0$ the solutions $u$ in Theorem \ref{mainthem1} are  nonnegative p-super-harmonic functions.\\
   
When $\BBR^N$ is replaced by a bounded domain $\Gw$, we have the following general result.   
\begin{theorem}\label{mainthem2}Let $q_1,q_2>0, q_1+q_2>p-1,0<q_2<\frac{N(p-1)}{N-1}$. Let $
\mu\in \mathfrak{M}(\Omega)$ be such that $\dist(\supp(\mu),\partial\Omega)>0$.  If for some 
	$C>0$  depending on $p$, $N$, $q_j$ and $\Gl_j$ (j=1,2), $\Omega$ and $\dist(\supp(\mu),\partial\Omega)$  there holds
	\begin{align}
	|\mu|(K)\leq C\text{Cap}_{\mathbf{G}_{\frac{q_1p+q_2}{q_1+q_2}},\frac{q_1+q_2}{q_1+q_2-p+1}}(K)\quad\text{for all compact }K\subset\mathbb{R}^N,
	\end{align}	
	then  problem \eqref{In-2}  admits a renormalized solution $u$ satisfying 
	 \bel{In-20}
	|u|\leq C_0\mathbf{W}_{1,p}[|\mu|]\quad\text{and }\; |\nabla u|\leq C_0\mathbf{W}_{\frac{1}{p},p}[|\mu|],
	\ee
	if $p>2$, and 
	\bel{In-21}
	|u|\leq C_0\left(\mathbf{I}_{p}[|\mu|]\right)^{\frac{1}{p-1}},~~~|\nabla u|\leq C_0\left(\mathbf{I}_{1}[|\mu|]\right)^{\frac{1}{p-1}},
	\ee
	if $\frac{3N-2}{2N-1}<p\leq 2$.  
	Moreover, if $\mu\geq 0$, then $u\geq 0$.     
\end{theorem}
 
 The key-stone of our method which combines sharp potential estimates and Schauder fixed point theorem is to reduce our problems (\ref{In-1})-(\ref{In-2}) to a system of nonlinear Wolff integral equations in the spirit of the method developed in \cite{VHV2} and \cite{QH-Ve}.  The proof of Theorem \ref{mainthem1} (and similarly for Theorem\ref{mainthem2}) is based upon the existence of a fixed point obtained by Schauder's theorem, of the 
 mapping $S$ which associates to  
 $$v\in E_{\Gl}=\left\{v:|v|\leq \Gl \left({\bf I}_p[|\gm_{n,k}|]\right)^{\frac{1}{p-1}}\text{ s.t. }|\nabla v|\leq \Gl \left({\bf I}_1[|\gm_{n,k}|]\right)^{\frac{1}{p-1}}\text{ in }B_{2k}\right\}
 $$
the solution $u=u_{n,k}=S(v)$, of  
   \bel{I-1-1}\BA {lll}
 -\text{div}\left(A(x,\nabla u)\right)=\chi_{_{B_k}}|v|^{q_1-1}v|\nabla v|^{q_2}+\gm_{n,k}\quad&\text{in }\, B_{2k},\\
 \phantom{-\text{div}\left(A(x,\nabla )\right)}
 u=0&\text{on }\,\prt B_{2k},
 \EA
 \ee
where  $\gm_{n,k}$ is a smooth approximation of $\gm$ with support in $B_{k}:=B_{k}(0)$ and $\Gl>0$, $k$, $n\in \BBN_*$ are parameters. In order to prove that the set 
$E_{\Gl}$ is invariant under $S$ we use a series of a priori estimates dealing with renormalized solutions of
 \bel{I-1-2}\BA {lll}
 -\text{div}\left(A(x,\nabla u)\right)=\gn\quad&\text{in }\,B_{2k},\\
 \phantom{-\text{div}\left(A(x,\nabla )\right)}
 u=0&\text{on }\,\prt B_{2k},
 \EA
 \ee
 where $\gn\in\mathfrak M_b(\Gw)$, and for our purpose $\gn=\chi_{_{B_k(0)}}|v|^{q_1-1}v|\nabla v|^{q_2}+\gm_{n,k}$. Then we use pointwise estimates satisfied by a renormalized solution, e.g. in the case $p\leq 2$, there holds a.e. in $B_{k}$
  \bel{I-1-3}\BA {lll}
|u(x)| \leq C\left({\bf I}_{p}\left[|\chi_{_{B_k}}|v|^{q_1-1}v|\nabla v|^{q_2}+\gm_{n,k}|\right](x)\right)^{\frac{1}{p-1}},
 \EA
 \ee  
 and
   \bel{I-1-4}\BA {lll}\displaystyle
|\nabla u(x)|\leq C\left({\bf I}_{1}\left[|\chi_{_{B_k}}|v|^{q_1-1}v|\nabla v|^{q_2}+\gm_{n,k}|\right](x)\right)^{\frac{1}{p-1}}, 
 \EA
 \ee 
 see Theorem \ref{renorm-th-2} and Corollary \ref{renorm-cor-1}. \\
 
 Using the fact that $v\in E_\Gl$ we derive
   \bel{I-1-5}\BA {lll}
|u(x)| \leq C\left(\Gl^{q_1+q_2}{\bf I}_p\left[\left({\bf I}_p[|\gm_{n,k}|]\right)^{\frac{q_1}{p-1}}\left({\bf I}_1[|\gm_{n,k}|]\right)^{\frac{q_2}{p-1}}\right](x)+{\bf I}_p[|\gm_{n,k}|](x)\right)^{\frac{1}{p-1}},
 \EA
 \ee  
 and
   \bel{I-1-6}\BA {lll}\displaystyle
|\nabla u(x)|\leq C\left(\Gl^{q_1+q_2}{\bf I}_1\left[\left({\bf I}_p[|\gm_{n,k}|]\right)^{\frac{q_1}{p-1}}\left({\bf I}_1[|\gm_{n,k}|]\right)^{\frac{q_2}{p-1}}\right](x)+{\bf I}_1[|\gm_{n,k}|](x)\right)^{\frac{1}{p-1}}.
 \EA
 \ee 
At this point we use the multplicative inequalities concerning the Riesz potential provided the measures satisfies some Lipschitz continuity estimate with respect to  $Cap_{{\bf I}_{\frac{pq_1+\gb pq_2}{q_1+q_2}},\frac{q_1+q_2}{q_1+q_2+1-p}}$:
    \bel{I-1-7}\BA {lll}\displaystyle
 {\bf I}_\ga\left[\left({\bf I}_p[|\gm_{n,k}|]\right)^{\frac{q_1}{p-1}}\left({\bf I}_1[|\gm_{n,k}|]\right)^{\frac{q_2}{p-1}}\right]\leq M {\bf I}_\ga[|\gm_{n,k}|]\quad\text{for }\ga=1\text{ or }p.
\EA
 \ee 
 These multiplicative inequalities are the key of our construction, since they imply that for suitable choice of $\Gl$, $E_\Gl$ is invariant. The compactness of $S$ being easy to prove we derive the existence of a solution to (\ref{I-1-1}).\smallskip
 
 One of the tools is a series of equivalence linking the Lipschitz continuity of $\gm$ with respect to some capacity with integral estimates of the Wolff (or Riesz) potential of the measure and even to a system of nonlinear Wolff integral equations as in \cite{QH-Ve}. In this spirit we prove the following:

 \begin{theorem}\label{th-pot-4-}
Let $1<p<N$, $0<\gb<1$, $q_1,q_2>0$ such that   $q_1+q_2>p-1$ and $q_2<\frac{N(p-1)}{N-\gb p}$. If $\gm$ is a nonnegative measure in $\BBR^N$, the  following statements are equivalent:\smallskip
 
 \nind (a) The inequality
 \bel{In-22}
\gm (K)\leq C_1Cap_{I_{\frac{pq_1+\gb pq_2}{q_1+q_2}},\frac{q_1+q_2}{q_1+q_2+1-p}}(K),
	\ee
 holds for any compact set $K\subset\BBR^N$, for some $C_1>0$.\smallskip
 
 \nind (b) The inequality 
  \bel{In-22-1}
\myint{K}{}\left({\bf W}_{\ga,p}[\gm]\right)^{q_1}\left({\bf W}_{\gb,p}[\gm]\right)^{q_2}dx\leq C_2Cap_{I_{\frac{pq_1+\gb pq_2}{q_1+q_2}},\frac{q_1+q_2}{q_1+q_2+1-p}}(K),
	\ee
 holds for any compact set $K\subset\BBR^N$ and some $C_2>0$.\smallskip
 
  \nind (c) The inequalities 
  \bel{In-22-2}\BA {lll}
{\bf W}_{1,p}\left[\left({\bf W}_{1,p}[\gm]\right)^{q_1}\left({\bf W}_{\frac 1p,p}[\gm]\right)^{q_2}\right]\leq C_3{\bf W}_{1,p}[\gm]\\[1mm]
{\bf W}_{\frac 1p,p}\left[\left({\bf W}_{1,p}[\gm]\right)^{q_1}\left({\bf W}_{\frac 1p,p}[\gm]\right)^{q_2}\right]\leq C_3{\bf W}_{\frac 1p,p}[\gm],
\EA	\ee
in $\mathbb{R}^N$ are verified for some $C_3>0$.
\smallskip

 \nind (d) The system of equations
  \bel{In-23}\BA {lll}
    U={\bf W}_{1,p}[U^{q_1}V^{q_2}]+\ge{\bf W}_{1,p}[\gm]\\[1mm]
V={\bf W}_{\frac 1p,p}[U^{q_1}V^{q_2}]+\ge{\bf W}_{\frac 1p,p}[\gm],
\EA\ee
in $\mathbb{R}^N$ has a solution $U,V\geq 0$ for $\ge>0$ small enough.\end{theorem}

 Actually the full statement is more complete and the above Theorem is a consequence of Theorem \ref{th-pot-4}. Furthermore it has an analogue in $\Gw$, see Theorem \ref{th-pot-5}.

\section{ Estimates on potential}
In the sequel $C$ denotes a generic constant depending essentially on some structural constants (i.e. the ones associated to the operator and reaction term) and the domain, the value of which may change from one occurence to another. Sometimes, in order to avoid confusion, we introduce notations $C_j$, $j=0,1,2...$. We also use the notation $\asymp$ to assert that the two quantities linked by this relation are comparable up to multiplication by constants of the previous type.
The following result is a general version of results of Phuc and Verbitsky \cite[Th 2.3]{22PhVe}. It connects the Lipschitz continuity of a positive measure in $\BBR^N$ with respect to some Riesz capacity to various integral or pointwise estimates of Wolff potentials of this measure. 
\begin{theorem} \label{th-pot-1}Let $1<p<N/\alpha$, $q>p-1$,  $\mu\in \mathfrak{M}^+(\mathbb{R}^N)$. Then, the following statements are equivalent:
\begin{description}
	\item[(a)]  The inequality 
	\begin{align}\label{pot-1}
	\mu(K)\leq C_1\text{Cap}_{\mathbf{I}_{\alpha p},\frac{q}{q-p+1}}(K),
	\end{align}
	holds for any compact set $K\subset\mathbb{R}^N$, for some $C_1>0$.
		\item[(b)]   The inequality 
		\begin{align}\label{pot-1*}
		\int_{K}\left(\mathbf{W}_{\alpha,p}[\mu](y)\right)^qdy\leq C_2\text{Cap}_{\mathbf{I}_{\alpha p},\frac{q}{q-p+1}}(K),
		\end{align}
		holds for any compact set $K\subset\mathbb{R}^N$, for some $C_2>0$.
			\item[(c)]   The inequality 
			\begin{align}\label{pot-2}
			\int_{\mathbb{R}^N}\left(\mathbf{W}_{\alpha,p}[\chi_{_{B_t(x)}}\mu](y)\right)^qdy\leq C_3 \mu(B_t(x)),
			\end{align}
			holds for any $x\in\mathbb{R}^N$ and $t>0$, for some $C_3>0$.
				\item[(d)]   The inequality 
				\begin{align}\label{pot-3*}
				\mathbf{W}_{\alpha,p}\left[\left(\mathbf{W}_{\alpha,p}[\mu]\right)^q\right]\leq C_4 \mathbf{W}_{\alpha,p}[\mu]<\infty, 
				\end{align}
				holds  almost everywhere in $\mathbb{R}^N$, for some $C_4>0$.
\end{description}
    \end{theorem}                                                                                       

\nind\Proof  {\it Step 1}: Proof of (a) $ \Leftrightarrow$ (b). By \cite[Theorem 1.1]{11HoJa}, (see also \cite[Theorem 2.3]{55VHV}) we have
\begin{align}\label{pot-3}
\int_{\mathbb{R}^N}\left(\mathbf{I}_{\alpha p}[\nu](y)\right)^{\frac{q}{p-1}}w(y)dy\asymp \int_{\mathbb{R}^N}\left(\mathbf{W}_{\alpha,p}[\nu](y)\right)^qw(y)dy~~\text{for all}~\nu\in\mathfrak{M}^+(\mathbb{R}^N),
\end{align}
where $w$ belongs to the Muckenhoupt class $\mathbf{A}_\infty$. So, thanks to \cite[Lemma 3.1]{MaVe} we obtain
\begin{align*}\displaystyle
\sup_{K \in \CK(\mathbb{R}^N)} \frac{\myint{K}{}\left(\mathbf{I}_{\alpha p}[\nu](y)\right)^{\frac{q}{p-1}}dy}{\text{Cap}_{\mathbf{I}_{\alpha p},\frac{q}{q-p+1}}(K)}\asymp \sup_{K \in \CK(\mathbb{R}^N)} \frac{\myint{K}{}\left(\mathbf{W}_{\alpha,p}[\nu](y)\right)^qdy}{\text{Cap}_{\mathbf{I}_{\alpha p},\frac{q}{q-p+1}}(K)}\qquad\text{for all}~\nu\in\mathfrak{M}^+(\mathbb{R}^N),
\end{align*} 
where $\CK(\BBR^N)$ denotes the set of compact subsets of $\BBR^N$. Moreover, by \cite[Theorem 2.1]{MaVe},
\begin{align*}
\sup_{K \in \CK(\mathbb{R}^N)} \frac{\nu (K)}{\text{Cap}_{\mathbf{I}_{\alpha p},\frac{q}{q-p+1}}(K)}\asymp\sup_{K \in \CK(\mathbb{R}^N)} \frac{\myint{K}{}\left(\mathbf{I}_{\alpha p}[\nu](y)\right)^{\frac{q}{p-1}}dy}{\text{Cap}_{\mathbf{I}_{\alpha p},\frac{q}{q-p+1}}(K)} \qquad\text{for all}~\nu\in\mathfrak{M}^+(\mathbb{R}^N).
\end{align*} 
From this we infer the equivalence between (a) and (b).\smallskip

\nind{\it Step 2}: Proof of (a) $ \Leftrightarrow  $ (c). By \cite[Theorem 2.1]{MaVe} (a) is equivalent to 
\begin{align*}
\int_{\mathbb{R}^N}\left(\mathbf{I}_{\alpha p}[\chi_{_{B_t(x)}}\mu](y)\right)^{\frac{q}{p-1}}dy\leq C \mu(B_t(x)),
\end{align*}
 for any ball $B_t(x)\subset\mathbb{R}^N$. It is equivalent to (c) because of \eqref{pot-3}.\smallskip

\nind{\it Step 3}: By Proposition \ref{prop-pot-2}, we obtain (c) $ \Rightarrow  $ (d).\smallskip

\nind\textbf{\it Step 4}:  Proof of (d) $ \Rightarrow  $ (b).  Set $d\nu(x)=\left(\mathbf{W}_{\alpha,p}[\mu](x)\right)^qdx$. Clearly, (d) implies 
\begin{align*}
\left(\mathbf{W}_{\alpha,p}[\nu](x)\right)^qdx\leq C d\nu(x).
\end{align*} Let $\mathbf{M}_\nu$ denote the centered Hardy-Littlewood maximal function defined for any $f\in L_{loc}^1(\mathbb{R}^N,d\nu)$  by 
               \[{\mathbf{M}_\nu }f(x) = \sup _{t > 0} \frac{1}{\nu (B_t(x))}\int_{B_t(x)}|f|d\nu  .\]
               If $E\subset\mathbb{R}^N$ is a Borel set, we have
               \begin{align*}
           \int_{\mathbb{R}^N}\left(\mathbf{W}_{\alpha,p}[\chi_{_E}\nu]\right)^{q}dx    \leq  \int_{\mathbb{R}^N}(\mathbf{M}_{\nu}\chi_{_E})^{\frac{q}{p-1}}\left(\mathbf{W}_{\alpha,p}[\nu]\right)^{q}dx\leq C\int_{\mathbb{R}^N}(\mathbf{M}_{\nu}\chi_{_E})^{\frac{q}{p-1}}d\nu.
               \end{align*}
               Since ${\mathbf{M}}_\nu$ is bounded on $L^s(\mathbb{R}^N,d\nu)$, $s>1$, by the Besicovitch's theorem, see e.g. \cite{22Fe}, we deduce that 
\begin{align*}
               \int_{\mathbb{R}^N}\left(\mathbf{W}_{\alpha,p}[\chi_{_E}\nu]\right)^{q}dx\leq C\nu(E),
               \end{align*}
               for any Borel set $E$.
Applying the equivalence of (a) and (c) with $\gm=\gn$, we derive (b). $\phantom{------}$
  \qeda
\medskip

The next result is the analogue of the previous one when the whole space is replaced by a ball. It can be proved in the same way, see also \cite[Proof of Theorem 2.3]{Ph2}.
\begin{theorem}\label{th-pot-2} Let $1<p<N/\alpha$, $q>p-1$, $\omega\in \mathfrak{M}_b^+(B_R(x_0))$ for some $R>0$ and $x_0\in\BBR^N$. Then, the following statements are equivalent:\smallskip

\nind (a) The inequality 
\begin{align}\label{pot-4}
\omega(K)\leq C_1\text{Cap}_{\mathbf{G}_{\alpha p},\frac{q}{q-p+1}}(K),
\end{align}
holds for any compact set $K\subset\mathbb{R}^N$, for some $C_1=C_1(R)>0$.\smallskip

\nind (b) The inequality 
\begin{align}\label{pot-4*}
\int_{K}\left(\mathbf{W}^{4R}_{\alpha,p}[\omega](y)\right)^qdy\leq C_2\text{Cap}_{\mathbf{G}_{\alpha p},\frac{q}{q-p+1}}(K),
\end{align}
holds for any compact set $K\subset\mathbb{R}^N$, for some $C_2=C_2(R)>0$.\smallskip

\nind (c) The inequality 
\begin{align}\label{pot-5}
\int_{\mathbb{R}^N}\left(\mathbf{W}^{4R}_{\alpha,p}[\chi_{_{B_t(x)}}\omega](y)\right)^qdy\leq C_3 \omega(B_t(x)),
\end{align}
holds for any $x\in\mathbb{R}^N$ and $t>0$, for some $C_3=C_3(R)>0$.\smallskip

\nind (d) The inequality 
\begin{align}\label{pot-6}
\mathbf{W}^{4R}_{\alpha,p}\left[\left(\mathbf{W}^{4R}_{\alpha,p}[\omega]\right)^q\right]\leq C_4 \mathbf{W}^{4R}_{\alpha,p}[\omega],
\end{align}
holds almost everywhere in $B_{2R}(x_0)$, for some $C_4=C_4(R)>0$.
\end{theorem}

The following stability result of the Lipschitz continuity of a measure with respect a capacity will be used several times in the sequel since we will approximate the initial data by smooth and truncated ones; its proof is easy, see e.g.  \cite[Lemma 2.7]{Ph1}. 
\begin{proposition}\label{prop-pot-1}
   Let $1<p<N/\alpha$ and $0<\beta<N/p$, $\mu\in \mathfrak{M}^+(\mathbb{R}^N)$, $\omega\in \mathfrak{M}_b^+(B_R(x_0))$ for some $R>0$ and $x_0\in\BBR^N$. Set $d\mu_n(x)=(\varphi_n*\mu)(x) dx$, $d\omega_n(x)=(\varphi_n*\omega)(x) dx$ where $\{\varphi_n\}$ is a sequence of mollifiers. Then,
   \begin{description}
   	\item[(i)]  If  inequality  \eqref{pot-1} in Theorem \ref{th-pot-1} holds with $q>\frac{(p-1)N}{N-\alpha p}$ and constant $C_1$, then
   	\begin{align}\label{pot-7}
   	\mu_n(K)\leq CC_1\text{Cap}_{\mathbf{I}_{\alpha p},\frac{q}{q-p+1}}(K)~~\text{for all }\, K\subset\mathbb{R}^N,~ n\in\mathbb{N}
   	\end{align}
   	for some $C=C(N,\alpha,p,q)>0$.
   	\item[(ii)] If  inequality    \eqref{pot-4} in Theorem  \ref{th-pot-2} holds with $q>p-1$ and constant $C_2$, then
   	\begin{align}\label{pot-8}
   	\omega_n(K)\leq CC_2\text{Cap}_{\mathbf{G}_{\alpha p},\frac{q}{q-p+1}}(K)~~\text{for all }\, K\subset\mathbb{R}^N,~ n\in\mathbb{N}
   	\end{align}
   	for some $C=C(N,\alpha,p,q)>0$.
   \end{description}
   \end{proposition}                                                                                       
   
   The next proposition is crucial as it gives pointwise estimates of interates of Wolff potentials of positive measures and connect them with the capacitary estimates of the Wolff potentials of the same measures.

\begin{proposition}\label{prop-pot-2}
   Let $1<p<N/\alpha$ and $0<\beta<N/p$, $\mu\in \mathfrak{M}^+(\mathbb{R}^N)$, $\omega\in \mathfrak{M}_b^+(B_R(x_0))$ for some $B_R(x_0)\subset\mathbb{R}^N$. Then, \smallskip

\nind (i) The inequality  \eqref{pot-2} in Theorem \ref{th-pot-1} with $q>\frac{(p-1)N}{N-\alpha p}$ implies that
   \begin{align}\label{pot-9}
   \mathbf{W}_{\beta,p}\left[\left(\mathbf{W}_{\alpha,p}[\mu]\right)^q\right]\leq C_1 \mathbf{W}_{\beta,p}[\mu]<\infty,
   \end{align}
   holds  almost everywhere in $\mathbb{R}^N$, for some $C_1>0$.\smallskip

\nind (ii)  The inequality  \eqref{pot-5} in Theorem  \ref{th-pot-2} with $q>p-1$ implies that
      \begin{align}\label{pot-10}
      \mathbf{W}_{\beta,p}^{4R}\left[\left(\mathbf{W}^{4R}_{\alpha,p}[\omega]\right)^q\right]\leq C_2 \mathbf{W}^{4R}_{\beta,p}[\omega],
      \end{align}
       holds  almost everywhere in $B_{2R}(x_0)$, for some $C_2>0$.
   \end{proposition}                                                                                       
\nind\Proof{\it Assertion (i)}. First we assume that $\mu$ has compact support. Let $x\in\BBR^N$ and $t>0$. For any $y\in B_t(x)$, 
   \begin{align*}
   \mathbf{W}_{\alpha,p}[\chi_{_{B_t(x)}}\mu](y)&\geq \int_{2t}^{+\infty}\left(\frac{\mu(B_t(x)\cap B_r(y))}{r^{N-\alpha p}}\right)^{\frac{1}{p-1}}\frac{dr}{r}\\&\geq 
   \int_{2t}^{+\infty}\left(\frac{\mu(B_t(x))}{r^{N-\alpha p}}\right)^{\frac{1}{p-1}}\frac{dr}{r}
   \\&\geq C\left(\frac{\mu(B_t(x))}{t^{N-\alpha p}}\right)^{\frac{1}{p-1}}.
   \end{align*}
From \eqref{pot-2} we have 
   \begin{align*}
   \mu(B_t(x))\geq C\int_{B_t(x)}\left(\mathbf{W}_{\alpha,p}[\chi_{_{B_t(x)}}\mu](y)\right)^qdy
   \geq Ct^N\left(\frac{\mu(B_t(x))}{t^{N-\alpha p}}\right)^{\frac{1}{p-1}}.
   \end{align*}
Hence, $ \mu(B_t(x))\leq C t^{N-\frac{\alpha pq}{q-p+1}}.$ Therefore
\begin{align}\label{pot-11}
\int_{r}^{\infty}\left(\frac{\mu(B_t(x))}{t^{N-\alpha p}}\right)^{\frac{1}{p-1}}\frac{dt}{t}\leq C r^{-\frac{\alpha p}{q-p+1}}.
\end{align}
Since, $B_t(y)\subset B_{2\max\{t,r\}}(x)$ for any $y\in B_r(x)$, we have
\begin{align*}
&\int_{B_r(x)}\left(\mathbf{W}_{\alpha,p}[\mu](y)\right)^qdy\leq C \int_{B_r(x)}\left(\int_{0}^{r}\left(\frac{\mu(B_t(y)\cap B_{2r}(x))}{t^{N-\alpha p}}\right)^{\frac{1}{p-1}}\frac{dt}{t}\right)^qdy
\\& ~~~~~~~~~~~~~~~~~ \phantom{------}+C  \int_{B_r(x)}\left(\int_{r}^{\infty}\left(\frac{\mu(B_t(y)\cap B_{2t}(x))}{t^{N-\alpha p}}\right)^{\frac{1}{p-1}}\frac{dt}{t}\right)^qdy
\\&~~~~~~~~~~~~~~\leq C \int_{B_{r}(x)}\left(\mathbf{W}_{\alpha,p}[\chi_{_{B_{2r}(x)}}\mu]\right)^qdy +Cr^N\left(\int_{r}^{\infty}\left(\frac{\mu(B_{2t}(x))}{t^{N-\alpha p}}\right)^{\frac{1}{p-1}}\frac{dt}{t}\right)^q
\\&\phantom{\int_{B_r(x)}\left(\mathbf{W}_{\alpha,p}[\mu](y)\right)^qdy}\leq C \mu(B_{2r}(x)) +Cr^N\left(\int_{r}^{\infty}\left(\frac{\mu(B_{2t}(x))}{t^{N-\alpha p}}\right)^{\frac{1}{p-1}}\frac{dt}{t}\right)^q.
\end{align*} 
Note that, in the  last inequality, we have  used \eqref{pot-2}.  Thus, 
\begin{align*}
&\mathbf{W}_{\beta,p}\left[\left(\mathbf{W}_{\alpha,p}[\mu]\right)^q\right](x)=\int_{0}^{\infty}\left(\frac{\int_{B_r(x)}\left(\mathbf{W}_{\alpha,p}[\mu](y)\right)^qdy}{r^{N-\beta p}}\right)^{\frac{1}{p-1}}\frac{dr}{r}\\& ~~~~\leq 
C\int_{0}^{\infty}\left(\frac{\mu(B_{2r}(x))}{r^{N-\beta p}}\right)^{\frac{1}{p-1}}\frac{dr}{r}+C\int_{0}^{\infty}r^{\frac{\beta p}{p-1}-1}\left(\int_{r}^{\infty}\left(\frac{\mu(B_{2t}(x))}{t^{N-\alpha p}}\right)^{\frac{1}{p-1}}\frac{dt}{t}\right)^{\frac{q}{p-1}}dr.
\end{align*}
Therefore, it remains to prove 
\begin{align*}
\int_{0}^{\infty}r^{\frac{\beta p}{p-1}-1}\left(\int_{r}^{\infty}\left(\frac{\mu(B_{2t}(x))}{t^{N-\alpha p}}\right)^{\frac{1}{p-1}}\frac{dt}{t}\right)^{\frac{q}{p-1}}dr\leq C\mathbf{W}_{\beta,p}[\mu](x).
\end{align*}
Notice that 
\begin{align*}
r^{\frac{\beta p}{p-1}}\left(\int_{r}^{\infty}\left(\frac{\mu(B_{2t}(x))}{t^{N-\alpha p}}\right)^{\frac{1}{p-1}}\frac{dt}{t}\right)^{\frac{q}{p-1}}\to 0 \text{ as }~t\to 0
\end{align*}
and 
\begin{align*}
r^{\frac{\beta p}{p-1}}\left(\int_{r}^{\infty}\left(\frac{\mu(B_{2t}(x))}{t^{N-\alpha p}}\right)^{\frac{1}{p-1}}\frac{dt}{t}\right)^{\frac{q}{p-1}}\leq C r^{\frac{\beta p}{p-1}-\frac{N-\alpha p}{p-1}\frac{q}{p-1}} (\mu (\mathbb{R}^N))^{\frac{q}{(p-1)^2}}\to 0
\end{align*}
 as $t\to \infty$, since $\frac{\beta p}{p-1}-\frac{N-\alpha p}{p-1}\frac{q}{p-1}<\frac{\beta p}{p-1}-\frac{N}{p-1}<0$.
 Hence, using integration by parts and inequality \eqref{pot-11}, we have
\begin{align*}
&\int_{0}^{\infty}r^{\frac{\beta p}{p-1}-1}\left(\int_{r}^{\infty}\left(\frac{\mu(B_{2t}(x))}{t^{N-\alpha p}}\right)^{\frac{1}{p-1}}\frac{dt}{t}\right)^{\frac{q}{p-1}}dr\\& =\frac{q}{\beta p}\int_0^\infty r^{\frac{\beta p}{p-1}}\left(\int_{r}^{\infty}\left(\frac{\mu(B_{2t}(x))}{t^{N-\alpha p}}\right)^{\frac{1}{p-1}}\frac{dt}{t}\right)^{\frac{q}{p-1}-1}\left(\frac{\mu(B_{2r}(x))}{r^{N-\alpha p}}\right)^{\frac{1}{p-1}}\frac{dr}{r}
\\&\leq C\int_0^\infty r^{\frac{\beta p}{p-1}}\left(r^{-\frac{\alpha p}{q-p+1}}\right)^{\frac{q}{p-1}-1}\left(\frac{\mu(B_{2r}(x))}{r^{N-\alpha p}}\right)^{\frac{1}{p-1}}\frac{dr}{r}
\\&=C\mathbf{W}_{\beta,p}[\mu](x).
\end{align*}
Next, we assume that $\mu$ is not necessarily compactly supported. From the previous step,  
 \begin{align*}
    \mathbf{W}_{\beta,p}\left[\left(\mathbf{W}_{\alpha,p}[\chi_{_{B_n(0)}}\mu]\right)^q\right]\leq C \mathbf{W}_{\beta,p}[\chi_{_{B_n(0)}}\mu]\leq C \mathbf{W}_{\beta,p}[\mu]<\infty ~\text{ a.e in }~\mathbb{R}^N.
    \end{align*}
Then we derive \eqref{pot-9} by Fatou's lemma. \smallskip

\nind{\it Assertion (ii)}. 
For any $x\in B_{2R}(x_0)$, $0<t<R/2$ and $y\in B_t(x)$, 
   \begin{align*}
   \mathbf{W}^{4R}_{\alpha,p}[\chi_{_{B_t(x)}}\omega](y)&\geq \int_{2t}^{4R}\left(\frac{\omega(B_t(x)\cap B_r(y))}{r^{N-\alpha p}}\right)^{\frac{1}{p-1}}\frac{dr}{r}
   \\&\geq C\left(\frac{\omega(B_t(x))}{t^{N-\alpha p}}\right)^{\frac{1}{p-1}}.
   \end{align*}
From \eqref{pot-5} we have 
   \begin{align*}
   \omega(B_t(x))\geq C\int_{B_t(x)}\left(\mathbf{W}^{4R}_{\alpha,p}[\chi_{_{B_t(x)}}\omega](y)\right)^qdy
   \geq Ct^N\left(\frac{\omega(B_t(x))}{t^{N-\alpha p}}\right)^{\frac{1}{p-1}}.
   \end{align*}
Hence, $\omega(B_t(x))\leq C t^{N-\frac{\alpha pq}{q-p+1}}$ for all $t\in (0,R/2)$ and $x\in B_{2R}(x_0)$. It implies 
\begin{align}\label{pot-12}
\int_{r}^{4R}\left(\frac{\omega(B_{2t}(x))}{t^{N-\alpha p}}\right)^{\frac{1}{p-1}}\frac{dt}{t}\leq C r^{-\frac{\alpha p}{q-p+1}} ~~\text{for all}~~x\in B_{2R}(x_0)\,,\; 0<r<4R.
\end{align}
Since $B_t(y)\subset B_{2\max\{t,r\}}(x)$ for any $0<r<4R$ and $y\in B_r(x)$, 
\begin{align*}
&\int_{B_r(x)}\left(\mathbf{W}^{4R}_{\alpha,p}[\omega](y)\right)^qdy\leq C \int_{B_r(x)}\left(\int_{0}^{r}\left(\frac{\omega(B_t(y)\cap B_{2r}(x))}{t^{N-\alpha p}}\right)^{\frac{1}{p-1}}\frac{dt}{t}\right)^qdy
\\&\phantom{\int_{B_r(x)}\left(\mathbf{W}^{4R}_{\alpha,p}[\omega](y)\right)^qdy} +C  \int_{B_r(x)}\left(\int_{r}^{4R}\left(\frac{\omega(B_t(y)\cap B_{2t}(x))}{t^{N-\alpha p}}\right)^{\frac{1}{p-1}}\frac{dt}{t}\right)^qdy
\\&~~~~~~~~~~~~~\leq C \int_{B_{r}(x)}\left(\mathbf{W}^{4R}_{\alpha,p}[\chi_{_{B_{2r}(x)}}\omega]\right)^qdy +Cr^N\left(\int_{r}^{4R}\left(\frac{\mu(B_{2t}(x))}{t^{N-\alpha p}}\right)^{\frac{1}{p-1}}\frac{dt}{t}\right)^q
\\&\phantom{\int_{B_r(x)}\left(\mathbf{W}^{4R}_{\alpha,p}[\omega](y)\right)^qdy}\leq C \mu(B_{2r}(x)) +Cr^N\left(\int_{r}^{4R}\left(\frac{\mu(B_{2t}(x))}{t^{N-\alpha p}}\right)^{\frac{1}{p-1}}\frac{dt}{t}\right)^q.
\end{align*} 
In the last inequality we have used \eqref{pot-5}. Thus,
as above, we only need to prove that 
\begin{align*}
\int_{0}^{4R}r^{\frac{\beta p}{p-1}-1}\left(\int_{r}^{4R}\left(\frac{\omega(B_{2t}(x))}{t^{N-\alpha p}}\right)^{\frac{1}{p-1}}\frac{dt}{t}\right)^{\frac{q}{p-1}}dr\leq C\mathbf{W}_{\beta,p}^{2R}[\omega](x).
\end{align*}
Using integration by parts and \eqref{pot-12}
\begin{align*}
&\int_{0}^{4R}r^{\frac{\beta p}{p-1}-1}\left(\int_{r}^{4R}\left(\frac{\omega(B_{2t}(x))}{t^{N-\alpha p}}\right)^{\frac{1}{p-1}}\frac{dt}{t}\right)^{\frac{q}{p-1}}dr\\& =\frac{q}{\beta p}\int_0^{4R} r^{\frac{\beta p}{p-1}}\left(\int_{r}^{4R}\left(\frac{\omega(B_{2t}(x))}{t^{N-\alpha p}}\right)^{\frac{1}{p-1}}\frac{dt}{t}\right)^{\frac{q}{p-1}-1}\left(\frac{\omega(B_{2r}(x))}{r^{N-\alpha p}}\right)^{\frac{1}{p-1}}\frac{dr}{r}
\\&\leq C\int_0^{4R} r^{\frac{\beta p}{p-1}}\left(r^{-\frac{\alpha p}{q-p+1}}\right)^{\frac{q}{p-1}-1}\left(\frac{\omega(B_{2r}(x))}{r^{N-\alpha p}}\right)^{\frac{1}{p-1}}\frac{dr}{r}
\\&\leq C\mathbf{W}^{8R}_{\beta,p}[\omega](x)\leq C\mathbf{W}^{4R}_{\beta,p}[\omega](x), 
\end{align*}
 since $\mathbf{W}^{8R}_{\beta,p}[\omega](x)\leq C \mathbf{W}^{4R}_{\beta,p}[\omega](x)$ for any $x\in B_{2R}(x_0)$, because $\supp \omega\subset B_{2R}(x_0)$.
\qeda\medskip 
  
 The next result is at the core of our construction since it connects the integral of product of Wolff potentials to some power to  the integral of a new Wolff potential. In this highly technical construction, the role of Hardy Littlewood maximal function plays an important role as well as classical tools from harmonic analysis such as the Vitali Covering Lemma.
\begin{theorem}\label{th-pot-3}Let $\alpha$, $\beta$, $q_1$, $q_2>0$,  $\alpha>\beta$, $1<p<\frac N\alpha$, $q_1+q_2>p-1$,  $q_2<\frac{N(p-1)}{N-\beta p}$ and $\frac{\alpha p q_1+\beta pq_2}{q_1+q_2}<N$. Then, there holds
 \bel{pot-15}\BA {lll}\displaystyle
\int_{\mathbb{R}^N}\left(\mathbf{M}_{\frac{\alpha p q_1+\beta pq_2}{q_1+q_2}}[\mu](x)\right)^{\frac{q_1+q_2}{p-1}}dx\asymp \int_{\mathbb{R}^N}\left(\mathbf{W}_{\frac{\alpha pq_1+\beta pq_2}{q_1+q_2},p}[\mu](x)\right)^{q_1+q_2}dx
\\\phantom{\int_{\mathbb{R}^N}\left(\mathbf{M}_{\frac{\alpha p q_1+\beta pq_2}{q_1+q_2}}[\mu](x)\right)^{\frac{q_1+q_2}{p-1}}dx}\displaystyle
\asymp \int_{\mathbb{R}^N}\left(\mathbf{W}_{\alpha,p}[\mu](x)\right)^{q_1}\left(\mathbf{W}_{\beta,p}[\mu](x)\right)^{q_2}dx,
 \EA\ee
 for any  $\mu \in\mathfrak{M}^+(\mathbb{R}^N)$, and
 \bel{pot-16}\BA {lll}\displaystyle
\int_{\mathbb{R}^N}\left(\mathbf{M}^{2R}_{\frac{\alpha p q+\beta pq_2}{q_1+q_2}}[\omega](x)\right)^{\frac{q_1+q_2}{p-1}}\asymp \int_{\mathbb{R}^N}\left(\mathbf{W}^{2R}_{\frac{\alpha pq_1+\beta pq_2}{q_1+q_2},p}[\omega](x)\right)^{q_1+q_2}dx\\
\phantom {\int_{\mathbb{R}^N}\left(\mathbf{M}^{2R}_{\frac{\alpha p q+\beta pq_2}{q_1+q_2}}[\omega](x)\right)^{\frac{q_1+q_2}{p-1}}}\displaystyle\asymp \int_{\mathbb{R}^N}\left(\mathbf{W}^{2R}_{\alpha,p}[\omega](x)\right)^{q_1}\left(\mathbf{W}_{\beta,p}^{2R}[\omega](x)\right)^{q_2}dx, 
  \EA\ee
  for any $R>0$ and $\omega\in \mathfrak{M}^+(\mathbb{R}^N)$ with $diam(\text{supp}~\omega)\leq R$. 
     \end{theorem}                                                                                       

 For proving this theorem we need several intermediate results. For any $\alpha\in (0,N)$, $s>0$, $R\in (0,\infty]$ we denote 
 \bel{pot-16-1}
  \mathbf{L}_{\alpha,s}^R[\mu](x)=\int_{0}^{R}\left(\frac{\mu(B_t(x))}{t^{N-\alpha}}\right)^s\frac{dt}{t},
  \ee
and $\mathbf{L}_{\alpha,s}[\mu]:=\mathbf{L}_{\alpha,s}^\infty[\mu]$ when $R=\infty$. We notice that $\mathbf{L}_{\alpha,s}^R$ is actually a Wolff potential since
 \bel{pot-16-2}
  \mathbf{L}_{\alpha,s}^R[\mu]= \mathbf{W}_{\frac{\ga s}{s+1},\frac{s+1}{ s}}^R[\mu]\quad\text{and }\;\; 
   \mathbf{L}_{\alpha,s}[\mu]= \mathbf{W}_{\frac{\ga s}{s+1},\frac{s+1}{ s}}[\mu].
  \ee

\begin{lemma}\label{lem-pot-1} Let $\alpha_1$, $\alpha_2$, $s_1$, $s_2>0$, $0<\alpha_2<\alpha_1<N$. There exist $C=C(N,\alpha_1,\alpha_2,s_1,s_2)>0$ and $\varepsilon_0=\varepsilon_0(N,\alpha_1,\alpha_2,s_1,s_2)>0$ such that for any $\mu\in\mathfrak{M}_+(\mathbb{R}^N)$, $R\in (0,\infty]$, $\varepsilon\in (0,\varepsilon_0)$ and $\lambda>0$, the inequality
\begin{align}\label{pot-17}
\left|\left\{\mathbf{L}_{\alpha_1,s_1}^{2R}[\mu]\mathbf{L}_{\alpha_2,s_2}^{2R}[\mu]>\varepsilon^{1/2}\lambda\right\}\right|<\infty,
\end{align}
implies
 \bel{pot-18}\BA {lll}
 \left|\left\{{\phantom {\myint{}{r^p}}}\!\!\!\!\!\!\!\!\!\!\!\mathbf{L}_{\alpha_1,s_1}^{R}[\mu]\mathbf{L}_{\alpha_2,s_2}^R[\mu]>a\lambda\right\}\cap\left\{\left(\mathbf{M}^{2R}_{\frac{\alpha_1 s_1+\alpha_2s_2}{s_1+s_2}}[\mu]\right)^{s_1+s_2}\leq\varepsilon\lambda\right\}\right|\\[3mm]\phantom{---------------}
 \leq  C\varepsilon^{\frac{ N}{2s_2(N-\alpha_2)}}\left|\left\{\mathbf{L}_{\alpha_1,s_1}^{2R}[\mu]\mathbf{L}_{\alpha_2,s_2}^{2R}[\mu]>\varepsilon^{1/2}\lambda\right\}\right|.
\EA \ee
 \end{lemma}                                                                                       

 To prove this, we need the following two lemmas:
\begin{lemma}\label{lem-pot-2} Let $0<\alpha<N$ and $s>0$. There exists $C=C(N,\alpha,s)$ such that 
 \begin{align}\label{011102141}
 \left| \left\{\mathbf{L}_{\alpha,s}[\omega]>\lambda\right\} \right|\leq C\left(\frac{(\omega(\mathbb{R}^N))^s}{\lambda}\right)^{\frac{N}{s(N-\alpha)}} \quad\text{for all }\;\lambda>0,
 \end{align}
 for any $\omega\in\mathfrak{M}_b^+(\mathbb{R}^N)$. 
 \end{lemma}
 \Proof  It is easy to see that  $\mathbf{L}_{\alpha,s}[\omega](x)\leq C\left(\mathbf{M}(\omega)(x)\right)^{\frac{s(N-\alpha)}{N}}(\omega(\mathbb{R}^N))^{\frac{\alpha s}{N}}$. Thus, thanks to  boundedness of the operator $\mathbf{M}$ from $\mathfrak{M}_b^+(\mathbb{R}^N)$ to $L^{1,\infty}(\mathbb{R}^d)$, we get \eqref{011102141}. The proof is complete.
 \qeda\\
 
 The next result is  a consequence of Vitali Covering Lemma.
\begin{lemma}\label{lem-pot-3}
                  Let $0<\varepsilon<1, R>0$ and $B:=B_R(x_0)$ for some $x_0\in \mathbb{R}^{N}$. Let $E\subset F\subset B$ be two measurable sets in $\mathbb{R}^{N}$ with $|E|<\varepsilon |B|$
                   and satisfying 
                   the following property: for all $x\in B$ and $r\in (0,R]$, we have
                  $B_r(x)\cap B\subset F$
                   provided $|E\cap B_r(x)|\geq \varepsilon |B_r(x)|.
                     $
                   Then $|E|\leq C\varepsilon |F|$ for some $C=C(N)$.             
                      \end{lemma}                                                                                       

 \nind{\it Proof of Lemma \ref{lem-pot-1}}. We only consider the case $R<\infty$, the case $R=\infty$ being similar. Let $\{B_{R}(x_j)\}$ be a cover of $\mathbb{R}^{N}$ such that, for some constant $M=M(N)>0$,
   $$
   \sum_j \chi_{_{B_{R/4}(x_j)}}(x)\leq M\qquad\text{for all }\;x\in\mathbb{R}^{N}.
  $$
It is sufficient to show that  there exist constants $c_1,c_2>0$ and $\varepsilon_0\in (0,1)$ depending on $N,\alpha_1,\alpha_2,s_1,s_2,p$ such that for any $B\in\{B_{R/4}(x_j)\}$, $\lambda>0$ and $\varepsilon\in (0,\varepsilon_0)$, there holds

  \bel{pot-19}\BA {llll}
\left|B\cap\left\{\mathbf{L}_{\alpha_1,s_1}^R[\mu]\mathbf{L}_{\alpha_1,s_1}^R[\mu]>a\lambda^{\!\!\!\!\!\!\!\!\!\!\!\!\phantom{p^{J^{J^{\int}}}}}\right\}\cap\left\{\left(\mathbf{M}^{2R}_{\frac{\alpha_1 s_1+\alpha_2s_2}{s_1+s_2}}[\mu]\right)^{s_1+s_2}\leq\varepsilon\lambda\right\}\right|\\[4mm]\phantom{------------}
\leq C\varepsilon^{\frac{ N}{2s_2(N-\alpha_2)}} \left|B\cap\left\{\mathbf{L}_{\alpha_1,s_1}^{2R}[\mu]\mathbf{L}_{\alpha_1,s_1}^{2R}[\mu]>\varepsilon^{1/2}\lambda^{\!\!\!\!\!\!\!\!\!\!\!\!\phantom{p^{J^{J^{\int}}}}}\right\}\right|,
 \EA \ee
  where 
 \begin{align*}
 a=1+\left(\frac{s_1+s_2}{s_1s_2(\alpha_1-\alpha_2)}\right)^2+2^{(N-\alpha_1)s_1+(N-\alpha_2)s_2+1}.
 \end{align*}
 Fix $\lambda>0$ and $0<\varepsilon<\min\{1/10,2^{-\frac{10s_1s_2(\alpha_1-\alpha_2)}{s_1+s_2}}\}$. We set 
 \begin{align*}
 E=B\cap\left\{\mathbf{L}_{\alpha_1,s_1}^R[\mu]\mathbf{L}_{\alpha_1,s_1}^R[\mu]>a\lambda^{\!\!\!\!\!\!\!\!\!\!\!\!\phantom{p^{J^{J^{\int}}}}}\right\}
  \cap\left\{\left(\mathbf{M}^{2R}_{\frac{\alpha_1 s_1+\alpha_2s_2}{s_1+s_2}}[\mu]\right)^{s_1+s_2}\leq\varepsilon\lambda\right\},
 \end{align*}
 and 
  \begin{align*}
  F=B\cap\left\{\mathbf{L}_{\alpha_1,s_1}^{2R}[\mu]\mathbf{L}_{\alpha_1,s_1}^{2R}[\mu]>\varepsilon^{1/2}\lambda\right\}.
  \end{align*}   Thanks to Lemma \ref{lem-pot-3} we will obtain \eqref{pot-19} provided we verify the following two claims:
                  \begin{equation}\label{pot-20}
                 |E|\leq C\varepsilon^{\frac{ N}{2s_2(N-\alpha_2)}}|B|,
                  \end{equation}
                  and, for any $x\in B$ and $0<r\leq R/4$, 
                  \begin{equation}\label{pot-21}
                  |E\cap B_r(x)|< C\varepsilon^{\frac{ N}{2s_2(N-\alpha_2)}}|B_r(x)|,
                  \end{equation}
                    whenever $B_r(x)\cap B\cap F^c\not =\emptyset$ and $E\cap B_r(x)\not = \emptyset.$
                     \smallskip
                     
\nind{\it Proof of \eqref{pot-20}}:   For any $x\in E$, we have 
                  \begin{align*}
                \mathbf{L}_{\alpha_1,s_1}^R[\mu](x)& \leq \int_{0}^{R}\left(t^{-\frac{\alpha_1 s_1+\alpha_2s_2}{s_1+s_2}+\alpha_1}\mathbf{M}^{2R}_{\frac{\alpha_1 s_1+\alpha_2s_2}{s_1+s_2}}[\mu](x)\right)^{s_1}\frac{dt}{t}\\& \leq \frac{s_1+s_2}{s_1s_2(\alpha_1-\alpha_2)}R^{\frac{s_1s_2(\alpha_1-\alpha_2)}{s_1+s_2}}(\varepsilon\lambda)^{\frac{s_1}{s_1+s_2}}.              
                  \end{align*}                 
Hence, the inequality $\mathbf{L}_{\alpha_1,s_1}^R[\mu](x)\mathbf{L}_{\alpha_2,s_2}^R[\mu](x)>\lambda$ implies 
                  \begin{align*}
                  \mathbf{L}_{\alpha_2,s_2}^R[\mu](x)>\frac{s_1s_2(\alpha_1-\alpha_2)}{s_1+s_2}R^{-\frac{s_1s_2(\alpha_1-\alpha_2)}{s_1+s_2}}\varepsilon^{-\frac{s_1}{s_1+s_2}}\lambda^{\frac{s_2}{s_1+s_2}}.
                  \end{align*}
            Clearly,   $ \mathbf{L}_{\alpha_2,s_2}^R[\mu]=\mathbf{L}_{\alpha_2,s_2}^R[\chi_{_{B_{2R}}(y_0)}\mu]$ in $B$ for any $y_0\in B$. Fix  $y_0\in E$,  we have 
            \begin{align*}
            |E|\leq \left|\left\{\mathbf{L}_{\alpha_2,s_2}[\chi_{_{B_{2R}(y_0)}}\mu]>\frac{s_1s_2(\alpha_1-\alpha_2)}{s_1+s_2}R^{-\frac{s_1s_2(\alpha_1-\alpha_2)}{s_1+s_2}}\varepsilon^{-\frac{s_1}{s_1+s_2}}\lambda^{\frac{s_2}{s_1+s_2}}\right\}\right|.
            \end{align*}
            Using \eqref{pot-18} from Lemma \ref{lem-pot-2} and the fact that  $\left(\mathbf{M}^{2R}_{\frac{\alpha_1 s_1+\alpha_2s_2}{s_1+s_2}}[\mu](y_0)\right)^{s_1+s_2}\leq\varepsilon\lambda$, we get
            \begin{align*}
            |E|&\leq C\left(\frac{(\mu(B_{2R}(y_0)))^{s_2}}{R^{-\frac{s_1s_2(\alpha_1-\alpha_2)}{s_1+s_2}}\varepsilon^{-\frac{s_1}{s_1+s_2}}\lambda^{\frac{s_2}{s_1+s_2}}}\right)^{\frac{N}{s_2(N-\alpha_2)}}
            \\& \leq  C\left(\frac{(\varepsilon\lambda)^{\frac{s_2}{s_1+s_2}}(2R)^{s_2N-s_2\frac{\alpha_1 s_1+\alpha_2s_2}{s_1+s_2}}}{R^{-\frac{s_1s_2(\alpha_1-\alpha_2)}{s_1+s_2}}\varepsilon^{-\frac{s_1}{s_1+s_2}}\lambda^{\frac{s_2}{s_1+s_2}}}\right)^{\frac{N}{s_2(N-\alpha_2)}}
            \\& =  C \varepsilon^{\frac{N}{s_2(N-\alpha_2)}}|B|
            \\& \leq   C \varepsilon^{\frac{N}{2s_2(N-\alpha_2)}}|B|.
            \end{align*}                
             We obtain \eqref{pot-20}.  \smallskip
                     
\nind{\it Proof of \eqref{pot-21}}: 
                    Take $x\in B$ and $0<r\leq R/4$.
                            Now assume that $B_r(x)\cap B\cap F^c\not= \emptyset$ and $E\cap B_r(x)\not = \emptyset$, then  there exists $x_1\in B_r(x)\cap B$ such that $\mathbf{L}_{\alpha_1,s_1}^{2R}[\mu](x_1)\mathbf{L}_{\alpha_2,s_2}^{2R}[\mu](x_1)\leq \varepsilon^{1/2}\lambda$.
                             We need to prove that
                             \begin{equation}\label{pot-22}
                                    |E\cap B_r(x)|< C\varepsilon^{\frac{ N}{2s_2(N-\alpha_2)}} |B_r(x)|. 
                                                     \end{equation}
                          To do this, we can write 
                        \begin{align*}
                        \mathbf{L}_{\alpha_1,s_1}^R[\mu](y)\mathbf{L}_{\alpha_1,s_1}^R[\mu](y)= \mathbf{T}_1(y) +\mathbf{T}_2(y)+\mathbf{T}_3(y)+\mathbf{T}_4(y),      
                          \end{align*}
                          where 
                         \begin{align*}
                         &\mathbf{T}_1(y)=\mathbf{L}_{\alpha_1,s_1}^{4r}[\mu](y)\mathbf{L}_{\alpha_2,s_2}^{4r}[\mu](y),\\&
                          \mathbf{T}_2(y)=\mathbf{L}_{\alpha_1,s_1}^{4r}[\mu](y)\int_{4r}^{R}\left(\frac{\mu(B_\rho(y))}{\rho^{N-\alpha_2}}\right)^{s_2}\frac{d\rho}{\rho},
                          \\&
                            \mathbf{T}_3(y)=\int_{4r}^{R}\left(\frac{\mu(B_\rho(y))}{\rho^{N-\alpha_1}}\right)^{s_1}\frac{d\rho}{\rho}\mathbf{L}_{\alpha_2,s_2}^{4r}[\mu](y)  ,\\&
                             \mathbf{T}_4(y)=\int_{4r}^{R}\left(\frac{\mu(B_\rho(y))}{\rho^{N-\alpha_1}}\right)^{s_1}\frac{d\rho}{\rho}\int_{4r}^{R}\left(\frac{\mu(B_\rho(y))}{\rho^{N-\alpha_2}}\right)^{s_2}\frac{d\rho}{\rho}.\end{align*}
                     For all $y\in E\cap B_r(x)$, we have
         \bel{pot-23}\BA {llll} \displaystyle
                 \mathbf{T}_2(y)\leq \int_{0}^{4r}\left(\rho^{-\frac{\alpha_1 s_1+\alpha_2s_2}{s_1+s_2}+\alpha_1}\mathbf{M}^R_{\frac{\alpha_1 s_1+\alpha_2s_2}{s_1+s_2}}[\mu](y)\right)^{s_1}\frac{d\rho}{\rho}\\[4mm] \phantom{-----------}\displaystyle
                 \times \int_{4r}^{R}\left(\rho^{-\frac{\alpha_1 s_1+\alpha_2s_2}{s_1+s_2}+\alpha_2}\mathbf{M}^R_{\frac{\alpha_1 s_1+\alpha_2s_2}{s_1+s_2}}[\mu](y)\right)^{s_2}\frac{d\rho}{\rho}\\[4mm] \phantom{\mathbf{T}_2(y)}
                  \displaystyle\leq 
                 \left(\frac{s_1+s_2}{s_1s_2(\alpha_1-\alpha_2)}\right)^2 r^{\frac{(\alpha_1-\alpha_2)s_1s_2}{s_1+s_2}}r^{-\frac{(\alpha_1-\alpha_2)s_1s_2}{s_1+s_2}}\left(\mathbf{M}^R_{\frac{\alpha_1 s_1+\alpha_2s_2}{s_1+s_2}}[\mu](y)\right)^{s_1+s_2}
                 \\[4mm]\phantom{\mathbf{T}_2(y)}
                  \displaystyle\leq \left(\frac{s_1+s_2}{s_1s_2(\alpha_1-\alpha_2)}\right)^2 \varepsilon\lambda,
\EA\ee 
also    
   \bel{pot-24}\BA {lll}
   \displaystyle           \mathbf{T}_4(y)\leq \int_{4r}^{R}\left(\frac{\mu(B_{2\rho}(x_1))}{\rho^{N-\alpha_1}}\right)^{s_1}\frac{d\rho}{\rho}\int_{4r}^{R}\left(\frac{\mu(B_{2\rho}(x_1))}{\rho^{N-\alpha_2}}\right)^{s_2}\frac{d\rho}{\rho}\\[4mm]\phantom{\mathbf{T}_4(y)}\displaystyle\leq 2^{(N-\alpha_1)s_1+(N-\alpha_2)s_2}
               \mathbf{L}_{\alpha_1,s_1}^{2R}[\mu](x_1)\mathbf{L}_{\alpha_2,s_2}^{2R}[\mu](x_1)\\[2mm]\phantom{\mathbf{T}_4(y)}\displaystyle\leq 2^{(N-\alpha_1)s_1+(N-\alpha_2)s_2}\lambda,\EA\ee
and
       \bel{pot-25}\BA {lll}\displaystyle
                      \mathbf{T}_3(y)\leq \int_{4r}^{R}\left(\frac{\mu(B_{2\rho}(x_1))}{\rho^{N-\alpha_1}}\right)^{s_1}\frac{d\rho}{\rho}\mathbf{L}_{\alpha_2,s_2}^{4r}[\mu](y)\\[4mm]\phantom{\mathbf{T}_4(y)}\displaystyle\leq 2^{s_1(N-\alpha_1)}
                       \mathbf{L}_{\alpha_1,s_1}^{2R}[\mu](x_1)\mathbf{L}_{\alpha_2,s_2}^{4r}[\mu](y),
\EA\ee
 Thus, 
  \begin{align*}
  |E\cap B_r(x)|&\leq \mathbf{Y}_1+\mathbf{Y}_2+\mathbf{Y}_3+\mathbf{Y}_4,
  \end{align*}
  where 
  \begin{align*}
 & \mathbf{Y}_1=\left|E\cap B_r(x)\cap\left\{\mathbf{T}_1>\lambda^{\phantom{I^J}}\!\!\!\!\!\right\}\cap\left\{\left(\mathbf{M}^{2R}_{\frac{\alpha_1 s_1+\alpha_2s_2}{s_1+s_2}}[\mu]\right)^{s_1+s_2}\leq\varepsilon\lambda\right\}\right|,\\&
   \mathbf{Y}_2=\left|E\cap B_r(x)\cap\left\{\mathbf{T}_2>\left(\frac{s_1+s_2}{s_1s_2(\alpha_1-\alpha_2)}\right)^2 \lambda^{\phantom{I^J}}\!\!\!\!\!\right\}\cap\left\{\left(\mathbf{M}^{2R}_{\frac{\alpha_1 s_1+\alpha_2s_2}{s_1+s_2}}[\mu]\right)^{s_1+s_2}\leq\varepsilon\lambda\right\}\right|,\\&
    \mathbf{Y}_3=\left|E\cap B_r(x)\cap\left\{\mathbf{T}_3>2^{s_1(N-\alpha_1)}\lambda^{\phantom{I^J}}\!\!\!\!\!\right\}\cap\left\{\left(\mathbf{M}^{2R}_{\frac{\alpha_1 s_1+\alpha_2s_2}{s_1+s_2}}[\mu]\right)^{s_1+s_2}\leq\varepsilon\lambda\right\}\right|,\\&
     \mathbf{Y}_4=\left|E\cap B_r(x)\cap\left\{\mathbf{T}_4>2^{(N-\alpha_1)s_1+(N-\alpha_2)s_2}\lambda^{\phantom{I^J}}\!\!\!\!\!\right\}\cap\left\{\left(\mathbf{M}^{2R}_{\frac{\alpha_1 s_1+\alpha_2s_2}{s_1+s_2}}[\mu]\right)^{s_1+s_2}\leq\varepsilon\lambda\right\}\right|.
  \end{align*}
        
          As in the proof of \eqref{pot-20}, it can be shown that
          \begin{align}\label{pot-25*}
          \mathbf{Y}_1\leq c_{11} \varepsilon^{\frac{N}{2s_2(N-\alpha_2)}}|B_r(x)|.
          \end{align}
          From \eqref{pot-23}-\eqref{pot-25}, we obtain 
                         $\mathbf{Y}_2=\mathbf{Y}_4=0$ and 
                         \begin{align*}
                         \mathbf{Y}_3&\leq \left| B_r(x)\cap\left\{\mathbf{L}_{\alpha_2,s_2}^{4r}[\mu]>\lambda(\mathbf{L}_{\alpha_1,s_1}^{2R}[\mu](x_1))^{-1}\right\}\right|
                         \\&=\left| B_r(x)\cap\left\{\mathbf{L}_{\alpha_2,s_2}^{4r}[\chi_{_{B_{6r}(x_1)}}\mu]>\lambda(\mathbf{L}_{\alpha_1,s_1}^{2R}[\mu](x_1))^{-1}\right\}\right|
                         \end{align*}
                         since $B_{4r}(y)\subset B_{6r}(x_1)$ for all $y\in B_r(x)$.  Using \eqref{pot-18} from Lemma \ref{lem-pot-2}, we get
                         \begin{align*}
                         \mathbf{Y}_3&\leq C\left(\frac{(\mu(B_{6r}(x_1)))^{s_2}}{\lambda(\mathbf{L}_{\alpha_1,s_1}^{2R}[\mu](x_1))^{-1}}\right)^{\frac{N}{s_2(N-\alpha_2)}}\\&\leq C\left(\frac{\mathbf{L}_{\alpha_1,s_1}^{2R}[\mu](x_1)\mathbf{L}_{\alpha_2,s_2}^{2R}[\mu](x_1)}{\lambda}\right)^{\frac{N}{s_2(N-\alpha_2)}}r^N
                         \\&\leq C\varepsilon^{\frac{ N}{2s_2(N-\alpha_2)}}|B_r(x)|.
                         \end{align*}
 Combining these inequalities, we infer \eqref{pot-22}.   
        \qeda 
        \medskip

\nind {\it Proof of Theorem \ref{th-pot-3}}. {\it Step 1: Proof of \eqref{pot-15}}.
By \cite[Theorem 2.3]{55VHV}, we have 
  \begin{align*}
  \int_{\mathbb{R}^N}\left(\mathbf{M}_{\frac{\alpha p q_1+\beta pq_2}{q_1+q_2}}[\mu](x)\right)^{\frac{q_1+q_2}{p-1}}dx\asymp \int_{\mathbb{R}^N}\left(\mathbf{W}_{\frac{\alpha pq_1+\beta pq_2}{q_1+q_2},p}[\mu](x)\right)^{q_1+q_2}dx.
  \end{align*}
  Next, we prove 
  \begin{align*}
  \int_{\mathbb{R}^N}\left(\mathbf{M}_{\frac{\alpha p q_1+\beta pq_2}{q_1+q_2}}[\mu](x)\right)^{\frac{q_1+q_2}{p-1}}dx\asymp \int_{\mathbb{R}^N}\left(\mathbf{W}_{\alpha,p}[\mu](x)\right)^{q_1}\left(\mathbf{W}_{\beta,p}[\mu](x)\right)^{q_2}dx.
  \end{align*}Since  for all $x\in \mathbb{R}^N$ there holds
  \begin{align*}
  &\left(\mathbf{W}_{\alpha,p}[\mu](x)\right)^{q_1}\left(\mathbf{W}_{\beta,p}[\mu](x)\right)^{q_2}\geq C\left(\mathbf{M}_{\frac{\alpha p q_1+\beta pq_2}{q_1+q_2}}[\mu](x)\right)^{\frac{q_1+q_2}{p-1}},\\ &\left(\mathbf{W}_{\alpha,p}[\mu](x)\right)^{q_1}\leq C\left(\mathbf{L}_{\alpha p, \frac{q_1}{(q_1+q_2)(p-1)}}[\mu](x)\right)^{q_1+q_2},\\&
  \left(\mathbf{W}_{\beta,p}[\mu](x)\right)^{q_2}\leq C\left(\mathbf{L}_{\beta p, \frac{q_2}{(q_1+q_2)(p-1)}}[\mu](x)\right)^{q_1+q_2},
  \end{align*}
It is  therefore enough to show that 
  
  \bel{pot-26}\BA {lll}\displaystyle
   \int_{\mathbb{R}^N}\left(\mathbf{L}_{\alpha p, \frac{q_1}{(q_1+q_2)(p-1)}}[\mu]\mathbf{L}_{\beta p, \frac{q_2}{(q_1+q_2)(p-1)}}[\mu]\right)^{q_1+q_2}dx
  \leq  C \int_{\mathbb{R}^N}\left(\mathbf{M}_{\frac{\alpha p q_1+\beta pq_2}{q_1+q_2}}[\mu]\right)^{\frac{q_1+q_2}{p-1}}dx.
   \EA \ee 
  Set $d\mu_n=\chi_{_{B_n(0)}}d\mu$, then we have  
  \begin{align*}
  \left|\left\{\mathbf{L}_{\alpha p, \frac{q_1}{(q_1+q_2)(p-1)}}[\mu_n]\mathbf{L}_{\beta p, \frac{q_2}{(q_1+q_2)(p-1)}}[\mu_n]>t\right\}\right|<\infty~~~\text{for all } ~t>0.
  \end{align*}
   Hence, by Lemma \eqref{lem-pot-1}, there exist positive constants $C$, $\varepsilon_0$, $a$ such that for any $\lambda>0, \varepsilon\in (0,\varepsilon_0)$,
    \begin{align*}
   & \left|\left\{\mathbf{L}_{\alpha p, \frac{q_1}{(q_1+q_2)(p-1)}}[\mu_n]\mathbf{L}_{\beta p, \frac{q_2}{(q_1+q_2)(p-1)}}[\mu_n]>a\lambda\right\}\right|\\[2mm]&~~~~~~~~~~~\leq  C\varepsilon^{\frac{ N(q_1+q_2)(p-1)}{2q_2(N-\beta p)}}\left|\left\{\mathbf{L}_{\alpha p, \frac{q_1}{(q_1+q_2)(p-1)}}[\mu_n]\mathbf{L}_{\beta p, \frac{q_2}{(q_1+q_2)(p-1)}}[\mu_n]>\varepsilon^{1/2}\lambda\right\}\right|
   \\&\phantom{----------------------} 
   +
   \left|\left\{\left(\mathbf{M}_{\frac{\alpha p q_1+\beta pq_2}{q_1+q_2}}[\mu]\right)^{\frac{1}{p-1}}>\varepsilon\lambda\right\}\right|.
    \end{align*}
    Multiplying by $\lambda^{q_1+q_2-1}$ and integrating over $(0,\infty)$, we get
    \begin{align*}
    &\int_{0}^{\infty}\lambda^{q_1+q_2}\left|\left\{\mathbf{L}_{\alpha p, \frac{q_1}{(q_1+q_2)(p-1)}}[\mu_n]\mathbf{L}_{\beta p, \frac{q_2}{(q_1+q_2)(p-1)}}[\mu_n]>a\lambda\right\}\right|\frac{d\lambda}{\lambda}
    \\[2mm]&~~~\leq  C\varepsilon^{\frac{ N(q_1+q_2)(p-1)}{2q_2(N-\beta p)}}\int_{0}^{\infty}\lambda^{q_1+q_2}\left|\left\{\mathbf{L}_{\alpha p, \frac{q_1}{(q_1+q_2)(p-1)}}[\mu_n]\mathbf{L}_{\beta p, \frac{q_2}{(q_1+q_2)(p-1)}}[\mu_n]>\varepsilon^{1/2}\lambda\right\}\right|\frac{d\lambda}{\lambda}\\[2mm]&\phantom{------------}+
      \int_{0}^{\infty}\lambda^{q_1+q_2}\left|\left\{\left(\mathbf{M}_{\frac{\alpha p q_1+\beta pq_2}{q_1+q_2}}[\mu]\right)^{\frac{1}{p-1}}>\varepsilon\lambda\right\}\right|\frac{d\lambda}{\lambda}.
    \end{align*}
    By a change of variable, we derive 
    \begin{align*}
    &\left(a^{-q_1-q_2}-C\varepsilon^{\frac{ N(q_1+q_2)(p-1)}{2q_2(N-\beta p)}-\frac{q_1+q_2}{2}}\right)\\&~~~~~~\times\int_{0}^{\infty}\lambda^{q_1+q_2}\left|\left\{\mathbf{L}_{\alpha p, \frac{q_1}{(q_1+q_2)(p-1)}}[\mu_n]\mathbf{L}_{\beta p, \frac{q_2}{(q_1+q_2)(p-1)}}[\mu_n]>\lambda\right\}\right|\frac{d\lambda}{\lambda}\\& ~~~~~\leq \varepsilon^{-q_1-q_2} \int_{0}^{\infty}\lambda^{q_1+q_2}\left|\left\{\left(\mathbf{M}_{\frac{\alpha p q_1+\beta pq_2}{q_1+q_2}}[\mu]\right)^{\frac{1}{p-1}}>\lambda\right\}\right|\frac{d\lambda}{\lambda}.
    \end{align*}
    Since $\frac{ N(q_1+q_2)(p-1)}{2q_2(N-\beta p)}-\frac{q_1+q_2}{2}>0$, there exists $\varepsilon_0>0$ such that for any $0<\ge\leq\ge_0$, there holds $a^{-q_1-q_2}-C\varepsilon^{\frac{ N(q_1+q_2)(p-1)}{2q_2(N-\beta p)}-\frac{q_1+q_2}{2}}>0$. Hence we obtain \eqref{pot-26} by Fatou's Lemma.\smallskip
    
    \nind{\it Step 2: Proof of \eqref{pot-16}}. By \cite[Theorem 2.3]{55VHV}, we have 
      \begin{align*}
      \int_{\mathbb{R}^N}\left(\mathbf{M}^{2R}_{\frac{\alpha p q_1+\beta pq_2}{q_1+q_2}}[\omega](x)\right)^{\frac{q_1+q_2}{p-1}}dx\asymp \int_{\mathbb{R}^N}\left(\mathbf{W}^{2R}_{\frac{\alpha pq_1+\beta pq_2}{q_1+q_2},p}[\omega](x)\right)^{q_1+q_2}dx.
      \end{align*}
      Next, we prove 
      \begin{align*}
      \int_{\mathbb{R}^N}\left(\mathbf{M}^{2R}_{\frac{\alpha p q_1+\beta pq_2}{q_1+q_2}}[\omega](x)\right)^{\frac{q_1+q_2}{p-1}}dx\asymp \int_{\mathbb{R}^N}\left(\mathbf{W}^{2R}_{\alpha,p}[\omega](x)\right)^{q_1}\left(\mathbf{W}^{2R}_{\beta,p}[\omega](x)\right)^{q_2}dx.
      \end{align*}
Let $x_0\in \mathbb{R}^N$ such that $\text{supp}(\omega)\subset B_R(x_0)$.  Since  for all $x\in \mathbb{R}^N$,
     \begin{align*}
     \left(\mathbf{W}^{4R}_{\alpha,p}[\omega](x)\right)^{q_1}\left(\mathbf{W}^{4R}_{\beta,p}[\omega](x)\right)^{q_2}\geq C\left(\mathbf{M}^{2R}_{\frac{\alpha p q_1+\beta pq_2}{q_1+q_2}}[\omega](x)\right)^{\frac{q_1+q_2}{p-1}},
     \end{align*}
     and for any $y\in B_{3R/2}(x_0)$,
     \begin{align*}
     \mathbf{W}^{4R}_{\alpha,p}[\omega](y)\leq C\mathbf{W}^{2R}_{\alpha,p}[\omega](y),~~\mathbf{W}^{4R}_{\beta,p}[\omega](y)\leq C\mathbf{W}^{2R}_{\beta,p}[\omega](y),
     \end{align*}
   we have, 
     \bel{pot-27}\BA {lll}\displaystyle
     \int_{\mathbb{R}^N}\left(\mathbf{M}^{2R}_{\frac{\alpha p q+\beta pq_2}{q_1+q_2}}[\omega](x)\right)^{\frac{q_1+q_2}{p-1}}\leq C \int_{B_{5R}(x_0)}\left(\mathbf{W}^{4R}_{\alpha,p}[\omega](x)\right)^{q_1}\left(\mathbf{W}_{\beta,p}^{4R}[\omega](x)\right)^{q_2}dx
    \\[3mm]\phantom{\int_{\mathbb{R}^N}\left(\mathbf{M}^{2R}_{\frac{\alpha p q+\beta pq_2}{q_1+q_2}}[\omega](x)\right)^{\frac{q_1+q_2}{p-1}}}\displaystyle
      \leq  C \int_{B_{3R/2}(x_0)}\left(\mathbf{W}^{2R}_{\alpha,p}[\omega](x)\right)^{q_1}\left(\mathbf{W}_{\beta,p}^{2R}[\omega](x)\right)^{q_2}dx
           \\[3mm]\phantom{-------------------}\displaystyle+  C R^N \left(\frac{\omega(\mathbb{R}^N)}{R^{N-\frac{\alpha p q_1+\beta pq_2}{q_1+q_2}}}\right)^{\frac{q_1+q_2}{p-1}}
           \\[3mm]\phantom{\int_{\mathbb{R}^N}\left(\mathbf{M}^{2R}_{\frac{\alpha p q+\beta pq_2}{q_1+q_2}}[\omega](x)\right)^{\frac{q_1+q_2}{p-1}}}\displaystyle
           \leq  C \int_{\mathbb{R}^N}\left(\mathbf{W}^{2R}_{\alpha,p}[\omega](x)\right)^{q_1}\left(\mathbf{W}_{\beta,p}^{2R}[\omega](x)\right)^{q_2}dx.
\EA\ee
     On the other hand, since there holds almost everywhere,
      \begin{align*}
          \left(\mathbf{W}^{2R}_{\alpha,p}[\omega](x)\right)^{q_1}\leq C\left(\mathbf{L}^{3R}_{\alpha p, \frac{q_1}{(q_1+q_2)(p-1)}}[\omega](x)\right)^{q_1+q_2},
          \left(\mathbf{W}^{2R}_{\beta,p}[\omega](x)\right)^{q_2}\leq C\left(\mathbf{L}^{3R}_{\beta p, \frac{q_2}{(q_1+q_2)(p-1)}}[\omega](x)\right)^{q_1+q_2},
          \end{align*}
  it is enough to prove that 
     \bel{pot-28}\displaystyle
      \int_{\mathbb{R}^N}\left(\mathbf{L}^{3R}_{\alpha p, \frac{q_1}{(q_1+q_2)(p-1)}}[\omega](x)\mathbf{L}^{3R}_{\beta p, \frac{q_2}{(q_1+q_2)(p-1)}}[\omega](x)\right)^{q_1+q_2}dx\leq  C \int_{\mathbb{R}^N}\left(\mathbf{M}^{2R}_{\frac{\alpha p q_1+\beta pq_2}{q_1+q_2}}[\omega](x)\right)^{\frac{q_1+q_2}{p-1}}dx.
       \ee 
    By Lemma \eqref{lem-pot-1} there exist positive constants $C$, $\varepsilon_0$ and $a$ such that for any $\lambda>0, \varepsilon\in (0,\varepsilon_0)$,
       \begin{align*}
      & \left|\left\{\mathbf{L}^{3R}_{\alpha p, \frac{q_1}{(q_1+q_2)(p-1)}}[\omega]\mathbf{L}^{3R}_{\beta p, \frac{q_2}{(q_1+q_2)(p-1)}}[\omega]>a\lambda\right\}\right|\\&~~~~~~~~~~~\leq  C\varepsilon^{\frac{ N(q_1+q_2)(p-1)}{2q_2(N-\beta p)}}\left|\left\{\mathbf{L}^{6R}_{\alpha p, \frac{q_1}{(q_1+q_2)(p-1)}}[\omega]\mathbf{L}^{6R}_{\beta p, \frac{q_2}{(q_1+q_2)(p-1)}}[\omega]>\varepsilon^{1/2}\lambda\right\}\right|\\&~~~~~~~~~~~~~~~ +
      \left|\left\{\left(\mathbf{M}^{6R}_{\frac{\alpha p q_1+\beta pq_2}{q_1+q_2}}[\omega]\right)^{\frac{1}{p-1}}>\varepsilon\lambda\right\}\right|.
       \end{align*}
       Multiplying by $\lambda^{q_1+q_2-1}$ and integrating over $(0,\infty)$, we obtain
       \begin{align*}
       & a^{-q_1-q_2}\int_{\mathbb{R}^N}\left(\mathbf{L}^{3R}_{\alpha p, \frac{q_1}{(q_1+q_2)(p-1)}}[\omega](x)\mathbf{L}^{3R}_{\beta p, \frac{q_2}{(q_1+q_2)(p-1)}}[\omega](x)\right)^{q_1+q_2}dx
       \\&~~~\leq  C\varepsilon^{\frac{ N(q_1+q_2)(p-1)}{2q_2(N-\beta p)}-\frac{q_1+q_2}{2}}\int_{\mathbb{R}^N}\left(\mathbf{L}^{6R}_{\alpha p, \frac{q_1}{(q_1+q_2)(p-1)}}[\omega](x)\mathbf{L}^{6R}_{\beta p, \frac{q_2}{(q_1+q_2)(p-1)}}[\omega](x)\right)^{q_1+q_2}dx\\&~~~~~~~~~~~~~ +
         \varepsilon^{-q_1-q_2} \int_{\mathbb{R}^N}\left(\mathbf{M}^{6R}_{\frac{\alpha p q_1+\beta pq_2}{q_1+q_2}}[\omega](x)\right)^{\frac{q_1+q_2}{p-1}}dx.
       \end{align*}
   Similarly as \eqref{pot-27}, we  can see that 
      \begin{align*}
     & \int_{\mathbb{R}^N}\left(\mathbf{L}^{6R}_{\alpha p, \frac{q_1}{(q_1+q_2)(p-1)}}[\omega](x)\mathbf{L}^{6R}_{\beta p, \frac{q_2}{(q_1+q_2)(p-1)}}[\omega](x)\right)^{q_1+q_2}dx\\&~~~~~~~\leq C \int_{\mathbb{R}^N}\left(\mathbf{L}^{2R}_{\alpha p, \frac{q_1}{(q_1+q_2)(p-1)}}[\omega](x)\mathbf{L}^{2R}_{\beta p, \frac{q_2}{(q_1+q_2)(p-1)}}[\omega](x)\right)^{q_1+q_2}dx,
      \end{align*}
      and
     \begin{align*}
      \int_{\mathbb{R}^N}\left(\mathbf{M}^{6R}_{\frac{\alpha p q_1+\beta pq_2}{q_1+q_2}}[\omega](x)\right)^{\frac{q_1+q_2}{p-1}}dx\leq \int_{\mathbb{R}^N}\left(\mathbf{M}^{2R}_{\frac{\alpha p q_1+\beta pq_2}{q_1+q_2}}[\omega](x)\right)^{\frac{q_1+q_2}{p-1}}dx,
      \end{align*}
      Therefore, since $\frac{ N(q_1+q_2)(p-1)}{2q_2(N-\beta p)}-\frac{q_1+q_2}{2}>0$, for some $\varepsilon>0$ small enough we infer  
      \eqref{pot-28}. 
$\phantom{ppp----------}$\qeda

\begin{lemma}\label{lem-pot-4} Let $\alpha>0,p>1$, $0<\alpha p<N$ and $0<\gamma<\frac{N(p-1)}{N-\alpha p}$. There exists a constant $C=C(N,\alpha,p,\gamma)$ such that for any $\mu\in \mathfrak{M}^+(\mathbb{R}^N)$, 
   \begin{align}\label{pot-29}
   \int_{B_r(x)}\left(\mathbf{W}^r_{\alpha,p}[\mu]\right)^{\gamma}dy\leq Cr^N\left(\frac{\mu(B_{2r}(x))}{r^{N-\alpha p}}\right)^{\frac{\gamma}{p-1}} ~~\text{for all $x\in\mathbb{R}^N$ and $r>0$}. 
   \end{align}
       \end{lemma}                                                                                       
\nind\Proof We have 
   \begin{align*}
   \int_{B_r(x)}\left(\mathbf{W}^r_{\alpha,p}[\mu]\right)^{\gamma}dy& \leq \int_{B_r(x)}\left(\mathbf{W}_{\alpha,p}[\chi_{_{B_{2r}(x)}}\mu]\right)^{\gamma}dy
  \\& =\gamma \int_{0}^{\infty} \lambda^{\gamma-1}|\left\{\mathbf{W}_{\alpha,p}[\chi_{_{B_{2r}(x)}}\mu]>\lambda\right\}\cap B_r(x)|d\lambda.
   \end{align*}
   By Lemma \ref{lem-pot-2}, we obtain 
   \begin{align*} &\int_{B_r(x)}\left(\mathbf{W}^r_{\alpha,p}[\mu]\right)^{\gamma}dy\leq \gamma \left(\frac{\mu(B_{2r}(x))}{r^{N-\alpha p}}\right)^{\frac{\gamma}{p-1}}|B_r(x)|\\&~~~~~~~~~~~~~~~~~~~~~~~~+\int_{\left(\frac{\mu(B_{2r}(x))}{r^{N-\alpha p}}\right)^{\frac{1}{p-1}}}^{\infty}\lambda^{\gamma-1}|\left\{\mathbf{W}_{\alpha,p}[\chi_{_{B_{2r}(x)}}\mu]>\lambda\right\}|d\lambda
  \\& ~~~~~~\leq Cr^N\left(\frac{\mu(B_{2r}(x))}{r^{N-\alpha p}}\right)^{\frac{\gamma}{p-1}}+ C\int_{\left(\frac{\mu(B_{2r}(x))}{r^{N-\alpha p}}\right)^{\frac{1}{p-1}}}^{\infty}\lambda^{\gamma-1}\left(\frac{(\mu(B_{2r}(x)))^{\frac{1}{p-1}}}{\lambda}\right)^{\frac{N(p-1)}{N-\alpha p }}d\lambda
  \\& ~~~~~~= Cr^N\left(\frac{\mu(B_{2r}(x))}{r^{N-\alpha p}}\right)^{\frac{\gamma}{p-1}}.
   \end{align*}
   which is the claim.
   \qeda\medskip  
   
   The next result is fundamental inasmuch it shows the equivalence between the capacitary estimates,  the potential inequalities used in our construction and the solvability of the system of nonlinear integral equations connected to (\ref{In-1}).

   \begin{theorem}\label{th-pot-4} Let $\alpha,\beta,q_1,q_2>0$, $\alpha>\beta$, $1<p< \min\{N/\alpha,N/\beta\}$, $q_1+q_2>p-1$, $q_2<\frac{N(p-1)}{N-\beta p}$ and $\frac{\alpha p q_1+\beta pq_2}{q_1+q_2}<N$ and $\gm\in \mathfrak{M}^+(\mathbb{R}^N)$. Then, the following statements are  equivalent: 
   	\begin{description}
   		\item[(a)] The inequality 
   		\begin{align}\label{pot-30}
   		\mu(K)\leq C_1\text{Cap}_{\mathbf{I}_{\frac{\alpha pq_1+\beta pq_2}{q_1+q_2}},\frac{q_1+q_2}{q_1+q_2-p+1}}(K),
   		\end{align}
   		holds for any compact set $K\subset\mathbb{R}^N$, for some $C_1>0$.
   		\item[(b)] The inequality 
   		\begin{align}\label{pot-31}
   		\int_{K}\left(\mathbf{W}_{\alpha,p}[\mu](x)\right)^{q_1}\left(\mathbf{W}_{\beta,p}[\mu](x)\right)^{q_2}dx\leq C_2\text{Cap}_{\mathbf{I}_{\frac{\alpha pq_1+\beta pq_2}{q_1+q_2}},\frac{q_1+q_2}{q_1+q_2-p+1}}(K),
   		\end{align}
   		holds for any compact set $K\subset\mathbb{R}^N$, for some $C_2>0$.
   		\item[(c)] The inequality 
   		\begin{align}\label{pot-32}
   		\int_{\mathbb{R}^N}\left(\mathbf{W}_{\frac{\alpha q_1+\beta q_2}{q_1+q_2},p}[\chi_{_{B_t(x)}}\mu](y)\right)^{q_1+q_2}dy\leq C_3 \mu(B_t(x)),
   		\end{align}
   		holds for any ball $B_t(x)\subset\mathbb{R}^N$, for some $C_3>0$.
   		\item[(d)]  The inequality 
   		\begin{align}\label{pot-33}
   		\int_{\mathbb{R}^N}\left(\mathbf{W}_{\alpha,p}[\chi_{_{B_t(x)}}\mu](y)\right)^{q_1}\left(\mathbf{W}_{\beta,p}[\chi_{_{B_t(x)}}\mu](y)\right)^{q_2}dy\leq C_4 \mu(B_t(x)),
   		\end{align}
   		holds for any ball $B_t(x)\subset\mathbb{R}^N$, for some $C_4>0$.
   		\item[(e)] The inequalities 
   		\begin{align}\label{pot-34}
   		&\mathbf{W}_{\alpha,p}\left[\left(\mathbf{W}_{\alpha,p}[\mu]\right)^{q_1}\left(\mathbf{W}_{\beta,p}[\mu]\right)^{q_2}\right]\leq C_5 \mathbf{W}_{\alpha,p}[\mu]<\infty \\[2mm]
   		&  \mathbf{W}_{\beta,p}\left[\left(\mathbf{W}_{\alpha,p}[\mu]\right)^{q_1}\left(\mathbf{W}_{\beta,p}[\mu]\right)^{q_2}\right]\leq C_5 \mathbf{W}_{\beta,p}[\mu]<\infty \label{pot-35} 
   		\end{align}
   		hold for some $C_5>0$.
   		\item[(f)] The system equation  
   		\begin{equation}\label{pot-36}
   		\begin{array}
   		[c]{l}%
   		U=\mathbf{W}_{\alpha,p}\left[U^{q_1}V^{q_2}\right] + \varepsilon\mathbf{W}_{\alpha,p}\left[\mu\right]\\[2mm]
   		V=\mathbf{W}_{\beta,p}\left[U^{q_1}V^{q_2}\right] +\varepsilon  \mathbf{W}_{\beta,p}\left[\mu\right] ,  
   		\end{array}
   		\end{equation}
   		in $\mathbb{R}^N$ has a nonnegative solution for some $\varepsilon>0$. 
   	\end{description}
            \end{theorem}
       
    \nind\Proof
   By Theorem \ref{th-pot-1} we have (a) $ \Leftrightarrow $ (c), by Theorem \ref{th-pot-3}, (c) $ \Leftrightarrow $ (d).  We now assume $(e)$. Put $\mathbf{T}[\gm]=\left(\mathbf{W}_{\alpha,p}[\mu](x)\right)^{q_1}\left(\mathbf{W}_{\beta,p}[\mu](x)\right)^{q_2}$ for any $\gm\in \mathfrak{M}^+(\mathbb{R}^N)$.
   It is easy to see that 
   \begin{align}\label{pot-13}
   \left(\mathbf{T}[\gm](x)\right)^{\gamma}\geq C \int_{0}^{\infty}\left(\frac{\gm(B_\rho(x))}{\rho^{N-\frac{\alpha q_1p+\beta q_2 p}{q_1+q_2}}}\right)^{\frac{\gamma(q_1+q_2)}{p-1}}\frac{d\rho}{\rho}=C\mathbf{W}_{\beta,s}[\mu](x)~~\text{for all} ~x\in\mathbb{R}^N
   \end{align}
   where $\gamma=\frac{p-1}{q_1}+\frac{p-1}{q_2}, \beta=\frac{\gamma(\alpha q_1p+\beta q_2p)}{\gamma(q_1+q_2)+p-1}$ and $s=\frac{\gamma(q_1+q_2)+p-1}{\gamma(q_1+q_2)}<1+\frac{1}{\gamma}$.  From \eqref{pot-34} and \eqref{pot-35}, we have 
   \begin{align*}
   \mathbf{T}\left[\mathbf{T}[\mu]\right]\leq C\mathbf{T}[\mu]<\infty ~~\text{almost everywhere}.
   \end{align*}
   Using \eqref{pot-13}, we obtain 
   \begin{align*}
   \left(\mathbf{W}_{\beta,s}\left[\mathbf{T}[\mu]\right]\right)^{\frac{1}{\gamma}}\leq C\mathbf{T}[\mu]<\infty ~~\text{almost everywhere}. 
   \end{align*}
   Applying $\mathbf{W}_{\beta,s}$ to both sides of the above inequality and using Theorem \ref{th-pot-1} with $\alpha=\beta,p=s,q=\frac{1}{\gamma}$ , we derive 
   \begin{align}
   \int_K\mathbf{T}[\mu](x)dx\leq C\text{Cap}_{\mathbf{I}_{\beta s},\frac{1}{1+\gamma-\gamma s}}(K),
   \end{align}
   for any compact set $K\subset\mathbb{R}^N$, which implies (b).  So, (e) $ \Rightarrow  $ (b). Next, assume (b), using  \eqref{pot-13} again, we derive from (b) that 
   \begin{align*}
   \int_K\left(\mathbf{W}_{\beta,s}[\mu](x)\right)^{\frac{1}{\gamma}}dx\leq C\text{Cap}_{\mathbf{I}_{\beta s},\frac{1}{1+\gamma-\gamma s}}(K),
   \end{align*}
   for any compact set $K\subset\mathbb{R}^N$. Thanks to Theorem \ref{th-pot-1}, we get (a). So,  (b) $ \Rightarrow  $ (a).
        
   \nind It remains to prove that {\it(i)}: (f) $ \Rightarrow  $ (a), {\it(ii)}:  (e) $ \Rightarrow  $ (f), \textbf{\it (iii)}: (a)+(c)+(d) $ \Rightarrow  $ (e).\smallskip
    
\nind  {\it(i)}: Assume that \eqref{pot-36}  has a nonnegative solution for some $\varepsilon>0$.  Set $d\nu(x)=U^{q_1}V^{q_2}dx+\varepsilon d\mu(x)$. Clearly
               \begin{align*}
               \left(\mathbf{W}_{\alpha,p}[\nu]\right)^{q_1}\left(\mathbf{W}_{\beta,p}[\nu]\right)^{q_2}\leq Cd\nu(x)~~\text{ in }~\mathbb{R}^N. 
               \end{align*}              
               If $E\subset\mathbb{R}^N$ is a Borel set, we have
               \begin{align*}
            \int_{\mathbb{R}^N}\left(\mathbf{W}_{\alpha,p}[\chi_{_E}\nu]\right)^{q_1}\left(\mathbf{W}_{\beta,p}[\chi_{_E}\nu]\right)^{q_2}dx   &\leq  \int_{\mathbb{R}^N}(\mathbf{M}_{\nu}\chi_{_E})^{\frac{q_1+q_2}{p-1}}\left(\mathbf{W}_{\alpha,p}[\nu]\right)^{q_1}\left(\mathbf{W}_{\beta,p}[\nu]\right)^{q_2}dx\\&\leq C\int_{\mathbb{R}^N}(\mathbf{M}_{\nu}\chi_{_E})^{\frac{q_1+q_2}{p-1}}d\nu.
               \end{align*}
               Since $M_\omega f$ is bounded on $L^s(\mathbb{R}^N,d\omega)$, $s>1$, we deduce from Fefferman's result \cite{22Fe} that
              \begin{align*}
                \int_{\mathbb{R}^N}\left(\mathbf{W}_{\alpha,p}[\chi_{_E}\nu]\right)^{q_1}\left(\mathbf{W}_{\beta,p}[\chi_{_E}\nu]\right)^{q_2}dx\leq C\nu(E),
                \end{align*}
            is verified     for any Borel set $E\subset\mathbb{R}^N$. Applying (a) $ \Leftrightarrow $ (c) to $\mu=\nu$, we derive that
                  \begin{align}\label{pot-37}
                     \nu(K)\leq C_1\text{Cap}_{\mathbf{I}_{\frac{\alpha pq_1+\beta pq_2}{q_1+q_2}},\frac{q_1+q_2}{q_1+q_2-p+1}}(K),
                     \end{align}
                     holds for any compact set $K\subset\mathbb{R}^N$. Since $\nu\geq \mu$, we obtain  (c). \smallskip
                     
 \nind {\it(ii)}: Suppose that \eqref{pot-34} and \eqref{pot-35} hold with constant $C_5>0$.Take $0<\varepsilon\leq \frac{1}{2(2C_5)^{\frac{p-1}{q_1+q_2-p+1}}}$.  Consider  the sequence $\{U_m,V_m\}_{m\geq 0}$ of nonnegative functions defined by $U_0=\mathbf{W}_{\alpha,p}[\mu], V_0=\mathbf{W}_{\beta,p}[\mu]$ and 
  \begin{equation*}
 \begin{array}
[c]{l}%
U_{m+1}=\mathbf{W}_{\alpha,p}\left[U_m^{q_1}V_m^{q_2}\right] + \varepsilon\mathbf{W}_{\alpha,p}\left[\mu\right]\\[2mm]
V_{m+1}=\mathbf{W}_{\beta,p}\left[U_m^{q_1}V_m^{q_2}\right] +\varepsilon  \mathbf{W}_{\beta,p}\left[\mu\right].                                                                                          
\end{array}
\end{equation*}
                   It is easy to see that $\{U_m,V_m\}_{m\geq 0}$ is well defined and satisfies
                   \begin{align*}
                   U_{m}\leq 2\varepsilon \mathbf{W}_{\alpha,p}[\mu], V_{m}\leq 2\varepsilon \mathbf{W}_{\alpha,p}[\mu]~~\text{for all }~m\geq 0. 
                   \end{align*}
Clearly $\{U_m\}, \{V_m\}$ are nondecreasing. Using the dominated convergence theorem, it follows that  $ (U(x),V(x)) := \mathop {\lim }\limits_{m \to \infty } (U_m(x),V_m(x))$
 is  a solution of \eqref{pot-36}.\smallskip
 
  \nind{\it(iii)}: Assume that statements (a), (c) and (d) hold true.  We first assume that $\mu$ has compact support. 
   From (a) we have
   \begin{align}\label{pot-38}
    \mu(B_r(x))\leq Cr^{N-\frac{\alpha pq_1+\beta pq_2}{q_1+q_2-p+1}}~~\text{for all }x\in \BBR^N \text{ and }r>0.
    \end{align}
    From (b)
     \begin{align*}
      \int_{B_r(x)}\left(\mathbf{W}^{r}_{\frac{\alpha q_1+\beta q_2}{q_1+q_2},p}[\mu](y)\right)^{q_1+q_2}dy\leq C_2 \mu(B_{2r}(x)) ~~\text{for all }x\in \BBR^N \text{ and }r>0.
      \end{align*} 
    Using H\"older's inequality and $\mathbf{W}^{r}_{\frac{\alpha q_1+\beta q_2}{q_1+q_2},p}[\mu]\geq r^{-\frac{(\alpha-\beta)pq_2}{(p-1)(q_1+q_2)}}\mathbf{W}^{r}_{\alpha ,p}[\mu]$, we obtain,  
    \begin{align}\label{pot-39}
       \int_{B_r(x)}\left(\mathbf{W}^{r}_{\alpha,p}[\mu](y)\right)^{q_1}dy\leq C r^{\frac{(\alpha-\beta)pq_1q_2+(p-1)Nq_2+(N-\beta p)(p-1)q_1}{(p-1)(q_1+q_2)}}\left(\frac{\mu(B_{2r}(x))}{r^{N-\beta p}}\right)^{\frac{q_1}{q_1+q_2}}, 
       \end{align}
 again for all $x\in\BBR^N$ and $r>0$.      From (c), 
       \begin{align}\label{pot-40}
         \int_{B_r(x)}\left(\mathbf{W}^r_{\alpha,p}[\mu](y)\right)^{q_1}\left(\mathbf{W}^r_{\beta,p}[\mu](y)\right)^{q_2}dy\leq C_3 \mu(B_{2r}(x)) ~~\text{for all }x\in \BBR^N \text{ and }r>0.
         \end{align}
   By Lemma \ref{lem-pot-4}, 
   \begin{align}\label{pot-41}
    \int_{B_r(x)}\left(\mathbf{W}^r_{\beta,p}[\mu]\right)^{q_2}dy\leq Cr^{N}\left(\frac{\mu(B_{2r}(x))}{r^{N-\beta p}}\right)^{\frac{q_2}{p-1}} ~~\text{for all }x\in \BBR^N \text{ and }r>0. 
    \end{align}
   We have, with $\eta=\alpha$ or $\eta=\beta$,
   \begin{align}\label{pot-42}
   \mathbf{W}_{\eta,p}\left[\left(\mathbf{W}_{\alpha,p}[\mu]\right)^{q_1}\left(\mathbf{W}_{\beta,p}[\mu]\right)^{q_2}\right](x)\leq C\sum_{i=1}^{4}\int_{0}^{\infty}\left(\frac{\mathbf{A}_i(x,r)}{r^{N-\eta p}}\right)^{\frac{1}{p-1}}\frac{dr}{r},
   \end{align}
   where 
   \begin{align*}
   &\mathbf{A}_1(x,r)=\int_{B_r(x)}\left(\mathbf{W}^r_{\alpha,p}[\mu](y)\right)^{q_1}\left(\mathbf{W}^r_{\beta,p}[\mu](y)\right)^{q_2}dy,\\&
   \mathbf{A}_2(x,r)=\int_{B_r(x)}\left(\mathbf{W}^r_{\alpha,p}[\mu](y)\right)^{q_1}\left(\int_{r}^{\infty}\left(\frac{\mu(B_t(y))}{t^{N-\beta p}}\right)^{\frac{1}{p-1}}\frac{dt}{t}\right)^{q_2}dy,
   \\&
   \mathbf{A}_3(x,r)=\int_{B_r(x)}\left(\int_{r}^{\infty}\left(\frac{\mu(B_t(y))}{t^{N-\alpha p}}\right)^{\frac{1}{p-1}}\frac{dt}{t}\right)^{q_1}\left(\mathbf{W}^r_{\beta,p}[\mu](y)\right)^{q_2}dy,
   \\&
   \mathbf{A}_4(x,r)=\int_{B_r(x)}\left(\int_{r}^{\infty}\left(\frac{\mu(B_t(y))}{t^{N-\alpha p}}\right)^{\frac{1}{p-1}}\frac{dt}{t}\right)^{q_1}\left(\int_{r}^{\infty}\left(\frac{\mu(B_t(y))}{t^{N-\beta p}}\right)^{\frac{1}{p-1}}\frac{dt}{t}\right)^{q_2}dy.
   \end{align*}
   Thanks to   \eqref{pot-40} we get
   \begin{align*}
   \mathbf{A}_1(x,r)\leq C\mu(B_{2r}(x)),
   \end{align*}
   which implies 
   \begin{align}\label{pot-43}
   \int_{0}^{\infty}\left(\frac{\mathbf{A}_1(x,r)}{r^{N-\eta p}}\right)^{\frac{1}{p-1}}\frac{dr}{r}\leq C \mathbf{W}_{\eta,p}[\mu](x).
   \end{align}
   Since $B_t(y)\leq B_{2t}(x)$ for any $y\in B_r(x)$, $t\geq r$  and thanks to \eqref{pot-39}, \eqref{lem-pot-4} we deduce
   \begin{align*}
   &\mathbf{A}_2(x,t)\leq \int_{B_r(x)}\left(\mathbf{W}^r_{\alpha,p}[\mu](y)\right)^{q_1}dy\left(\int_{r}^{\infty}\left(\frac{\mu(B_{2t}(x))}{t^{N-\beta p}}\right)^{\frac{1}{p-1}}\frac{dt}{t}\right)^{q_2}\\&~\leq C r^{\frac{(\alpha-\beta)pq_1q_2+(p-1)Nq_2+(N-\beta p)(p-1)q_1}{(p-1)(q_1+q_2)}}\left(\frac{\mu(B_{2r}(x))}{r^{N-\beta p}}\right)^{\frac{q_1}{q_1+q_2}} \left(\int_{r}^{\infty}\left(\frac{\mu(B_{2t}(x))}{t^{N-\beta p}}\right)^{\frac{1}{p-1}}\frac{dt}{t}\right)^{q_2}
   \\&~ \leq C r^{\frac{(\alpha-\beta)pq_1q_2+(p-1)Nq_2+(N-\beta p)(p-1)q_1}{(p-1)(q_1+q_2)}} \left(\int_{r}^{\infty}\left(\frac{\mu(B_{2t}(x))}{t^{N-\beta p}}\right)^{\frac{1}{p-1}}\frac{dt}{t}\right)^{q_2+\frac{q_1(p-1)}{q_1+q_2}},
   \end{align*}
   then
   \begin{align*}
 &\mathbf{A}_3(x,t)\leq \left(\int_{r}^{\infty}\left(\frac{\mu(B_{2t}(x))}{t^{N-\alpha p}}\right)^{\frac{1}{p-1}}\frac{dt}{t}\right)^{q_1}\int_{B_r(x)}\left(\mathbf{W}^r_{\beta,p}[\mu](y)\right)^{q_2}dy\\
&\phantom{ \mathbf{A}_3(x,t)}\leq C \left(\int_{r}^{\infty}\left(\frac{\mu(B_{2t}(x))}{t^{N-\alpha p}}\right)^{\frac{1}{p-1}}\frac{dt}{t}\right)^{q_1}r^{N}\left(\frac{\mu(B_{2r}(x))}{r^{N-\beta p}}\right)^{\frac{q_2}{p-1}}
   \\
   &\phantom{ \mathbf{A}_3(x,t)}
   \leq C r^{N} \left(\int_{r}^{\infty}\left(\frac{\mu(B_{2t}(x))}{t^{N-\alpha p}}\right)^{\frac{1}{p-1}}\frac{dt}{t}\right)^{q_1}\left(\int_{r}^{\infty}\left(\frac{\mu(B_{2t}(x))}{t^{N-\beta p}}\right)^{\frac{1}{p-1}}\frac{dt}{t}\right)^{q_2},
   \end{align*}
   and finally 
   \begin{align*}
   &\mathbf{A}_4(x,t)\leq C r^{N} \left(\int_{r}^{\infty}\left(\frac{\mu(B_{2t}(x))}{t^{N-\alpha p}}\right)^{\frac{1}{p-1}}\frac{dt}{t}\right)^{q_1}\left(\int_{r}^{\infty}\left(\frac{\mu(B_{2t}(x))}{t^{N-\beta p}}\right)^{\frac{1}{p-1}}\frac{dt}{t}\right)^{q_2}.
   \end{align*}
   
 \nind \textbf{I-} From the estimate of $\mathbf{A}_2$ we derive 
   \begin{align*}
    &\int_{0}^{\infty}\left(\frac{\mathbf{A}_2(x,r)}{r^{N-\eta p}}\right)^{\frac{1}{p-1}}\frac{dr}{r}\leq C \int_{0}^{\infty}r^{\frac{(\alpha-\beta)pq_1q_2+(p-1)Nq_2+(N-\beta p)(p-1)q_1}{(p-1)^2(q_1+q_2)}-\frac{N-\eta p}{p-1}}\\&~~~~~~~~~~~~~~~~~~~~~~~~~~~~~~~~\times\left(\int_{r}^{\infty}\left(\frac{\mu(B_{2t}(x))}{t^{N-\beta p}}\right)^{\frac{1}{p-1}}\frac{dt}{t}\right)^{\frac{q_2}{p-1}+\frac{q_1}{q_1+q_2}}\frac{dr}{r}.
   \end{align*}
Since  $\frac{\alpha p q_1+\beta pq_2}{q_1+q_2}<N$, it follows that
   \begin{align*}
   0&<\kappa:=\frac{(\alpha-\beta)pq_1q_2+(p-1)Nq_2+(N-\beta p)(p-1)q_1}{(p-1)^2(q_1+q_2)}-\frac{N-\eta p}{p-1}
   \\&~~~~~~~<\frac{N-\beta p}{p-1}\left(\frac{q_2}{p-1}+\frac{q_1}{q_1+q_2}\right).
   \end{align*}
   Hence, 
   \begin{align*}
   r^{\kappa}\left(\int_{r}^{\infty}\left(\frac{\mu(B_{2t}(x))}{t^{N-\beta p}}\right)^{\frac{1}{p-1}}\frac{dt}{t}\right)^{\frac{q_2}{p-1}+\frac{q_1}{q_1+q_2}}\to 0 ~\text{ as}~t\to 0,
   \end{align*}
   and therefore
   \begin{align*}
     & r^{\kappa}\left(\int_{r}^{\infty}\left(\frac{\mu(B_{2t}(x))}{t^{N-\beta p}}\right)^{\frac{1}{p-1}}\frac{dt}{t}\right)^{\frac{q_2}{p-1}+\frac{q_1}{q_1+q_2}}\\& ~~~~~~~~~~~~~~~~~~~~~~\leq C
     r^{\kappa-\frac{N-\beta p}{p-1}\left(\frac{q_2}{p-1}+\frac{q_1}{q_1+q_2}\right)}\left(\mu(\mathbb{R}^N)\right)^{\frac{q_2}{(p-1)^2}+\frac{q_1}{(p-1)(q_1+q_2)}},
 \end{align*}
  a quantity which converges to $0$ when $t\to 0$.  Hence, by integration be parts,
    we obtain 
   \begin{align*}
   &\int_{0}^{\infty}\left(\frac{\mathbf{A}_2(x,r)}{r^{N-\eta p}}\right)^{\frac{1}{p-1}}\frac{dr}{r}\leq C \int_{0}^{\infty}r^{\kappa}\left(\int_{r}^{\infty}\left(\frac{\mu(B_{2t}(x))}{t^{N-\beta p}}\right)^{\frac{1}{p-1}}\frac{dt}{t}\right)^{\frac{q_2}{p-1}+\frac{q_1}{q_1+q_2}}\frac{dr}{r}
   \\&~~~~~~~~~\leq C\int_{0}^{\infty}r^{\kappa}\left(\int_{r}^{\infty}\left(\frac{\mu(B_{2t}(x))}{t^{N-\beta p}}\right)^{\frac{1}{p-1}}\frac{dt}{t}\right)^{\frac{q_2}{p-1}+\frac{q_1}{q_1+q_2}-1}\left(\frac{\mu(B_{2r}(x))}{r^{N-\beta p}}\right)^{\frac{1}{p-1}}\frac{dr}{r}
   \\&~~~~~~~~~= C\int_{0}^{\infty}r^{\frac{(\alpha-\beta)pq_1q_2+(p-1)Nq_2+(N-\beta p)(p-1)q_1}{(p-1)^2(q_1+q_2)}-\frac{N-\beta p}{p-1}}\\&~~~~~~~~~~~~~~~~~\times\left(\int_{r}^{\infty}\left(\frac{\mu(B_{2t}(x))}{t^{N-\beta p}}\right)^{\frac{1}{p-1}}\frac{dt}{t}\right)^{\frac{q_2}{p-1}+\frac{q_1}{q_1+q_2}-1}\left(\frac{\mu(B_{2r}(x))}{r^{N-\eta p}}\right)^{\frac{1}{p-1}}\frac{dr}{r}.
   \end{align*} 
   Observing that we have  from \eqref{pot-38}, 
   \begin{align*}
   r^{\frac{(\alpha-\beta)pq_1q_2+(p-1)Nq_2+(N-\beta p)(p-1)q_1}{(p-1)^2(q_1+q_2)}-\frac{N-\beta p}{p-1}}\left(\int_{r}^{\infty}\left(\frac{\mu(B_{2t}(x))}{t^{N-\beta p}}\right)^{\frac{1}{p-1}}\frac{dt}{t}\right)^{\frac{q_2}{p-1}+\frac{q_1}{q_1+q_2}-1}\leq C,
   \end{align*}
we derive
   \begin{align}\label{pot-44}
   \int_{0}^{\infty}\left(\frac{\mathbf{A}_2(x,r)}{r^{N-\eta p}}\right)^{\frac{1}{p-1}}\frac{dr}{r}\leq C \mathbf{W}_{\eta,p}[\mu](x).
   \end{align}

 \nind \textbf{II-} From the estimate of $\mathbf{A}_3$ and $\mathbf{A}_4$, we have, as above,  by integration be parts,
   \begin{align*}
  & \int_{0}^{\infty}\left(\frac{\mathbf{A}_3(x,r)}{r^{N-\eta p}}\right)^{\frac{1}{p-1}}\frac{dr}{r}+\int_{0}^{\infty}\left(\frac{\mathbf{A}_4(x,r)}{r^{N-\eta p}}\right)^{\frac{1}{p-1}}\frac{dr}{r}\\& ~~~\leq C \int_{0}^{\infty}r^{\frac{\eta p}{p-1}}\left(\int_{r}^{\infty}\left(\frac{\mu(B_{2t}(x))}{t^{N-\alpha p}}\right)^{\frac{1}{p-1}}\frac{dt}{t}\right)^{\frac{q_1}{p-1}}\left(\int_{r}^{\infty}\left(\frac{\mu(B_{2t}(x))}{t^{N-\beta p}}\right)^{\frac{1}{p-1}}\frac{dt}{t}\right)^{\frac{q_2}{p-1}}\frac{dr}{r}
  \\& ~~~=C \int_{0}^{\infty}\mathbf{D}_1(x,r)\left(\frac{\mu(B_{2r}(x))}{r^{N-\eta p}}\right)^{\frac{1}{p-1}}\frac{dr}{r}+C \int_{0}^{\infty}\mathbf{D}_2(x,r)\left(\frac{\mu(B_{2r}(x))}{r^{N-\eta p}}\right)^{\frac{1}{p-1}}\frac{dr}{r},
   \end{align*}
   where
   \begin{align*}
 & \mathbf{D}_1(x,r)=r^{\frac{\alpha p}{p-1}}\left(\int_{r}^{\infty}\left(\frac{\mu(B_{2t}(x))}{t^{N-\alpha p}}\right)^{\frac{1}{p-1}}\frac{dt}{t}\right)^{\frac{q_1}{p-1}-1}\left(\int_{r}^{\infty}\left(\frac{\mu(B_{2t}(x))}{t^{N-\beta p}}\right)^{\frac{1}{p-1}}\frac{dt}{t}\right)^{\frac{q_2}{p-1}},\\
  &\mathbf{D}_2(x,r)=r^{\frac{\beta p}{p-1}}\left(\int_{r}^{\infty}\left(\frac{\mu(B_{2t}(x))}{t^{N-\alpha p}}\right)^{\frac{1}{p-1}}\frac{dt}{t}\right)^{\frac{q_1}{p-1}}\left(\int_{r}^{\infty}\left(\frac{\mu(B_{2t}(x))}{t^{N-\beta p}}\right)^{\frac{1}{p-1}}\frac{dt}{t}\right)^{\frac{q_2}{p-1}-1}.
   \end{align*}
Clearly,
   \begin{align*}
  & \int_{r}^{\infty}\left(\frac{\mu(B_{2t}(x))}{t^{N-\beta p}}\right)^{\frac{1}{p-1}}\frac{dt}{t}\leq r^{-\frac{(\alpha-\beta)p}{p-1}}\int_{r}^{\infty}\left(\frac{\mu(B_{2t}(x))}{t^{N-\alpha p}}\right)^{\frac{1}{p-1}}\frac{dt}{t},\\
  &\int_{r}^{\infty}\left(\frac{\mu(B_{2t}(x))}{t^{N-\alpha p}}\right)^{\frac{1}{p-1}}\frac{dt}{t}\leq C r^{\frac{\alpha p}{p-1}-\frac{\alpha pq_1+\beta pq_2}{(p-1)(q_1+q_2-p+1)}}.
   \end{align*}
   We derive 
   \begin{align*}
   \mathbf{D}_1(x,r)&\leq C r^{\frac{\alpha p}{p-1}}\left(\int_{r}^{\infty}\left(\frac{\mu(B_{2t}(x))}{t^{N-\alpha p}}\right)^{\frac{1}{p-1}}\frac{dt}{t}\right)^{\frac{q_1}{p-1}-1}\left(r^{-\frac{(\alpha-\beta)p}{p-1}}\int_{r}^{\infty}\left(\frac{\mu(B_{2t}(x))}{t^{N-\alpha p}}\right)^{\frac{1}{p-1}}\frac{dt}{t}\right)^{\frac{q_2}{p-1}}
  \\&= C_1 r^{\frac{\alpha p}{p-1}-\frac{(\alpha-\beta)pq_2}{(p-1)^2}}\left(\int_{r}^{\infty}\left(\frac{\mu(B_{2t}(x))}{t^{N-\alpha p}}\right)^{\frac{1}{p-1}}\frac{dt}{t}\right)^{\frac{q_1+q_2}{p-1}-1}
  \\& \leq  C_2 r^{\frac{\alpha p}{p-1}-\frac{(\alpha-\beta)pq_2}{(p-1)^2}}\left(r^{\frac{\alpha p}{p-1}-\frac{\alpha pq_1+\beta pq_2}{(p-1)(q_1+q_2-p+1)}}\right)^{\frac{q_1+q_2}{p-1}-1}=C_3.
   \end{align*}
   Next, we estimate $\mathbf{D}_2(x,r)$. If $\frac{q_1}{p-1}\geq 1$, similarly as for estimate of $\mathbf{D}_1(x,r)$ we obtain $\mathbf{D}_2(x,r)\leq C$. If   $\frac{q_1}{p-1}<1$, we have
   \begin{align*}
   \mathbf{D}_2(x,r)&=\frac{q_1}{p-1}r^{\frac{\beta p}{p-1}}\int_{r}^{\infty}\left(\frac{\mu(B_{2t}(x))}{t^{N-\alpha p}}\right)^{\frac{1}{p-1}}\left(\int_{t}^{\infty}\left(\frac{\mu(B_{2s}(x))}{s^{N-\alpha p}}\right)^{\frac{1}{p-1}}\frac{ds}{s}\right)^{\frac{q_1}{p-1}-1}\frac{dt}{t}\\
   & \phantom{--------------}\times\left(\int_{r}^{\infty}\left(\frac{\mu(B_{2t}(x))}{t^{N-\beta p}}\right)^{\frac{1}{p-1}}\frac{dt}{t}\right)^{\frac{q_2}{p-1}-1}
   \\[2mm]& \leq r^{\frac{\beta p}{p-1}}\int_{r}^{\infty}\left(\frac{\mu(B_{2t}(x))}{t^{N-\alpha p}}\right)^{\frac{1}{p-1}}\left(\int_{t}^{\infty}\left(\frac{\mu(B_{2s}(x))}{s^{N-\alpha p}}\right)^{\frac{1}{p-1}}\frac{ds}{s}\right)^{\frac{q_1}{p-1}-1}\\
   & \phantom{--------------}\times\left(\int_{t}^{\infty}\left(\frac{\mu(B_{2s}(x))}{s^{N-\beta p}}\right)^{\frac{1}{p-1}}\frac{ds}{s}\right)^{\frac{q_2}{p-1}-1}\frac{dt}{t}.
   \end{align*}
   On the other hand, 
   \begin{align*}
   &\left(\frac{\mu(B_{2t}(x))}{t^{N-\alpha p}}\right)^{\frac{1}{p-1}}=t^{\frac{(\alpha-\beta)p(p-1-q_2)}{(p-1)^2}}\left(\left(\frac{\mu(B_{2t}(x))}{t^{N-\alpha p}}\right)^{\frac{1}{p-1}}\right)^{\frac{q_2}{p-1}}\left(\left(\frac{\mu(B_{2t}(x))}{t^{N-\beta p}}\right)^{\frac{1}{p-1}}\right)^{1-\frac{q_2}{p-1}}
 \\& ~~\leq C  t^{\frac{(\alpha-\beta)p(p-1-q_2)}{(p-1)^2}}\left(\int_{t}^{\infty}\left(\frac{\mu(B_{2s}(x))}{s^{N-\alpha p}}\right)^{\frac{1}{p-1}}\frac{ds}{s}\right)^{\frac{q_2}{p-1}}\left(\int_{t}^{\infty}\left(\frac{\mu(B_{2s}(x))}{s^{N-\beta p}}\right)^{\frac{1}{p-1}}\frac{ds}{s}\right)^{1-\frac{q_2}{p-1}}\!\!\!\!\!\!,
   \end{align*}
   therefore, 
   \begin{align*}
   \mathbf{D}_2(x,r)&\leq C  r^{\frac{\beta p}{p-1}}\int_{r}^{\infty}t^{\frac{(\alpha-\beta)p(p-1-q_2)}{(p-1)^2}}\left(\int_{t}^{\infty}\left(\frac{\mu(B_{2s}(x))}{s^{N-\alpha p}}\right)^{\frac{1}{p-1}}\frac{ds}{s}\right)^{\frac{q_1+q_2}{p-1}-1}\frac{dt}{t}
  \\& \leq  C_1  r^{\frac{\beta p}{p-1}}\int_{r}^{\infty}t^{\frac{(\alpha-\beta)p(p-1-q_2)}{(p-1)^2}}\left(r^{\frac{\alpha p}{p-1}-\frac{\alpha pq_1+\beta pq_2}{(p-1)(q_1+q_2-p+1)}}\right)^{\frac{q_1+q_2}{p-1}-1}\frac{dt}{t}
  \\& = C_2.
   \end{align*}
   Hence, 
   \begin{align}\label{pot-45}
   \int_{0}^{\infty}\left(\frac{\mathbf{A}_3(x,r)}{r^{N-\eta p}}\right)^{\frac{1}{p-1}}\frac{dr}{r}+\int_{0}^{\infty}\left(\frac{\mathbf{A}_4(x,r)}{r^{N-\eta p}}\right)^{\frac{1}{p-1}}\frac{dr}{r}\leq C \mathbf{W}_{\eta,p}[\mu](x). 
   \end{align}
   Combining \eqref{pot-42} with  \eqref{pot-43}, \eqref{pot-44} and \eqref{pot-45}  we obtain \begin{align*}
   \mathbf{W}_{\eta,p}\left[\left(\mathbf{W}_{\alpha,p}[\mu]\right)^{q_1}\left(\mathbf{W}_{\beta,p}[\mu]\right)^{q_2}\right]\leq C \mathbf{W}_{\eta,p}[\mu]<\infty,
   \end{align*}
for $\eta=\alpha$ or $\beta$, provided $\mu$ has compact support in $\mathbb{R}^N.$
Next, we assume that $\mu$ may not have compact support. Since the above constants noted $C$ are independent of $\gm$, for $n\in\BBN^*$, we set 
$\gm_n=\chi_{_{B_n(0)}}\mu$
    \begin{align*}
       \mathbf{W}_{\eta,p}\left[\left(\mathbf{W}_{\alpha,p}[\gm_n]\right)^{q_1}\left(\mathbf{W}_{\beta,p}[\gm_n]\right)^{q_2}\right]\leq C \mathbf{W}_{\eta,p}[\gm_n]\leq C \mathbf{W}_{\eta,p}[\mu]<\infty<C',
       \end{align*}
 for $\eta=\alpha$ or $\beta$. Then we infer (e) by Fatou's lemma.
    \qeda\medskip
    
    An important step for proving relative compactness in nonlinear problems is the convergence of the nonlinear terms and their equi-integrability is one of the 
    key tool for such a task. 
     
    \begin{lemma} \label{equi-integ-condi}Let $\mu$ be satisfying \eqref{pot-30} with compact support in $\mathbb{R}^N$. Set $\mu_n=\varphi_n\star\mu$. Then, 
    	\begin{align}
    	\left(\mathbf{W}_{\alpha,p}[\mu_n](x)\right)^{q_1}\left(\mathbf{W}_{\beta,p}[\mu_n](x)\right)^{q_2}
    	\end{align}
    	 is equi-integrable in  $B_t(0)$ for all $t>1.$
    \end{lemma}
\Proof Since $\supp \mu_n\subset B_{t_0}(0)$ for some $t_0>0$ and 
$$
\left(\mathbf{W}^{2T}_{\alpha,p}[\mu_n]\right)^{q_1}\leq C\left(\mathbf{L}^{3T}_{\alpha p, \frac{q_1}{(q_1+q_2)(p-1)}}[\mu_n]\right)^{q_1+q_2},
	\left(\mathbf{W}^{2T}_{\beta,p}[\mu_n]\right)^{q_2}\leq C\left(\mathbf{L}^{3T}_{\beta p, \frac{q_2}{(q_1+q_2)(p-1)}}[\mu_n]\right)^{q_1+q_2},$$
	 it suffices to show that $$ \left(\mathbf{L}^{2(t_0+t)}_{\alpha p, \frac{q_1}{(q_1+q_2)(p-1)}}[\mu_n]\mathbf{L}^{2(t_0+t)}_{\beta p, \frac{q_2}{(q_1+q_2)(p-1)}}[\mu_n]\right)^{q_1+q_2}	$$ is equi-integrable in  $B_{t}(0)$. Since $\left(\mathbf{I}_{\frac{\alpha pq_1+p\beta q_2}{q_1+q_2}}^{2(t_0+t)}[\mu_n]\right)^{\frac{q_1+q_2}{p-1}}\leq C \left(\mathbf{I}_{\frac{\alpha pq_1+p\beta q_2}{q_1+q_2}}^{4(t_0+t)}[\mu]\star\varphi_n\right)^{\frac{q_1+q_2}{p-1}}$, so $\left(\mathbf{I}_{\frac{\alpha pq_1+p\beta q_2}{q_1+q_2}}^{2(t_0+t)}[\mu_n]\right)^{\frac{q_1+q_2}{p-1}}$ is equi-integrable in  $B_t(0)$ for any $t>t_0$.
	Thus, by \cite[Proposition 1.27]{luigi} we can find a nondecreasing function $\Phi: [0,\infty)\to [0,\infty)$ such that  $\Phi(\lambda)/\lambda\to\infty$ as  $\lambda\to \infty$, and $\phi(2^j\lambda)\leq j\phi(\lambda)$ for all $\lambda>0,j\in\mathbb{N}$ and  $\Phi'(\lambda)=\phi(\lambda)$
	\begin{align*}
	\int_{0}^{\infty}\phi(\lambda)\left|\left\{\left(\mathbf{I}_{\frac{\alpha pq_1+p\beta q_2}{q_1+q_2}}^{2(t_0+t)}[\mu_n]\right)^{\frac{q_1+q_2}{p-1}}>\lambda\right\}\right|d\lambda\leq 1.
	\end{align*} 
	On the other hand, by Lemma \ref{lem-pot-1}, there exists $C>0$ and $\varepsilon_0>0$ such that 
	\begin{align}\label{es7}
	\nonumber	&\left|\left\{\mathbf{L}^{2(t_0+t)}_{\alpha p, \frac{q_1}{(q_1+q_2)(p-1)}}[\mu_n]\mathbf{L}^{2(t_0+t)}_{\beta p, \frac{q_2}{(q_1+q_2)(p-1)}}[\mu_n]>a\lambda, \left(\mathbf{I}_{\frac{\alpha pq_1+p\beta q_2}{q_1+q_2}}^{2(t_0+t)}[\mu_n]\right)^{\frac{1}{p-1}}\leq \varepsilon \lambda\right\}\right|\\
	&\hspace{2cm}\leq C\varepsilon^{\frac{N(q_1+q_2)(p-1)}{2q_2(N-\beta p)}} 	\left|\left\{\mathbf{L}^{2(t_0+t)}_{\alpha p, \frac{q_1}{(q_1+q_2)(p-1)}}[\mu_n]\mathbf{L}^{2(t_0+t)}_{\beta p, \frac{q_2}{(q_1+q_2)(p-1)}}[\mu_n]>\varepsilon^{1/2}\lambda\right\}\right|
	\end{align}
	for any $\varepsilon\in (0,\varepsilon_0)$ and $t>0$, for some $a>1$. This gives 
	\begin{align}\label{es7'}
	\nonumber	&\left|\left\{\left(\mathbf{L}^{2(t_0+t)}_{\alpha p, \frac{q_1}{(q_1+q_2)(p-1)}}[\mu_n]\mathbf{L}^{2(t_0+t)}_{\beta p, \frac{q_2}{(q_1+q_2)(p-1)}}[\mu_n]\right)^{q_1+q_2}>a\lambda, \left(\mathbf{I}_{\frac{\alpha pq_1+p\beta q_2}{q_1+q_2}}^{2(t_0+t)}[\mu_n]\right)^{\frac{q_1+q_2}{p-1}}\leq \varepsilon \lambda\right\}\right|\\
	&\hspace{2cm}\leq C\varepsilon^{\frac{N(p-1)}{2q_2(N-\beta p)}} 	\left|\left\{\left(\mathbf{L}^{2(t_0+t)}_{\alpha p, \frac{q_1}{(q_1+q_2)(p-1)}}[\mu_n]\mathbf{L}^{2(t_0+t)}_{\beta p, \frac{q_2}{(q_1+q_2)(p-1)}}[\mu_n]\right)^{q_1+q_2}>\varepsilon^{1/2}\lambda\right\}\right|
	\end{align}
	for any $\varepsilon\in (0,\varepsilon_0)$ and $t>0$, for some $a>1$.
	 It is easy to obtain from the above two inequality  that 
	\begin{align*}
	&\int_{0}^{\infty}\phi(\lambda)\left|\left\{\left(\mathbf{L}^{2(t_0+t)}_{\alpha p, \frac{q_1}{(q_1+q_2)(p-1)}}[\mu_n]\mathbf{L}^{2(t_0+t)}_{\beta p, \frac{q_2}{(q_1+q_2)(p-1)}}[\mu_n]\right)^{q_1+q_2}>\lambda\right\}\right|d\lambda\\
	&\leq C\varepsilon^{\frac{N(p-1)}{2q_2(N-\beta p)}}  	\int_0^\infty \phi(\lambda)\left|\left\{\left(\mathbf{L}^{2(t_0+t)}_{\alpha p, \frac{q_1}{(q_1+q_2)(p-1)}}[\mu_n]\mathbf{L}^{2(t_0+t)}_{\beta p, \frac{q_2}{(q_1+q_2)(p-1)}}[\mu_n]\right)^{q_1+q_2}>a^{-1}\varepsilon^{1/2}\lambda\right\}\right|\, d\lambda\\&+ 
	C	\int_{0}^{\infty}\phi(\lambda)\left|\left\{\left(\mathbf{I}_{\frac{\alpha pq_1+p\beta q_2}{q_1+q_2}}^{2(t_0+t)}[\mu_n]\right)^{\frac{q_1+q_2}{p-1}}>\varepsilon \lambda\right\}\right|d\lambda\\&
	\leq C\varepsilon^{\frac{N(p-1)}{2q_2(N-\beta p)}-1/2}  \int_{0}^{\infty}\phi(a\varepsilon^{-1/2}\lambda)\left|\left\{\left(\mathbf{L}^{2(t_0+t)}_{\alpha p, \frac{q_1}{(q_1+q_2)(p-1)}}[\mu_n]\mathbf{L}^{2(t_0+t)}_{\beta p, \frac{q_2}{(q_1+q_2)(p-1)}}[\mu_n]\right)^{q_1+q_2}>\lambda\right\}\right|d\lambda
	\\&+ C	\int_{0}^{\infty}\phi(\varepsilon\lambda)\left|\left\{\left(\mathbf{I}_{\frac{\alpha pq_1+p\beta q_2}{q_1+q_2}}^{2(t_0+t)}[\mu_n]\right)^{\frac{q_1+q_2}{p-1}}> \lambda\right\}\right|d\lambda
	\end{align*}
	Since $\phi(\varepsilon\lambda),\phi(a\varepsilon^{-1/2}\lambda)\leq C|\log(\varepsilon)|\phi(\lambda)$ for any $\lambda>0$,  $\varepsilon<<1$ and $\frac{N(p-1)}{2q_2(N-\beta p)}-1/2>0$, so it is easy to get that 
	\begin{align*}
	\int_{0}^{\infty}\phi(\lambda)\left|\left\{\left(\mathbf{L}^{2(t_0+t)}_{\alpha p, \frac{q_1}{(q_1+q_2)(p-1)}}[\mu_n]\mathbf{L}^{2(t_0+t)}_{\beta p, \frac{q_2}{(q_1+q_2)(p-1)}}[\mu_n]\right)^{q_1+q_2}>\lambda\right\}\right|d\lambda\leq  C
	\end{align*}
	Hence, $$ \left(\mathbf{L}^{2(t_0+t)}_{\alpha p, \frac{q_1}{(q_1+q_2)(p-1)}}[\mu_n]\mathbf{L}^{2(t_0+t)}_{\beta p, \frac{q_2}{(q_1+q_2)(p-1)}}[\mu_n]\right)^{q_1+q_2}	$$ is equi-integrable in  $B_{t}(0)$. The proof is complete.
\qeda\medskip

    The next statement is the analogue of Theorem \ref{th-pot-4} in a bounded domain.    

\begin{theorem}\label{th-pot-5}
   Let $\alpha,\beta,q_1,q_2>0$, $\alpha>\beta$, $1<p<\frac N\ga$, $q_1+q_2>p-1$, $q_2<\frac{N(p-1)}{N-\beta p}$,  $\omega\in \mathfrak{M}_b^+(B_R(x_0))$ for some $B_R(x_0)\subset\mathbb{R}^N$, extended by $0$ in $B^c_R(x_0)$. Then, the following statements are  equivalent:
   \begin{description}
   	\item[(a)] The inequality 
   	\begin{align}\label{pot-46}
   	\omega(K)\leq C_1\text{Cap}_{\mathbf{G}_{\frac{\alpha pq_1+\beta pq_2}{q_1+q_2}},\frac{q_1+q_2}{q_1+q_2-p+1}}(K),
   	\end{align}
   	holds for any compact set $K\subset\mathbb{R}^N$, for some $C_1=C_1(R)>0$.
   	\item[(b)] The inequality 
   	\begin{align}\label{pot-47}
   	\int_{K}\left(\mathbf{W}^{4R}_{\alpha,p}[\omega](x)\right)^{q_1}\left(\mathbf{W}^{4R}_{\beta,p}[\omega](x)\right)^{q_2}dx\leq C_2\text{Cap}_{\mathbf{G}_{\frac{\alpha pq_1+\beta pq_2}{q_1+q_2}},\frac{q_1+q_2}{q_1+q_2-p+1}}(K),
   	\end{align}
   	holds for any compact set $K\subset\mathbb{R}^N$, for some $C_2=C_2(R)>0$.
   	   	\item[(c)] The inequality 
   	   	\begin{align}\label{pot-48}
   	   	\int_{\mathbb{R}^N}\left(\mathbf{W}^{4R}_{\frac{\alpha q_1+\beta q_2}{q_1+q_2},p}[\chi_{_{B_t(x)}}\omega](y)\right)^{q_1+q_2}dy\leq C_3 \omega(B_t(x)),
   	   	\end{align}
   	   	holds for any ball $B_t(x)\subset\mathbb{R}^N$, for some $C_3=C_3(R)>0$.
   	   	\item[(d)] The inequality 
   	   	\begin{align}\label{pot-49}
   	   	\int_{\mathbb{R}^N}\left(\mathbf{W}^{4R}_{\alpha,p}[\chi_{_{B_t(x)}}\omega](y)\right)^{q_1}\left(\mathbf{W}^{4R}_{\beta,p}[\chi_{_{B_t(x)}}\omega](y)\right)^{q_2}dy\leq C_4 \omega(B_t(x)),
   	   	\end{align}
   	   	holds for any ball $B_t(x)\subset\mathbb{R}^N$, for some $C_4=C_4(R)>0$.
   	   	   	\item[(e)]  The system of inequalities 
   	   	   	\bel{pot-50}\BA {lll}
   	   	   	(i)\qquad\qquad\mathbf{W}^{4R}_{\alpha,p}\left[\left(\mathbf{W}^{4R}_{\alpha,p}[\omega]\right)^{q_1}\left(\mathbf{W}^{4R}_{\beta,p}[\omega]\right)^{q_2}\right]\leq C_5 \mathbf{W}^{4R}_{\alpha,p}[\omega]\\[2mm]
   	   	   	(ii)\qquad\qquad\mathbf{W}^{4R}_{\beta,p}\left[\left(\mathbf{W}^{4R}_{\alpha,p}[\omega]\right)^{q_1}\left(\mathbf{W}^{4R}_{\beta,p}[\omega]\right)^{q_2}\right]\leq C_5 \mathbf{W}^{4R}_{\beta,p}[\omega],
   	   	   	\EA\ee
   	   	   	holds in $B_{2R}(x_0)$  for some $C_5=C_5(R)>0$.
   \end{description}
       \end{theorem}

\nind\Proof 
By Theorem \ref{th-pot-2} we have (a) $ \Leftrightarrow $ (c); by Theorem \ref{th-pot-3}, (c) $ \Leftrightarrow $ (d). As in the proof  of Theorem \ref{th-pot-4}, we can see that (e) $ \Rightarrow  $ (a) and (e) $ \Rightarrow  $ (b). Since 
   \begin{align*}
   \left(\mathbf{W}^{4R}_{\alpha,p}[\omega](x)\right)^{q_1\gamma}\left(\mathbf{W}^{4R}_{\beta,p}[\omega](x)\right)^{q_2\gamma}&\geq C\int_{0}^{4R}\left(\frac{\omega(B_r(x))}{r^{N-\alpha p}}\right)^{\frac{q_1\gamma}{p-1}}\left(\frac{\omega(B_r(x))}{r^{N-\beta p}}\right)^{\frac{q_2\gamma}{p-1}}\frac{dr}{r}\\[1mm]&=C\mathbf{W}^{4R}_{\alpha_0,p_0}[\omega](x)\quad\text{ for all }\; x\in B_{2R}(x_0),
   \end{align*}
   where $\gamma=\frac{1}{q_1}+\frac{1}{q_2},$ $\alpha_0=\frac{\gamma(\alpha p q_1+\beta p q_2)}{\gamma(q_1+q_2)+p-1}$ and $p_0=\frac{\gamma(q_1+q_2)+p-1}{\gamma(q_1+q_2)}$, then, (b) implies that
    \begin{align}\label{pot-50*}
        \int_{K}\left(\mathbf{W}^{4R}_{\alpha_0,p_0}[\omega](x)\right)^{\frac{1}{\gamma}}dx\leq C\text{Cap}_{\mathbf{G}_{\frac{\alpha pq_1+\beta pq_2}{q_1+q_2}},\frac{q_1+q_2}{q_1+q_2-p+1}}(K)=\text{Cap}_{\mathbf{G}_{\alpha_0p_0},\frac{\frac{1}{\gamma}}{\frac{1}{\gamma}-p_0+1}}(K),
        \end{align}
        is verified for any compact set $K\subset\mathbb{R}^N$. Therefore (a) follows  by Theorem \ref{th-pot-2}.\smallskip
        
    It remains to prove (a)+(c)+(d) $ \Rightarrow  $ (e).
From (a) we have
\begin{align}\label{pot-51}
 \omega(B_r(x))\leq Cr^{N-\frac{\alpha pq_1+\beta pq_2}{q_1+q_2-p+1}}~~\text{for all $x\in\mathbb{R}^N$ and $r>0$}.
 \end{align}
 From (b)
  \begin{align*}
   \int_{B_r(x)}\left(\mathbf{W}^{r}_{\frac{\alpha q_1+\beta q_2}{q_1+q_2},p}[\omega](y)\right)^{q_1+q_2}dy\leq C_2 \omega(B_{2r}(x)) ~~\text{for all $x\in\mathbb{R}^N$ and $0<r<8R$}.
   \end{align*} 
 Using H\"older's inequality and $\mathbf{W}^{r}_{\frac{\alpha q_1+\beta q_2}{q_1+q_2},p}[\omega]\geq r^{-\frac{(\alpha-\beta)pq_2}{(p-1)(q_1+q_2)}}\mathbf{W}^{r}_{\alpha ,p}[\omega]$, we get 
 \begin{align}\label{pot-52}
    \int_{B_r(x)}\left(\mathbf{W}^{r}_{\alpha,p}[\omega](y)\right)^{q_1}dy\leq C r^{\frac{(\alpha-\beta)pq_1q_2+(p-1)Nq_2+(N-\beta p)(p-1)q_1}{(p-1)(q_1+q_2)}}\left(\frac{\omega(B_{2r}(x))}{r^{N-\beta p}}\right)^{\frac{q_1}{q_1+q_2}}, 
    \end{align}
    for all $x\in\mathbb{R}^N$ and  $0<r<8R.$\\
    From (c), 
    \begin{align}\label{pot-53}
      \int_{B_r(x)}\left(\mathbf{W}^r_{\alpha,p}[\omega](y)\right)^{q_1}\left(\mathbf{W}^r_{\beta,p}[\omega](y)\right)^{q_2}dy\leq C_3 \omega(B_{2r}(x)) ~~\text{for all $x\in\mathbb{R}^N$ and } 0<r<8R.
      \end{align}
By Lemma \ref{lem-pot-4}, 
\begin{align}\label{pot-54}
 \int_{B_r(x)}\left(\mathbf{W}^r_{\beta,p}[\mu]\right)^{q_2}dy\leq Cr^{N}\left(\frac{\mu(B_{2r}(x))}{r^{N-\beta p}}\right)^{\frac{q_2}{p-1}} ~~\text{for all $x\in\mathbb{R}^N$ and } 0<r<8R.
 \end{align}
Next we have 
for $\eta=\alpha$ or $\eta=\beta$ and almost all $x\in B_{2R}(x_0)$,
\begin{align}\label{pot-55}
\mathbf{W}^{4R}_{\eta,p}\left[\left(\mathbf{W}^{4R}_{\alpha,p}[\mu]\right)^{q_1}\left(\mathbf{W}^{4R}_{\beta,p}[\mu]\right)^{q_2}\right](x)\leq C\sum_{i=1}^{4}\int_{0}^{4R}\left(\frac{\mathbf{A}_i(x,r)}{r^{N-\eta p}}\right)^{\frac{1}{p-1}}\frac{dr}{r},
\end{align}
where 
\begin{align*}
&\mathbf{A}_1(x,r)=\int_{B_r(x)}\left(\mathbf{W}^r_{\alpha,p}[\mu](y)\right)^{q_1}\left(\mathbf{W}^r_{\beta,p}[\mu](y)\right)^{q_2}dy,\\&
\mathbf{A}_2(x,r)=\int_{B_r(x)}\left(\mathbf{W}^r_{\alpha,p}[\mu](y)\right)^{q_1}\left(\int_{r}^{4R}\left(\frac{\mu(B_t(y))}{t^{N-\beta p}}\right)^{\frac{1}{p-1}}\frac{dt}{t}\right)^{q_2}dy,
\\&
\mathbf{A}_3(x,r)=\int_{B_r(x)}\left(\int_{r}^{4R}\left(\frac{\mu(B_t(y))}{t^{N-\alpha p}}\right)^{\frac{1}{p-1}}\frac{dt}{t}\right)^{q_1}\left(\mathbf{W}^r_{\beta,p}[\mu](y)\right)^{q_2}dy,
\\&
\mathbf{A}_4(x,r)=\int_{B_r(x)}\left(\int_{r}^{4R}\left(\frac{\mu(B_t(y))}{t^{N-\alpha p}}\right)^{\frac{1}{p-1}}\frac{dt}{t}\right)^{q_1}\left(\int_{r}^{4R}\left(\frac{\mu(B_t(y))}{t^{N-\beta p}}\right)^{\frac{1}{p-1}}\frac{dt}{t}\right)^{q_2}dy.
\end{align*}
Thanks to   \eqref{pot-40} there holds
\begin{align*}
\mathbf{A}_1(x,r)\leq C\mu(B_{2r}(x)),
\end{align*}
which implies 
\begin{align}\label{pot-56}
\int_{0}^{4R}\left(\frac{\mathbf{A}_1(x,r)}{r^{N-\eta p}}\right)^{\frac{1}{p-1}}\frac{dr}{r}\leq C \mathbf{W}^{8R}_{\eta,p}[\mu](x)\leq C \mathbf{W}^{4R}_{\eta,p}[\mu](x)~~\text{for all } x\in B_{2R}(x_0).
\end{align}
Since $B_t(y)\leq B_{2t}(x)$ for any $y\in B_r(x)$ and $t\geq r$ , and thanks to \eqref{pot-39} and \eqref{lem-pot-4} we deduce that there holds, for $0<r<4R$ and $ x\in B_{2R}(x_0)$,
\begin{align*}
&\mathbf{A}_2(x,r)\leq \int_{B_r(x)}\left(\mathbf{W}^r_{\alpha,p}[\mu](y)\right)^{q_1}dy\left(\int_{r}^{4R}\left(\frac{\mu(B_{2t}(x))}{t^{N-\beta p}}\right)^{\frac{1}{p-1}}\frac{dt}{t}\right)^{q_2}\\
&\phantom{\mathbf{A}_2(x,r)}\leq C r^{\frac{(\alpha-\beta)pq_1q_2+(p-1)Nq_2+(N-\beta p)(p-1)q_1}{(p-1)(q_1+q_2)}}\!\left(\frac{\mu(B_{2r}(x))}{r^{N-\beta p}}\right)^{\frac{q_1}{q_1+q_2}}\! \left(\int_{r}^{4R}\!\left(\frac{\mu(B_{2t}(x))}{t^{N-\beta p}}\right)^{\frac{1}{p-1}}\!\frac{dt}{t}\right)^{q_2}
\\&\phantom{\mathbf{A}_2(x,r)} \leq C r^{\frac{(\alpha-\beta)pq_1q_2+(p-1)Nq_2+(N-\beta p)(p-1)q_1}{(p-1)(q_1+q_2)}} \left(\int_{r}^{4R}\left(\frac{\mu(B_{2t}(x))}{t^{N-\beta p}}\right)^{\frac{1}{p-1}}\frac{dt}{t}\right)^{q_2+\frac{q_1(p-1)}{q_1+q_2}}.
\end{align*}
Next
\begin{align*}
&\mathbf{A}_3(x,r)\leq \left(\int_{r}^{4R}\left(\frac{\mu(B_{2t}(x))}{t^{N-\alpha p}}\right)^{\frac{1}{p-1}}\frac{dt}{t}\right)^{q_1}\int_{B_r(x)}\left(\mathbf{W}^r_{\beta,p}[\mu](y)\right)^{q_2}dy
\\&\phantom{\mathbf{A}_3(x,r)}\leq C \left(\int_{r}^{4R}\left(\frac{\mu(B_{2t}(x))}{t^{N-\alpha p}}\right)^{\frac{1}{p-1}}\frac{dt}{t}\right)^{q_1}r^{N}\left(\frac{\mu(B_{2r}(x))}{r^{N-\beta p}}\right)^{\frac{q_2}{p-1}}
\\&\phantom{\mathbf{A}_3(x,r)}\leq C r^{N} \left(\int_{r}^{4R}\left(\frac{\mu(B_{2t}(x))}{t^{N-\alpha p}}\right)^{\frac{1}{p-1}}\frac{dt}{t}\right)^{q_1}\left(\int_{r}^{4R}\left(\frac{\mu(B_{2t}(x))}{t^{N-\beta p}}\right)^{\frac{1}{p-1}}\frac{dt}{t}\right)^{q_2},
\end{align*}
and
\begin{align*}
&\mathbf{A}_4(x,r)\leq C r^{N} \left(\int_{r}^{4R}\left(\frac{\mu(B_{2t}(x))}{t^{N-\alpha p}}\right)^{\frac{1}{p-1}}\frac{dt}{t}\right)^{q_1}\left(\int_{r}^{4R}\left(\frac{\mu(B_{2t}(x))}{t^{N-\beta p}}\right)^{\frac{1}{p-1}}\frac{dt}{t}\right)^{q_2}.
\end{align*}
As in the proof of  Theorem \ref{th-pot-4}, we easily  obtain 
\begin{align*}
\int_{0}^{4R}\left(\frac{\mathbf{A}_2(x,r)}{r^{N-\eta p}}\right)^{\frac{1}{p-1}}\frac{dr}{r}\leq C \mathbf{W}^{4R}_{\eta,p}[\mu](x)~~\text{for all } x\in B_{2R}(x_0),
\end{align*}
and
\begin{align*}
\int_{0}^{4R}\left(\frac{\mathbf{A}_3(x,r)}{r^{N-\eta p}}\right)^{\frac{1}{p-1}}\frac{dr}{r}+\int_{0}^{4R}\left(\frac{\mathbf{A}_4(x,r)}{r^{N-\eta p}}\right)^{\frac{1}{p-1}}\frac{dr}{r}\leq C \mathbf{W}^{4R}_{\eta,p}[\mu](x)~~\text{for all } x\in B_{2R}(x_0).
\end{align*}
Combining these inequalities with \eqref{pot-55} and \eqref{pot-56}, we get (e). 
 \qeda    \medskip

   \section{Renormalized solutions}
   Let $\Omega$ be a bounded domain in $\mathbb{R}^N$.   If $\mu\in\mathfrak{M}_b(\Omega)$, we denote by $\mu^+$ and $\mu^-$ respectively its positive and negative parts in the Jordan decomposition. We denote by $\mathfrak{M}_0(\Omega)$ the space of diffuse measures in $\Omega$ and by $\mathfrak{M}_s(\Omega)$ the space of measures in $\Omega$ which are singular with respect to the $\text{Cap}_{\mathbf{G}_1,p}$ which means that their support is set of zero $\text{Cap}_{\mathbf{G}_1,p}$-capacity. Classically, any $\mu\in\mathfrak{M}_b(\Omega)$ can be written in a unique way under the form $\mu=\mu_0+\mu_s$ where $\mu_0\in \mathfrak{M}_0(\Omega)\cap \mathfrak{M}_b(\Omega)$ and $\mu_s\in \mathfrak{M}_s(\Omega)$.
   It is well known  that any  $\mu_0\in \mathfrak{M}_0(\Omega)\cap\mathfrak{M}_b(\Omega)$ can be written under the form $\mu_0=f-\text{div} ~g$ where $f\in L^1(\Omega)$ and $g\in L^{p'}(\Omega,\mathbb{R}^N)$.
    
    For $k>0$ and $s\in\mathbb{R}$ we set $T_k(s)=\max\{\min\{s,k\},-k\}$. If $u$ is a measurable function defined  in $\Omega$, finite a.e. and such that $T_k(u)\in W^{1,p}_{loc}(\Omega)$ for any $k>0$, there exists a measurable function $v:\Omega\mapsto \mathbb{R}^N$ such that $\nabla T_k(u)=\chi_{_{|u|\leq k}}v$ 
    a.e. in $\Omega$ and for all $k>0$. We define the gradient of $u$ by $v=\nabla u$ almost everywhere. We recall the definition of a renormalized solution given in \cite{22DMOP}.
    
    \begin{definition} Let $A:\mathbb{R}^N\mapsto \mathbb{R}^N$ satisfy \eqref{In-3}. Let $\mu=\mu_0+\mu_s\in\mathfrak{M}_b(\Omega)$. A measurable  function $u$ defined in $\Omega$ and finite a.e. is called a renormalized solution of 
    
   \begin{equation}
   \label{renorm-1}\BA{lll}
     - \text{div}(A(x,\nabla u)) = \mu \qquad&\text{in }\;\Omega,  \\ 
    \phantom{- \text{div}(A(x,\nabla ))}
    u = 0\qquad&\text{on }\;\partial \Omega,  \\ 
\EA
   \end{equation}
     if $T_k(u)\in W^{1,p}_0(\Omega)$ for any $k>0$, $|{\nabla u}|^{p-1}\in L^r(\Omega)$ for any $0<r<\frac{N}{N-1}$, and $u$ has the property that for any $k>0$ there exist $\lambda_k^+$ and $\lambda_k^-$ belonging to $\mathfrak{M}_{b}^+\cap\mathfrak{M}_0(\Omega)$, respectively concentrated on the sets $u=k$ and $u=-k$, with the property that 
     $\mu_k^+\rightharpoonup\mu_s^+$, $\mu_k^-\rightharpoonup\lambda_s^-$ in the narrow topology of measures and such that
    \[
    \int_{\{|u|<k\}}A(x,\nabla u).\nabla\varphi
    dx=\int_{\{|u|<k\}}{\varphi d}{\mu_{0}}+\int_{\Omega}\varphi d\lambda_{k}%
    ^{+}-\int_{\Omega}\varphi d\lambda_{k}^{-},%
    \]
     for every $\varphi\in W^{1,p}_0(\Omega)\cap L^{\infty}(\Omega)$.
   
    \end{definition}
\begin{proposition}\label{renorm-prop-1}{\rm\cite{TrWa}} If $\mu \in \mathfrak{M}_0(\Omega)$, then problem \eqref{renorm-1} has a unique renormalized solution. 
    \end{proposition}                                                                                       
    We recall the next two important results which are proved in \cite[Th 4.1, Sec 5.1]{22DMOP}.
\begin{theorem}
    \label{recall} Let $\{\mu_n\}\subset \mathfrak{M}_b(\Omega)$ be a sequence such that $\sup_n|{\mu_n}|(\Omega)<\infty$ and let $\{u_n\}$ be renormalized solutions of 
    \begin{equation}
    \label{renorm-2}\begin{array}{ll}
     - \operatorname{\text{div}}A(x,\nabla u_n)) = \mu_n \qquad&\text{in }\;\Omega  \\ 
    \phantom{- \text{div}(A(x,\nabla ))}
    u_n = 0\qquad&\text{on }\;\partial \Omega.  \\ 
\EA
    \end{equation}
    Then, up to a subsequence, $\{u_n\}$ converges a.e. to a solution $u$ of $-\text{div}(A(x,\nabla u))=\mu$ in the sense of distributions in $\Omega$, for some measure $\mu \in \mathfrak{M}_b(\Omega)$, and for every $k>0$, $k^{-1}\norm{\nabla T_k(u)}_{L^p}^{p}\leq M$ for some $M>0$.
        \end{theorem}                                                                                       

    The following fundamental stability result of \cite{22DMOP} extends Theorem \ref{recall}.
\begin{theorem}
    \label{stab} Let $\mu=\mu_0+\mu_s^+-\mu_s^-\in \mathfrak{M}_b(\Omega)$, with $\mu_0=f-\text{\text{div}}\, g\in \mathfrak M_0(\Omega)$, $\mu_s^+,\mu_s^-\in \mathfrak{M}_s^+(\Omega)$. Assume there are sequences
    $\{f_n\}\subset L^1(\Omega)$, $\{g_n\}\subset (L^{p'}(\Omega))^N$, $\{\eta_n^1\}, \{\eta_n^2\}\subset \mathfrak{M}_b^+(\Omega)$ such that $f_n\rightharpoonup f$ weakly in $L^1(\Omega)$, $g_n\to g$ in $L^{p'}(\Omega)$ and $\text{\text{div}}\, g_n$ is bounded in $\mathfrak{M}_b(\Omega)$, $\eta_n^1\rightharpoonup \mu_s^+$ and $\eta_n^2\rightharpoonup \mu_s^-$ in the narrow topology. If $\mu_n=f_n-\text{\text{div}}\,  g_n+\eta_n^1-\eta_n^2$
    and $u_n$ is a renormalized solution of (\ref{renorm-2}), then, up to a subsequence, $u_n$ converges a.e. to a renormalized solution $u$ of \eqref{renorm-1}. Furthermore, $T_k(u_n)\to T_k(u)$ in $W^{1,p}_0(\Omega)$ for any $k>0$. 
    \end{theorem} 
    
      The next result is proved  in \cite[Th 3.2]{22Bi3}. Therein it plays an important role in study the  stability of the renormalized solutions of 
      the following problem with absorption,
    \begin{equation} \label{14}\begin{array}{ll}
    -\text{div} (A(x,\nabla u_{n_k}))+|u_{n_k}|^{q-1}u_{n_k}=\gm_{n_k}\qquad&\text{in }\Omega,\\[1mm]
    \phantom{ -------------,}
    u_{n_k}=0&\text{on }\partial\Omega.\\
\end{array}
    \end{equation}

\begin{theorem} \label{renorm-th-1}Let $\{n_k\}_k$ be an increasing sequence in $\mathbb{N}$, $q>p-1$,  $\{\mu_{n_k}\}_{k}$ be a sequence in $\mathfrak{M}(\mathbb{R}^N)$ such that \begin{equation*}
    \mathop {sup}\limits_{k\geq k_0} \left| {{\mu _{n_k}}} \right|\left( B_{n_{k_0}}(0) \right) <  + \infty ~~\text{for all}~~k_0\in\mathbb{N}.
    \end{equation*}
    Let $u_{n_k}$ be a renormalized solution of  \eqref{renorm-1} with data $\mu_{n_k}$ and $\Omega=B_{n_k}(0)$ such that $\{|u_{n_k}|^q\}_{k\geq k_0}$ is bounded in $L^1(B_{n_{k_0}}(0))$ for any $k_0$. Then, there exist subsequence of $\{u_{n_k}\}_{k}$, still denoted by $\{u_{n_k}\}_{k}$ a measure $\gm$ and measurable function $u$ such that $\gm_{n_k}\rightharpoonup \gm$ in the weak sense of measures, $u_{n_k}\to u$, $\nabla u_{n_k}\to \nabla u$ a.e in $\mathbb{R}^N$. Moreover, $|\nabla u_{n_k}|^{p-2}\nabla u_{n_k} \to |\nabla u|^{p-2}\nabla u$ in $L^s_{loc}(\mathbb{R}^N)$ for all $ 0\le s<\frac{N}{N-1}$ and $u$ satisfies \eqref{renorm-1} in the sense of distributions in $\mathbb{R}^N$.
   \end{theorem}

\begin{theorem}{\rm\cite{22PhVe,55VHV}}
   \label{renorm-th-2}Let $\Omega$ be a bounded domain of $\mathbb{R}^N$. Then there exists a constant $C=C(N,p,\Lambda_1,\Lambda_2)> 1$ such that if $\mu\in \mathfrak{M}_b(\Omega)$ and $u$ is a renormalized solution of problem \eqref{renorm-1} there holds
   \begin{equation}
   \label{renorm-3} |u(x)|\leq C{\bf W}^{2R}_{1,p}[|\mu|](x)
    \quad\text{a.e. in }\Omega,
   \end{equation}
   where $R=diam(\Omega).$ Moreover, if $\mu\geq 0$ and $u\geq 0$ then, 
   \begin{align}\label{041120146}
   u(x)\geq \frac{1}{C}\mathbf{W}_{1,p}^{\frac{d(x,\partial\Omega)}{4}}[\mu](x)~~ \text{a.e in }~\Omega.
   \end{align}
   \end{theorem}
\begin{theorem}{\rm\cite{DuMin,KuMin1,QH2}}
      \label{renorm-th-3} Let $\Omega$ be a bounded domain of $\mathbb{R}^N$.  Then there exists a constant $C=C(N,p,\Lambda_1,\Lambda_2,diam(\Omega))>0$ such that if $\mu \in C_b(\Omega)$ and $u$ is a solution of problem \eqref{renorm-1} there holds
      \begin{equation}
      \label{renorm-4}|\nabla u(x)|\leq C\left({\bf I}^{r}_{1}[|\mu|](x)\right)^{\frac{1}{p-1}}+C\left(\fint_{B_r(x)}|\nabla u|^{\gamma_0}dx\right)^{1/\gamma_0},
      \end{equation}
      for any $B_r(x)\subseteq \Omega$ and for some $\gamma_0\in (0,\frac{N(p-1)}{N-1})$.  Moreover, if $A(x,\xi)=A(\xi)$ for any $(x,\xi)\in \mathbb{R}^N\times\mathbb{R}^N$, then the constant in \eqref{renorm-4} does not depend on $diam(\Omega)$.
      \end{theorem}
\begin{corollary}
      \label{renorm-cor-1} Let $\Omega$ be a bounded domain of $\mathbb{R}^N$,  $R=diam(\Omega)$, $\mu\in \mathfrak{M}_b(\Omega)$. Then there exists a constant $C=C(N,p,\Lambda_1,\Lambda_2)>0$ and a renormalized solution $u$ of problem \eqref{renorm-1} such that
            \begin{equation}
            \label{renorm-5}|\nabla u(x)|\leq C\left(\frac{R}{\delta}\right)^{N}\left({\bf I}^{2R}_{1}[|\mu|](x)\right)^{\frac{1}{p-1}},
            \end{equation}
            for any $x\in \Omega$ such that $d(x,\partial\Omega)>\delta$ with $\delta\in (0,R/2)$.   Moreover, if $A(x,\xi)=A(\xi)$ for any $(x,\xi)\in \mathbb{R}^N\times\mathbb{R}^N$, then the constant in \eqref{renorm-5} does not depend on $R$.
      \end{corollary}
      \Proof We can choose $\mu_n\in C_c^\infty(\Omega)$ such that $\mu_n$ converges to $\mu$ in the sense of theorem \ref{stab} and $|\mu_n|\leq \varphi_n*|\mu|$, where $\{\varphi_n\}$ is a sequence of mollifiers in $\mathbb{R}^N$. Let $u_n$ be solutions of problem \eqref{renorm-1} with data $\mu_n$. Fixed $\delta\in (0,R/2)$, by Theorem \ref{renorm-th-3}, we have 
      \begin{equation*}
          |\nabla u_n(x)|\leq C\left({\bf I}^{\delta/2}_{1}[|\mu_n|](x)\right)^{\frac{1}{p-1}}+C\fint_{B_{\delta/2}(x)}|\nabla u_n|dx,
            \end{equation*}
            for any $x \in \Omega, d(x,\partial\Omega)>\delta$.
          Notice that (see e.g. \cite{22DMOP})
          \begin{align*}
          |\{|\nabla u_n |>s\}|\leq C \frac{(|\mu_n |(\Omega))^{\frac{N}{N-1}}}{s^{\frac{N(p-1)}{N-1}}}~~\text{for all}~s>0.
          \end{align*} 
It leads to 
          \begin{align*}
          \int_{\Omega}|\nabla u_n(x)|dx\leq C R^N\left(\frac{|\mu_n|(\Omega)}{R^{N-1}}\right)^{\frac{1}{p-1}}.
          \end{align*}
          Thus, 
          \begin{align*}
          |\nabla u_n(x)|&\leq C\left(\frac{R}{\delta}\right)^{N}\left({\bf I}^{2R}_{1}[|\mu_n|](x)\right)^{\frac{1}{p-1}}\\&
          \leq C\left(\frac{R}{\delta}\right)^{N}\left((\varphi_n*{\bf I}^{2R}_{1}[|\mu|])(x)\right)^{\frac{1}{p-1}},
          \end{align*}          
                      for any $x \in \Omega, d(x,\partial\Omega)>\delta$. \\
     On the other hand, by theorem \ref{stab}, there exists a subsequence of $\{u_n\}$ converging to a renormalized solution $u$ of \eqref{renorm-1} with data $\mu_n$. 
     Therefore, $u$ satisfies \eqref{renorm-5} since $\varphi_n*{\bf I}^{2R}_{1}[|\mu|]\to{\bf I}^{2R}_{1}[|\mu|]$ almost everywhere. 
      
      \qeda
    \section{Proof of the main results}

{\it Proof of Theorem \ref{mainthem1}}.\smallskip

   \nind{\it Step 1: Case $\frac{3N-2}{2N-1}<p\leq 2$}. Let $\mu_{n,k}\in C_c^\infty(B_{2k}(0))$ for $k\in \mathbb{N}$ such that $\mu_{n,k}$ converges to $\chi_{_{B_{k}(0)}}\mu$ in the sense of theorem \ref{stab} with $\Omega=B_{2k}(0)$ and $|\mu_{n,k}|\leq \varphi_n*(\chi_{_{B_{k}(0)}}|\mu|)$, where $\{\varphi_n\}$ is a sequence of mollifiers in $\mathbb{R}^N$. Thanks to Proposition \ref{prop-pot-1}, 
    \begin{align}\label{cap-appthm1}
   |\mu_{n,k}|(K)\leq C'C\text{Cap}_{\mathbf{I}_{\frac{q_1p+q_2}{q_1+q_2}},\frac{q_1+q_2}{q_1+q_2-p+1}}(K)\quad\text{for all compact }K\subset\mathbb{R}^N
   \end{align} 
   We will prove that if $C$ in \eqref{cap-appthm1} is small enough, then for any $k\geq1,n\in \mathbb{N}$ the problem 
   \begin{equation}\label{pr-3}
   \BA{lll}%
   -\text{div}(A(\nabla u_{n,k}))=\chi_{_{B_k(0)}}|u_{n,k}|^{q_1-1}u|\nabla u_{n,k}|^{q_2}+\mu_{n,k}\quad&\text{in }\;B_{2k}(0),\\[1mm] 
   \phantom{-\text{div}(A(\nabla ))}
   u_{n,k}=0&\text{on }\;\partial B_{2k}(0),                                                                                          
   \EA
   \end{equation}
   has a renormalized solution  satisfying 
   \begin{align}\label{esstep1thm1}|u_{n,k}|\leq C_0\left(\mathbf{I}_{p}[|\mu_{n,k}|]\right)^{\frac{1}{p-1}},~~~|\nabla u_{n,k}|\leq C_0\left(\mathbf{I}_{1}[|\mu_{n,k}|]\right)^{\frac{1}{p-1}}~~\text{in}~~B_k(0).
   \end{align}
   By Theorem  \eqref{th-pot-4}, we need to prove  that there exists $M>0$ such that, if for $\alpha=1$ and $p$, the following inequalities hold, 
       \begin{align}\label{pr-1}\mathbf{I}_{\alpha}\left[\left(\mathbf{I}_{p}[|\mu_{n,k}|]\right)^{\frac{q_1}{p-1}}\left(\mathbf{I}_{1}[|\mu_{n,k}|]\right)^{\frac{q_2}{p-1}}\right]\leq M \mathbf{I}_{\alpha}[|\mu_{n,k}|]<\infty~~\text{ almost everywhere},
       \end{align}
then problem \eqref{pr-3}  has a renormalized solution  satisfying \eqref{esstep1thm1}.

              For any $k\in\mathbb{N}$,
     we set \begin{equation*}
                   \mathbf{E}_\Lambda=\left\{u:|u|\leq \Lambda \left(\mathbf{I}_p[|\mu_{n,k}|]\right)^{\frac{1}{p-1}},~~|\nabla u|\leq \Lambda \left(\mathbf{I}_1[|\mu_{n,k}|]\right)^{\frac{1}{p-1}}~\text{ in }~B_{k}(0)\right\}.
                   \end{equation*}
                   Since $\mu_{n,k}\in C_c^\infty(B_{2k}(0))$, $ \mathbf{E}_\Lambda\subset W^{1,\infty}_0(B_{2k}(0))$. 
     Clearly, $\mathbf{E}_\Lambda$ is convex and closed under the strong topology of $W^{1,1}_0(B_{2k}(0))$. Moreover, if $u\in \mathbf{E}_\Lambda$, then $|u|^{q_1}|\nabla u|^{q_2}\in L^1(B_{k}(0))$. \\                           
                   We consider the map $S:\mathbf{E}_\Lambda\mapsto W^{1,1}_0(\Omega)$ defined for each $v\in \mathbf{E}_\Lambda$ by $S(v)=u$, where $u\in W^{1,1}_0(\Omega)$  is the unique renormalized solution of   
                   \begin{equation}\label{pr-4}
\BA {lll}
-\text{div}(A(\nabla u))=\chi_{_{B_k(0)}}|v|^{q_1-1}v|\nabla v|^{q_2}+\mu_{n,k}\quad&\text{in }\;B_{2k}(0),\\ [1mm]
\phantom{-\text{div}(A(x,\nabla ))}
u=0\quad&\text{on }\;\partial B_{2k}(0).
\EA
 \end{equation}                                                                  
Notice that $\mathbf{W}_{1,p}[|\mu_{n,k}|]\leq C \left(\mathbf{I}_p[|\mu_{n,k}|]\right)^{\frac{1}{p-1}}$,
     this is due to the fact that, since $1<p<2$, then 
$$W_{1,p}[\mu](x)\sim \sum a_n^{\frac 1{p-1}} \leq \left( \sum a_n\right)^{\frac 1{p-1}} \sim (I_{p}[\mu])^{1/(p-1)}\text{ with }a_n=\myfrac{\mu(B(x,2^n))}{2^{n(N-p)}}$$
By Theorem \ref{renorm-th-2} and Corollary \ref{renorm-cor-1}, we have 
       \begin{align*}
      & |u|\leq C \left(\mathbf{I}_p[\chi_{_{B_k(0)}}|v|^{q_1}|\nabla v|^{q_2}dx+|\mu_{n,k}|]\right)^{\frac{1}{p-1}},\\&|\nabla u|\leq C \left(\mathbf{I}_1[\chi_{_{B_k(0)}}|v|^{q_1}|\nabla v|^{q_2}dx+|\mu_{n,k}|]\right)^{\frac{1}{p-1}},
       \end{align*}
       in $B_k(0)$. By the definition of $v$ we get 
        \begin{align*}
          & |u|\leq C \left(\Lambda^{q_1+q_2}\mathbf{I}_p[\left(\mathbf{I}_p[|\mu_{n,k}|]\right)^{\frac{q_1}{p-1}}\left(\mathbf{I}_1[|\mu_{n,k}|]\right)^{\frac{q_2}{p-1}}]+\mathbf{I}_p[|\mu_{n,k}|]\right)^{\frac{1}{p-1}},\\&|\nabla u|\leq C \left(\Lambda^{q_1+q_2}\mathbf{I}_1[\left(\mathbf{I}_p[|\mu_{n,k}|]\right)^{\frac{q_1}{p-1}}\left(\mathbf{I}_1[|\mu_{n,k}|]\right)^{\frac{q_2}{p-1}}]+\mathbf{I}_1[|\mu_{n,k}|]\right)^{\frac{1}{p-1}},
           \end{align*}
           in $B_k(0)$. 
          Using \eqref{pr-1}, we obtain
     \begin{align*}
            & |u|\leq C \left(\Lambda^{q_1+q_2}M\mathbf{I}_p[|\mu_{n,k}|]+\mathbf{I}_p[|\mu_{n,k}|]\right)^{\frac{1}{p-1}}=C \left(\Lambda^{q_1+q_2}M+1\right)^{\frac{1}{p-1}}\left(\mathbf{I}_p[|\mu_{n,k}|]\right)^{\frac{1}{p-1}},\\&|\nabla u|\leq C \left(\Lambda^{q_1+q_2}M\mathbf{I}_1[|\mu_{n,k}|]+\mathbf{I}_1[|\mu_{n,k}|]\right)^{\frac{1}{p-1}}=C \left(\Lambda^{q_1+q_2}M+1\right)^{\frac{1}{p-1}}\left(\mathbf{I}_1[|\mu_{n,k}|]\right)^{\frac{1}{p-1}},
             \end{align*}
             in $B_k(0)$.\
             We choose  $$\Lambda=C\left(\frac{q_1+q_2}{q_1+q_2-p+1}\right)^{\frac{1}{p-1}},~~M=\Lambda^{-q_1-q_2}\left(\left(\frac{q_1+q_2}{q_1+q_2-p+1}\right)^{\frac{1}{p-1}}-1\right),$$ then $C \left(\Lambda^{q_1+q_2}M+1\right)^{\frac{1}{p-1}}=\Lambda$ and $u\in \mathbf{E}_\Lambda$. Hence, $S(\mathbf{E}_\Lambda)\subset \mathbf{E}_\Lambda$.  \smallskip
             
             Next we show that $S:\mathbf{E}_\Lambda\mapsto \mathbf{E}_\Lambda$ is continuous.  Let $\{v_n\}$ be a sequence in $\mathbf{E}_\Lambda$ such that $v_m$ converges strongly in $W_0^{1,1}(B_{2k}(0))$ to a function $v\in\mathbf{E}_\Lambda$. Set $u_m=S(v_m)$. We need to show that $u_m\to S(v)$ in  $W_0^{1,1}(B_{2k}(0))$. We have 
\begin{equation*}
\BA {lll}
-\text{div}(A(\nabla u_m))=\chi_{_{B_k(0)}}|v_m|^{q_1-1}v_m|\nabla v_m|^{q_2}+\mu_{n,k}\quad&\text{in }\;B_{2k}(0),\\ [1mm]
\phantom{-\text{div}(A(\nabla))}
u_m=0\quad&\text{on }\;\partial B_{2k}(0),
\EA
 \end{equation*}  
and \begin{align*}
               |u_m|, |v_m|\leq \Lambda\left(\mathbf{I}_{p}[|\mu_{n,k}|]\right)^{\frac{1}{p-1}},~~~|\nabla u_m|, |\nabla v_m|\leq \Lambda\left(\mathbf{I}_{1}[|\mu_{n,k}|]\right)^{\frac{1}{p-1}} ~~\text{ in }  B_k(0).
               \end{align*}
Since   $\left(\mathbf{I}_p[|\mu_{n,k}|]\right)^{\frac{q_1}{p-1}}\left(\mathbf{I}_1[|\mu_{n,k}|]\right)^{\frac{q_2}{p-1}}\in L^{1}_{loc}(\mathbb{R}^N)$, we obtain

$$\chi_{_{B_k(0)}}|v_m|^{q_1-1}v_m|\nabla v_m|^{q_2}\to \chi_{_{B_k(0)}}|v|^{q_1-1}v|\nabla v|^{q_2}\quad\text{as }\;n\to\infty.$$ 
Applying Theorem \ref{stab}, we derive that $u_n\to S(v)$ in    $W_0^{1,1}(B_{2k}(0))$ as $n\to\infty$. Similarly, we can prove that $S(\mathbf{E}_\Lambda)$ is pre-compact  under the strong topology of $W^{1,1}_0(B_{2k}(0))$. \smallskip

                            Thus, by Schauder Fixed Point Theorem, $S$ has a fixed point on $\mathbf{E}_\Lambda$. 
                     This means that for any $k,n\in \mathbb{N}$, problem \eqref{pr-3} has a renormalized solution $u_{n,k}$ satisfying \eqref{esstep1thm1}. \\
                     By Lemma \ref{equi-integ-condi}, $\left\{\left(\mathbf{I}_p[|\mu_{n,k}|]\right)^{\frac{q_1}{p-1}}\left(\mathbf{I}_1[|\mu_{n,k}|]\right)^{\frac{q_2}{p-1}}\right\}_n$
                      is equi-integrable in  $B_{2k}(0)$. Thus, by a standard compactness argument, we get that $u_{n,k}$ converges to a renormalized solution   $u_k$ of  \begin{equation}\label{pr-3b}
                      \BA{lll}%
                      -\text{div}(A(\nabla u_{k}))=\chi_{_{B_k(0)}}|u_{k}|^{q_1-1}u_k|\nabla u_{k}|^{q_2}+\chi_{_{B_{k}(0)}}\mu\quad&\text{in }\;B_{2k}(0),\\[1mm] 
                      \phantom{   -\text{div}(A(\nabla))}
                      u_{k}=0&\text{on }\;\partial B_{2k}(0),                                                                                          
                      \EA
                      \end{equation}
                      which satisfies 
                      \begin{align}\label{esstep1thm1b}|u_{k}|\leq C_0\left(\mathbf{I}_{p}[|\mu|]\right)^{\frac{1}{p-1}},~~~|\nabla u_{k}|\leq C_0\left(\mathbf{I}_{1}[|\mu|]\right)^{\frac{1}{p-1}}~~\text{in}~~B_k(0).
                      \end{align}
                      Finally, thanks to Theorem \ref{renorm-th-1},  there exists a subsequence of $\{u_k\}_{k}$, still denoted by $\{u_k\}_{k}$  and $u\in W^{1,1}_{loc}(\mathbb{R}^N)$ such that $u_k$ converges to $u$ and  $\nabla u_k$ converges to $\nabla u$ almost everywhere.
                      Since \begin{align*}
                                 \chi_{_{B_k(0)}} |u_k|\leq C_0\left(\mathbf{I}_{p}[|\mu|]\right)^{\frac{1}{p-1}},~~~\chi_{_{B_k(0)}}|\nabla u_k|\leq C_0\left(\mathbf{I}_{1}[|\mu|]\right)^{\frac{1}{p-1}} \quad\text{ for all }\; k,
                                  \end{align*} 
                     and  $\left(\mathbf{I}_p[|\mu|]\right)^{\frac{q_1}{p-1}}\left(\mathbf{I}_1[|\mu|]\right)^{\frac{q_2}{p-1}}\in L^{1}_{loc}(\mathbb{R}^N)$, thus  $\chi_{_{B_k(0)}} |u_k|^{q_1-1}u_k|\nabla u_k|^{q_2}\to |u|^{q_1-1}u|\nabla u|^{q_2}$ in $L^1_{loc}(\mathbb{R}^N)$.
This implies that, $u$ is a  solution of problem \eqref{In-1} with $g(x,u,\nabla u)=|u|^{q_1}u|\nabla u|^{q_2}$ in the sense of distributions in $\BBR^N$ and it satisfies  \begin{align}\label{esstep1thm1c}|u|\leq C_0\left(\mathbf{I}_{p}[|\mu|]\right)^{\frac{1}{p-1}},~~~|\nabla u|\leq C_0\left(\mathbf{I}_{1}[|\mu|]\right)^{\frac{1}{p-1}}~~\text{in}~~B_k(0).
\end{align}

\nind{\it Step 2: Case $p> 2$}. In order to obtain the result, we will use 
            \begin{align*}\mathbf{W}_{\alpha,p}\left[\left(\mathbf{W}_{1,p}[|\mu_{n,k}|]\right)^{q_1}\left(\mathbf{W}_{\frac{1}{p},p}[|\mu_{n,k}|]\right)^{q_2}\right]\leq M \mathbf{W}_{\alpha,p}[|\mu_{n,k}|]<\infty~~\text{ almost everywhere},
                   \end{align*}
with $\alpha=1$ and $\alpha=1/p$, instead of \eqref{pr-1}; and 
\begin{equation*}
\mathbf{F}_\Lambda=\left\{u\in W^{1,1}_0(B_{2k}(0)):|u|\leq \Lambda \mathbf{W}_{1,p}[|\mu_{n,k}|],~~|\nabla u|\leq \Lambda \mathbf{W}_{\frac{1}{p},p}[|\mu_{n,k}|]~\text{ in }~B_{k}(0)\right\},
\end{equation*}  
 instead of $\mathbf{E}_\Lambda$.  We omit the details. The proof is complete.
        \qeda\medskip\\\\
\nind       {\it Proof of Theorem \ref{mainthem2}}.\smallskip

\nind{\it Step 1: Case $\frac{3N-2}{2N-1}<p\leq 2$}. Let $\mu_{n}\in C_c^\infty(\Omega)$ such that $\mu_{n}$ converges to $\mu$ in the sense of theorem \ref{stab} and $|\mu_n|\leq \varphi_n*(|\mu|)$, where $\{\varphi_n\}$ is a sequence of mollifiers in $\mathbb{R}^N$. Thanks to Proposition \ref{prop-pot-1}, 
\begin{align}\label{cap-appthm2}
|\mu_{n}|(K)\leq C'C\text{Cap}_{\mathbf{G}_{\frac{q_1p+q_2}{q_1+q_2}},\frac{q_1+q_2}{q_1+q_2-p+1}}(K)\quad\text{for all compact }K\subset\mathbb{R}^N
\end{align} 
We will prove that if $C$ in \eqref{cap-appthm2} is small enough, then for any $n\in \mathbb{N}$ the problem 
\begin{equation}\label{pr-4b}
\BA{lll}%
-\text{div}(A(x,\nabla u_{n}))=|u_{n}|^{q_1-1}u|\nabla u_{n}|^{q_2}+\mu_{n}\quad&\text{in }\;\Omega\\[1mm] 
\phantom{-\text{div}(A(x,\nabla ))}
u_{n}=0&\text{on }\;\Omega,                                                                                          
\EA
\end{equation}
has a renormalized solution  satisfying 
\begin{align}\label{esstep1thm2}|u_{n}|\leq C_0\left(\mathbf{I}_{p}^{4R}[|\mu_{n}|]\right)^{\frac{1}{p-1}},~~~|\nabla u_{n}|\leq C_0\left(\mathbf{I}_{1}^{4R}[|\mu_{n}|]\right)^{\frac{1}{p-1}}~~\text{in}~~\Omega.
\end{align}
By Theorem \eqref{th-pot-5}, we need to prove  that there exists $M>0$ such that, if for $\alpha=1$ and $p$, the following inequalities hold, 
\begin{align}\label{pr-2}\mathbf{I}_{\alpha}\left[\left(\mathbf{I}_{p}^{4R}[|\mu_{n,k}|]\right)^{\frac{q_1}{p-1}}\left(\mathbf{I}_{1}^{4R}[[|\mu_{n,k}|]\right)^{\frac{q_2}{p-1}}\right]\leq M \mathbf{I}_{\alpha}^{4R}[|\mu_{n,k}|]<\infty~~\text{ almost everywhere},
\end{align}
then problem \eqref{pr-4b}  has a renormalized solution  satisfying \eqref{esstep1thm2}.

We have to  prove that there exists $M>0$ such that if $\alpha=1$ and $p$, there holds
             \begin{align}\label{pr-6}\mathbf{I}^{4R}_{\alpha}\left[\left(\mathbf{I}^{4R}_{p}[|\mu_n|]\right)^{\frac{q_1}{p-1}}\left(\mathbf{I}_{1}^{4R}[|\mu_n|]\right)^{\frac{q_2}{p-1}}\right]\leq M \mathbf{I}^{4R}_{\alpha}[|\mu_n|]<\infty~~\text{ almost everywhere in }~B_{2R}(x_0),
                   \end{align}                           
                    For $n\in \mathbb{N}$ fixed, we set \begin{equation*}
                     \mathbf{E}_\Lambda=\left\{u\in W^{1,1}_0(\Omega):|u|\leq \Lambda\left(\mathbf{I}_{p}^{4R}[|\mu_n|]\right)^{\frac{1}{p-1}},~|\nabla u|\leq \Lambda\left(\mathbf{I}_{1}^{4R}[|\mu_n|]\right)^{\frac{1}{p-1}}~\text{ in }~ \Omega\right\}.
                     \end{equation*}
                     Clearly, $\mathbf{E}_\Lambda$ is closed under the strong topology of $W^{1,1}_0(\Omega)$, convex and $|u|^{q_1}|\nabla u|^{q_2}\in L^\infty(\Omega)$ for any $u\in \mathbf{E}_\Lambda $.         
                     We consider the map $S: \mathbf{E}_\Lambda\mapsto W^{1,1}_0(\Omega)$ defined for each $v\in  \mathbf{E}_\Lambda $ by $S(v)=u$, where $u\in W^{1,1}_0(\Omega)$  is the unique renormalized solution of 
                             \begin{equation*} \begin{array}{lll}
                                            -\text{div}(A(x,\nabla u)) =|v|^{q_1-1}v|\nabla v|^{q_2} +\mu_n ~~\text{in}~\Omega \\[1mm]     
                                            \phantom{-\text{div}(A(x,\nabla ))}                
                                                                           u = 0\quad ~\text{on }~\partial\Omega,\\ 
                                                                           \end{array} \end{equation*}
                     We will show that $S(\mathbf{E}_\Lambda)$ is subset of  $\mathbf{E}_\Lambda$ for some $\Lambda>0$ small enough.\\
                     For $v\in \mathbf{E}_\Lambda$ and $u=S(v)$, we have 
                     \begin{equation*}
                    |v|\leq \Lambda\left(\mathbf{I}_{p}^{4R}[|\mu_n|]\right)^{\frac{1}{p-1}},~|\nabla v|\leq \Lambda\left(\mathbf{I}_{1}^{4R}[|\mu|_n]\right)^{\frac{1}{p-1}}.
                     \end{equation*}
                     In particular,
                     \begin{align*}
                    ||v||_{L^\infty(\Omega_{d/2})}\leq C_1\Lambda  d^{-\frac{N-p}{p-1}}(|\mu_n|(\Omega))^{1/(p-1)},~~ |||\nabla v|||_{L^\infty(\Omega_{d/2})}\leq C_1\Lambda  d^{-\frac{N-1}{p-1}}(|\mu_n|(\Omega))^{1/(p-1)},
                     \end{align*} where $\Omega_{d/2}=\{x\in\Omega: d(x,\partial\Omega)\leq d/2\}$.\\
                     From \eqref{pr-6} with $\alpha =1$ and $p$ we derive
                     \begin{align*}
                     \mathbf{I}_{p}^{4R}[||v|^{q_1-1}v|\nabla v|^{q_2} +\mu_n|]&\leq \Lambda^{q_1+q_2}\mathbf{I}_{p}^{4R}[\left(\mathbf{I}_{p}^{4R}[|\mu_n|]\right)^{\frac{q_1}{p-1}}\left(\mathbf{I}_{1}^{4R}[|\mu_n|]\right)^{\frac{q_2}{p-1}}]+ \mathbf{I}_{p}^{4R}[|\mu_n|]
                     \\& \leq \left(\Lambda^{q_1+q_2}M+1\right)\mathbf{I}_{p}^{4R}[|\mu_n|],             \end{align*}
                     and
                      \begin{align*}
                                    \mathbf{I}_{1}^{4R}[||v|^{q_1-1}v|\nabla v|^{q_2} +\mu_n|]&\leq \Lambda^{q_1+q_2}\mathbf{I}_{1}^{4R}[\left(\mathbf{I}_{p}^{4R}[|\mu|]\right)^{\frac{q_1}{p-1}}\left(\mathbf{I}_{1}^{4R}[|\mu_n|]\right)^{\frac{q_2}{p-1}}]+ \mathbf{I}_{1}^{4R}[|\mu_n|]
                                    \\& \leq \left(\Lambda^{q_1+q_2}M+1\right)\mathbf{I}_{1}^{4R}[|\mu_n|].             \end{align*} 
                 By Theorem \ref{renorm-th-2} and $\mathbf{W}^{4R}_{1,p}[|\mu_n|]\lesssim \left(\mathbf{I}^{4R}_p[|\mu_n|]\right)^{\frac{1}{p-1}}$ in $\Omega$, we have 
                 \bel{pr-8}\BA {lll}
                 |u|\leq  C_2\left(\mathbf{I}_{p}^{4R}[||v|^{q_1-1}v|\nabla v|^{q_2} +\mu_n|]\right)^{\frac{1}{p-1}} \\
                 \phantom{|u|}\leq C_2\left(\Lambda^{q_1+q_2}M+1\right)^{\frac{1}{p-1}}\left(\mathbf{I}_{p}^{4R}[|\mu_n|]\right)^{\frac{1}{p-1}}~~\text{ in }~\Omega.
                 \EA\ee           
                     From Corollary \ref{renorm-cor-1}, we derive
                     \begin{align}\nonumber |\nabla u(x)|&\leq C_3\left(\frac{R}{d}\right)^{N}\left(\mathbf{I}_{1}^{4R}[||v|^{q_1-1}v|\nabla v|^{q_2} +\mu_n|](x)\right)^{\frac{1}{p-1}}
                     \\& \leq C_4\left(\Lambda^{q_1+q_2}M+1\right)^{\frac{1}{p-1}}\left(\mathbf{I}_{1}^{4R}[|\mu_n|](x)\right)^{\frac{1}{p-1}}\label{pr-8*}
                     \end{align} 
                                 for any $x\in \Omega$ verifying $d(x,\partial\Omega)>d/4$. By the standard regularity results for quasilinear equations, we deduce 
       \begin{align*}
       |||\nabla u|||_{L^\infty(\Omega_{d/4})}&\leq C_5\left(||u||_{L^\infty(\Omega_{d/2})}+|||v|^{q_1}|\nabla v|^{q_2}||_{L^\infty(\Omega_{d/2})}^{1/(p-1)}\right), 
       \end{align*}
         where $C_5=C_5(N,p,\Omega)$.\smallskip
         
\nind (a) Estimate of $|||v|^{q_1}|\nabla v|^{q_2}||_{L^\infty(\Omega_{d/2})}^{1/(p-1)}$. From \eqref{pr-6} , we have $|\mu_n|(\Omega)\leq C_6 M^{\frac{p-1}{q_1+q_2-p+1}}$. Thus
                    \begin{align*}
                    |||v|^{q_1}|\nabla v|^{q_2}||_{L^\infty(\Omega_{d/2})}^{1/(p-1)}&\leq ||v||_{L^\infty(\Omega_{d/2})}^{\frac{q_1}{p-1}} |||\nabla v|||_{L^\infty(\Omega_{d/2})}^{\frac{q_2}{p-1}}
                    \\& \leq 
                    \left(C_1\Lambda  d^{-\frac{N-p}{p-1}}(|\mu_n|(\Omega))^{1/(p-1)}\right)^{\frac{q_1}{p-1}}\left(C_1\Lambda  d^{-\frac{N-1}{p-1}}(|\mu_n|(\Omega))^{1/(p-1)}\right)^{\frac{q_2}{p-1}}
                    \\&\leq C_7\Lambda^{\frac{q_1}{p-1}+\frac{q_2}{p-1}}(|\mu_n|(\Omega))^{\frac{q_1+q_2}{(p-1)^2}}
                    \\&\leq C_8\Lambda^{\frac{q_1}{p-1}+\frac{q_2}{p-1}}M^{\frac{1}{p-1}} \inf_{x\in \Omega}\left(\mathbf{I}_{1}^{4R}[|\mu_n|](x)\right)^{\frac{1}{p-1}},
                    \end{align*}
                    where $C_8=C_8(N,p,\alpha,q_1,q_2,d,R)$.\smallskip
         
\nind (b) Estimate of $||u||_{L^\infty(\Omega_{d/2})}$.
                    By \eqref{pr-8} we have  
                   \begin{align*}
                   ||u||_{L^\infty(\Omega_{d/2})}&\leq C_2\left(\Lambda^{q_1+q_2}M+1\right)^{\frac{1}{p-1}}\left(||\mathbf{I}_{p}^{4R}[|\mu_n|]||_{L^\infty(\Omega_{d/2})}\right)^{\frac{1}{p-1}}
                   \\&\leq C_9\left(\Lambda^{q_1+q_2}M+1\right)^{\frac{1}{p-1}}d^{-\frac{N-p}{p-1}}\left(|\mu_n|(\Omega)\right)^{\frac{1}{p-1}}\\&
                   \leq  C_{10}\left(\Lambda^{q_1+q_2}M+1\right)^{\frac{1}{p-1}}\inf_{x\in \Omega}\left(\mathbf{I}_{1}^{4R}[|\mu_n|](x)\right)^{\frac{1}{p-1}}.
                   \end{align*}
                 Therefore, 
                 \begin{align*}
                  |||\nabla u|||_{L^\infty(\Omega_{d/4})}&\leq C_{11}\left(\Lambda^{\frac{q_1}{p-1}+\frac{q_2}{p-1}}M^{\frac{1}{p-1}}+\left(\Lambda^{q_1+q_2}M+1\right)^{\frac{1}{p-1}}\right)\inf_{x\in \Omega}\left(\mathbf{I}_{1}^{4R}[|\mu_n|](x)\right)^{\frac{1}{p-1}}, 
                  \end{align*}     
                      where $C_{11}=C_{11}(N,p,\alpha,q_1,q_2,d,R,\Omega)$.\\
                     Combining this with \eqref{pr-8*} we get for all $x\in\Omega$,
                     \begin{align}\nonumber
                     |\nabla u(x)|&\leq C_4 \left(\Lambda^{q_1+q_2}M+1\right)^{\frac{1}{p-1}}\left(\mathbf{I}_{1}^{4R}[|\mu_n|](x)\right)^{\frac{1}{p-1}}\\& +C_{11}\left(\Lambda^{\frac{q_1}{p-1}+\frac{q_2}{p-1}}M^{\frac{1}{p-1}}+\left(\Lambda^{q_1+q_2}M+1\right)^{\frac{1}{p-1}}\right)\left(\mathbf{I}_{1}^{4R}[|\mu_n|](x)\right)^{\frac{1}{p-1}}. \label{pr-8**}
                     \end{align} 
                   We can find $M,\Lambda>0$ such that 
                     \begin{align*}
                     & C_2\left(\Lambda^{q_1+q_2}M+1\right)^{\frac{1}{p-1}}\leq \Lambda,\\&
                    C_4 \left(\Lambda^{q_1+q_2}M+1\right)^{\frac{1}{p-1}}+C_{11}\left(\Lambda^{\frac{q_1}{p-1}+\frac{q_2}{p-1}}M^{\frac{1}{p-1}}+\left(\Lambda^{q_1+q_2}M+1\right)^{\frac{1}{p-1}}\right)\leq \Lambda. 
                     \end{align*}             
                     Thus, from \eqref{pr-8} and \eqref{pr-8**} we obtain  $S(\mathbf{E}_\Lambda)\subset \mathbf{E}_\Lambda$. Moreover,  it can be shown that the map $S:\mathbf{E}_\Lambda\mapsto \mathbf{E}_\Lambda$ is continuous and $S(\mathbf{E}_\Lambda)$ is pre-compact  under the strong topology of $W^{1,1}_0(\Omega)$. Then by Schauder Fixed Point Theorem, $S$ has a fixed point on $\mathbf{E}_\Lambda$. This means problem \eqref{pr-4b}  has a renormalized solution  satisfying \eqref{esstep1thm2}. \smallskip
                     
\nind{\it Step 2: The case $p>2$}.  To obtain the result, we will use 
\begin{align*}\mathbf{W}_{\alpha,p}^{4R}\left[\left(\mathbf{W}_{1,p}^{4R}[|\mu_n|]\right)^{q_1}\left(\mathbf{W}_{\frac{1}{p},p}[|\mu_n|]\right)^{q_2}\right]\leq M \mathbf{W}_{\alpha,p}^{4R}[|\mu|]<\infty~~\text{ almost everywhere in}~\Omega
\end{align*}
with $\alpha=1$ and $\alpha=1/p$, instead of \eqref{pr-1}; and 
\begin{equation*}
\mathbf{F}_\Lambda=\left\{u\in W^{1,1}_0(\Omega):|u|\leq \Lambda \mathbf{W}_{1,p}^{4R}[|\mu_n|],~~|\nabla u|\leq \Lambda \mathbf{W}_{\frac{1}{p},p}^{4R}[|\mu_n|]~\text{ in }~\Omega\right\},
\end{equation*}  
instead of $\mathbf{E}_\Lambda$.  We omit the details. The proof is complete.
 \qeda\medskip

\end{document}